\documentclass{amsart}
\usepackage{amssymb,amsmath,amsthm,graphicx,booktabs,array,url}
\newtheorem{theorem}{Theorem}[section]
\newtheorem{lemma}[theorem]{Lemma}
\newtheorem{proposition}[theorem]{Proposition}
\newtheorem{corollary}[theorem]{Corollary}
\theoremstyle{definition}\newtheorem{definition}[theorem]{Definition}
\theoremstyle{remark}\newtheorem{remark}[theorem]{Remark}
\numberwithin{equation}{section}
\newcommand{\Om}{\Omega}\newcommand{\R}{\mathbb{R}}\newcommand{\N}{\mathbb{N}}
\DeclareMathOperator{\dist}{dist}\DeclareMathOperator{\supp}{supp}\DeclareMathOperator{\tr}{tr}\DeclareMathOperator{\Div}{div}\DeclareMathOperator{\curl}{curl}\DeclareMathOperator{\sym}{sym}

\newcommand{\recR}[1]{record~R#1 of \cite{lbb3drepo}}
\newcommand{\docD}[1]{document~D#1 of \cite{lbb3drepo}}
\newcommand{\M}{\mathcal M}\newcommand{\Mz}{\mathcal M_0}
\newcommand{\kess}{\kappa_{\rm e}}

\newcommand{\FichKappaZero}{\frac{169}{200}}

\newcommand{\FichBetaLow}{0.278388}
\newcommand{\FichBetaLowSq}{0.077499878544}%
\newcommand{\FichBetaSqLow}{\frac{31}{400}}
\newcommand{\FichBetaSqLowDec}{0.0775}
\newcommand{\FichWidth}{1.15\cdot10^{-1}}%
\newcommand{\FichGap}{0.042546100960}%

\newcommand{\FichLamQ}{0.967969333}%
\newcommand{\FichLamS}{0.440939565}%
\newcommand{\FichMarginQ}{0.032}

\newcommand{\FichBindClass}{{\rm Q0}}
\newcommand{\FichNQ}{492}
\newcommand{\FichNS}{202}
\newcommand{\FichNQone}{532}%
\newcommand{\FichMergedQ}{40}
\newcommand{\FichMergedS}{43}
\newcommand{\FichDofQ}{1\,251\,541}
\newcommand{\FichDofS}{425\,221}
\newcommand{\FichPointsQ}{4.24\cdot10^9}%
\newcommand{\FichPointsS}{1.41\cdot10^9}%

\newcommand{\FichRho}{\frac1{128}}
\newcommand{\FichHmax}{\frac1{16}}

\newcommand{\FichTetsTwo}{223\,608}
\newcommand{\FichVertsTwo}{41\,009}
\newcommand{\FichTetsSix}{74\,536}
\newcommand{\FichVertsSix}{14\,191}
\newcommand{\FichBalls}{64}
\newcommand{\FichOrbits}{16}
\newcommand{\FichBallsTwo}{37}
\newcommand{\FichTaylor}{3}
\newcommand{\FichRadii}{\{\frac12,\frac14\}}

\newcommand{\FichFacesSix}{152\,460}
\newcommand{\FichBdSix}{6\,776}

\newcommand{\FichFacesTwo}{454\,476}
\newcommand{\FichBdTwo}{14\,520}

\newcommand{\FichPiecesQ}{8\,118\,432}
\newcommand{\FichDistinctQ}{4\,304\,841}
\newcommand{\FichPiecesS}{2\,706\,144}
\newcommand{\FichDistinctS}{1\,434\,947}
\newcommand{\FichRules}{\{1,2,3,4,5,7,8\}}
\newcommand{\FichRuleMax}{8}

\newcommand{\FichKpartial}{0.81991}%

\title[Verified inf-sup bounds on rectangular polyhedra]{Verified two-sided bounds for the inf-sup constant of the divergence\\ on rectangular polyhedra}
\author{Xuefeng Liu}
\address{Department of Information and Sciences, Tokyo Woman's Christian University, 2-6-1 Zempukuji, Suginami-ku, Tokyo 167-8585, Japan}
\email{xfliu@lab.twcu.ac.jp}

\begin{document}
\begin{abstract}
We prove two-sided bounds for the inf-sup constant $\beta(\Om)$ of the divergence, the constant in the Ladyzhenskaya--Babu\v{s}ka--Brezzi condition of the
Stokes problem. The domains are the cube $(-1,1)^3$, the box $(0,2)\times(0,1)^2$, the L-shaped prism $((-1,1)^2\setminus[0,1)^2)\times(0,1)$ and the Fichera
domain $(-1,1)^3\setminus[0,1]^3$. The bounds are $0.340587\le\beta\le0.3411711$ for the cube, $0.254950\le\beta\le0.267469$ for the box,
$0.193649\le\beta\le0.208165$ for the L-shaped prism and $\FichBetaLow\le\beta\le0.314282$ for the Fichera domain. The lower bounds use the Friedrichs--Velte
formulation $\beta^2=(1-\kappa)/2$ over pairs $(h,g)$ with $\nabla h=\curl g$, $\Div g=0$. A multiplier identity, from a divergence-free quadratic tensor of such
pairs, controls them near the boundary through a matrix-valued weight; Taylor functionals of $h$ on balls control the interior. The criterion is one matrix
inequality per symmetry class, verified in exact and ball arithmetic. The upper bounds come from polynomial, edge-singular and vertex-singular trial pairs with
rigorously enclosed Rayleigh quotients. Every convex rectangular vertex generates essential spectrum of the Cosserat operator at or below $0.116398$, strictly
below the interval $[\frac12-\frac1\pi,\frac12+\frac1\pi]$ of right-angle edges. For the cube this gives $0.116\le\beta^2\le\inf\sigma_{\rm ess}\le0.116398$; for
the other three domains $\beta^2$ lies strictly below this vertex value. Open are whether $\beta^2$ of the cube equals $\inf\sigma_{\rm ess}$ (an isolated
eigenvalue below it is not excluded), and whether the verified conical multiplier bound and the enclosed angular Rayleigh quotient at the vertex,
which agree to $2.1\cdot10^{-6}$, are equal.
\end{abstract}
\maketitle

\section{Introduction}
Let $\Om\subset\R^3$ be a bounded domain and $L^2_0(\Om):=\{q\in L^2(\Om):\int_\Om q=0\}$. The inf-sup constant
\begin{equation}\label{A:beta}
\beta(\Om):=\inf_{q\in L^2_0(\Om)\setminus\{0\}}\ \sup_{v\in H^1_0(\Om)^3\setminus\{0\}}\frac{(q,\Div v)}{|v|_1\,\|q\|}
\end{equation}
enters the stability of the Stokes and Navier--Stokes equations. An explicit lower bound for it is needed in computer-assisted existence proofs in which the
pressure is estimated \cite{LiuNakaoOishi}. Analytic lower bounds of Horgan--Payne type are known for star-shaped domains \cite{HorganPayne,CostabelDauge},
\cite{Payne2007}, \cite[Sect.~3.2]{Zsuppan2020}, but they are far from the true value. Numerical methods give upper bounds or rest on uncheckable
hypotheses \cite{Gallistl19,Gallistl21}.

On polyhedra the Cosserat operator $\mathcal S$ (Definition~\ref{A:def:S}), the bottom of whose spectrum is $\beta^2$ (Lemma~\ref{A:lem:S}), has essential spectrum: every edge of
opening $\omega$ contributes the interval $[\frac12-\frac{|\sin\omega|}{2\omega},\frac12+\frac{|\sin\omega|}{2\omega}]$ to $\sigma_{\rm ess}(\mathcal S)$ \cite[Thm.~4.1]{CCDL} (see
Theorem~\ref{B-thm:edge}). For the cube this gives $\beta^2\le\frac12-\frac1\pi$. Whether the bottom of $\sigma_{\rm ess}$ on cuboids comes from edges or corners was left open in
\cite[Sect.~5.2]{CCDL}, where a corner contribution is conjectured from the computations (see Section~\ref{sec:ess} and Remark~\ref{B-rem:fe}).

In a companion paper \cite{companion} we proposed, for planar polygons, a certificate for lower bounds of $\beta$ that is unconditional and convergent. It is based on the
Bergman-space formulation of Friedrichs and a multiplier identity for holomorphic functions. Here we extend the method to three dimensions and apply it to rectangular polyhedra,
by which we mean bounded polyhedral domains, convex or not, whose faces are parallel to the coordinate planes. For the cube, the box, the L-shaped prism and the Fichera domain
we obtain the following bounds.
\begin{theorem}[main results]\label{thm:main}
\begin{center}
\footnotesize
\begin{tabular}{@{}lll@{}}
\toprule
domain $\Om$ & certified bounds & $\inf\sigma_{\rm ess}(\mathcal S_\Om)\in$\\
\midrule
$(-1,1)^3$ & $0.340587\le\beta\le0.3411711$ & $[0.116,\,0.116397718977]$\\
$(0,2)\times(0,1)^2$ & $0.254950\le\beta\le0.267469$ & $[0.065,\,0.116397718977]$\\
$((-1,1)^2\setminus[0,1)^2)\times(0,1)$ & $0.193649\le\beta\le0.208165$ & $[0.0375,\,0.116397718977]$\\
$(-1,1)^3\setminus[0,1]^3$ & $\FichBetaLow\le\beta\le0.314282$ & $[\FichBetaSqLowDec,\,0.116397718977]$\\
\bottomrule
\end{tabular}
\end{center}
Upper bounds: Theorem~\ref{thm:upper} and Propositions~\ref{oth:prop:boxup}, \ref{oth:prop:lpup}, \ref{oth:prop:fichup}. Lower bounds: Theorems~\ref{thm:lower},
\ref{oth:thm:box}, \ref{oth:thm:lp}, \ref{oth:thm:fich}. Essential spectrum: Corollaries~\ref{cor:ess} and~\ref{B-cor:cubeess}.
For the last three domains $\beta^2\le0.0988<\frac12(1-L)$, where $\frac12(1-L)\ge0.116396$ is the value generated by a convex rectangular vertex
(Corollary~\ref{B-cor:cubeess}(d)).
\end{theorem}
The following three arguments prove these bounds.
\begin{enumerate}
\item \emph{Lower bounds for $\beta$ (Sections~\ref{sec:fv}--\ref{sec:cert}, \ref{sec:lower}, \ref{sec:other}).} By Theorem~\ref{thm:fv}, $\beta^2=(1-\kappa)/2$, where
$\kappa$ is the supremum of a quotient over pairs $\nabla h=\curl g$, $\Div g=0$. A conservation law $\Div T=0$
(Lemma~\ref{lem:T}; an identity of Payne \cite[(2.9)]{Payne2007}) gives $\int(h^2-|g|^2)=\int F^{\mathsf T}W(m)F$ for every Lipschitz $m$ vanishing on $\partial\Om$,
without regularity of $\Om$ or traces of $F$ (Lemma~\ref{lem:mult}). A matrix weight near the boundary and Taylor functionals of $h$ on interior balls reduce
$\kappa\le\kappa_0$ to one test $\lambda_{\max}(G)\le1$ per symmetry class (Theorem~\ref{thm:rel}, Proposition~\ref{oth:prop:classes}), verified in exact and ball
arithmetic (Section~\ref{sec:verif}).
\item \emph{Upper bounds for $\beta$ (Sections~\ref{sec:upper}, \ref{sec:other}).} Every pair gives a lower bound for $\kappa$; we use explicit polynomial, edge- and
vertex-singular pairs.
\item \emph{Essential spectrum (Section~\ref{sec:ess}).} Vertex-singular pairs whose homogeneity degree tends to $-\frac32$ concentrate at a convex rectangular vertex,
and the normalised pressures converge weakly to zero. Hence such a vertex generates essential spectrum of $\mathcal S$ at or below $0.116397718977<\frac12-\frac1\pi$
(Theorem~\ref{thm:ess}, Corollary~\ref{cor:ess}). A single trial bounds only $\min\sigma(\mathcal S)$; this sequence bounds $\inf\sigma_{\rm ess}(\mathcal S)$.
\end{enumerate}
Section~\ref{sec:vertex} encloses the angular trial constant and the pointwise multiplier constant at the vertex, $0.767204562046\le L\le U_*\le0.767206641725$
(Theorem~\ref{thm:vexact}, Proposition~\ref{prop:UL}). Appendix~\ref{A:app:lit} contains auxiliary proofs, Appendix~\ref{app:design} design data, and the repository
\cite{lbb3drepo} literature proofs, further tables and a reading guide.

\emph{Notation.} Throughout, $\Om\subset\R^3$ is a bounded domain (open and connected). $\|\cdot\|$ and $(\cdot,\cdot)$ are the norm and inner product of $L^2(\Om)$ (for scalar, vector or matrix fields), $|v|_1:=\|\nabla v\|$, and
summation over repeated indices is used. $\varepsilon_{ijk}$ is the Levi-Civita symbol, so that $(a\times b)_i=\varepsilon_{ijk}a_jb_k$ and
$(\curl v)_i=\varepsilon_{ijk}\partial_jv_k$. For a vector field $m$, $(Dm)_{ij}:=\partial_jm_i$. For symmetric matrices $A\preceq B$ means that $B-A$ is positive
semidefinite; $e_0:=(1,0,0,0)^{\mathsf T}$; vectors in $\R^4$ are ordered as $F=(h,g_1,g_2,g_3)$. The main symbols, with the place of their definition, are listed in
\docD{1}.

\section{The Friedrichs--Velte formulation}\label{sec:fv}
This section reduces bounds for $\beta$ to bounds for a quotient over pairs of harmonic fields (Theorem~\ref{thm:fv}).
\subsection{Pairs}
The space of pairs and its subspace of pairs with mean-free first component are
\begin{equation}\label{A:M}
\begin{aligned}
\M&:=\big\{F=(h,g)\in L^2(\Om)\times L^2(\Om)^3:\ \nabla h=\curl g,\ \Div g=0\ \text{in }\mathcal D'(\Om)\big\},\\
\Mz&:=\Big\{F\in\M:\ \int_\Om h=0\Big\}.
\end{aligned}
\end{equation}
The system $\nabla h=\curl g$, $\Div g=0$ is the three-dimensional counterpart of the Cauchy--Riemann
equations; it is the Moisil--Teodorescu system of quaternionic analysis \cite[Rem.~1.2]{CostabelForms}. We put
\begin{equation}\label{kappa}
|F|^2:=h^2+|g|^2,\qquad \kappa(\Om):=\sup_{F\in\Mz\setminus\{0\}}\frac{\|h\|^2-\|g\|^2}{\|h\|^2+\|g\|^2}.
\end{equation}
\begin{lemma}[regularity of pairs]\label{A:lem:reg}
Every component of $F=(h,g)\in\M$ is harmonic, hence real-analytic in $\Om$.
\end{lemma}
\begin{proof}
In $\mathcal D'(\Om)$, $\Delta h=\Div\curl g=0$ and $\Delta g=\nabla\Div g-\curl\curl g=-\curl\nabla h=0$; apply Weyl's lemma.
\end{proof}
\begin{lemma}[div--curl identity {\cite[(1.1)]{CostabelForms}}]\label{A:lem:divcurl}
For every $v\in H^1_0(\Om)^3$, $|v|_1^2=\|\Div v\|^2+\|\curl v\|^2$. No regularity of $\partial\Om$ is needed.
\end{lemma}
\begin{lemma}[duality]\label{A:lem:dual}
For every $F=(h,g)\in\M$,
\begin{equation}\label{dual0}(h,\Div v)=-(g,\curl v)\qquad\text{for all }v\in H^1_0(\Om)^3 .\end{equation}
\end{lemma}
\begin{proof}
For $v\in C_c^\infty(\Om)^3$, $(h,\Div v)=-\langle\nabla h,v\rangle=-\langle\curl g,v\rangle=-(g,\curl v)$, since $\curl$ is formally self-adjoint. Both sides are continuous in $v$ for the
$H^1$ norm, and $C_c^\infty(\Om)^3$ is dense in $H^1_0(\Om)^3$.
\end{proof}

\subsection{The Cosserat operator}
By Lemma~\ref{A:lem:divcurl}, $|(q,\Div v)|\le\|q\|\,\|\Div v\|\le\|q\|\,|v|_1$. Hence the Riesz representation theorem in $H^1_0(\Om)^3$, with the inner
product $(\nabla\cdot,\nabla\cdot)$, gives the following definition.
\begin{definition}[Cosserat operator]\label{A:def:S}
For $q\in L^2_0(\Om)$ let $v_q\in H^1_0(\Om)^3$ be the solution of $(\nabla v_q,\nabla v)=(q,\Div v)$ for all $v\in H^1_0(\Om)^3$, and
\begin{equation}\label{A:S}
\mathcal Sq=\mathcal S_\Om q:=\Div v_q .
\end{equation}
\end{definition}
\begin{lemma}[{\cite[(1.5)]{CCDL}}]\label{A:lem:S}\label{B-lem:S}
$\mathcal S$ maps $L^2_0(\Om)$ into itself, is bounded and self-adjoint, $(\mathcal Sq,q)=|v_q|_1^2$, $0\le\mathcal S\le I$, and
\[\beta(\Om)^2=\inf_{q\in L^2_0(\Om),\ \|q\|=1}(\mathcal Sq,q)=\min\sigma(\mathcal S).\]
\end{lemma}
Proofs of Lemmas~\ref{A:lem:divcurl} and~\ref{A:lem:S} are reproduced in \docD{2}.

\subsection{The Friedrichs--Velte formulation}
\begin{theorem}[{\cite[Thm.~2.1]{CostabelForms}, \cite[Thm.~1]{Velte1998}}]\label{thm:fv}
$\beta(\Om)^2=(1-\kappa(\Om))/2$.
\end{theorem}
The identity is due to Velte \cite[Thm.~1]{Velte1998} for simply connected domains with $C^2$ boundary and to Costabel \cite[Thm.~2.1]{CostabelForms} for every
bounded open set; we include a short proof. Step~4 needs no eigenfunction, since $\beta^2$ may lie in $\sigma_{\rm ess}(\mathcal S)$.
\begin{proof}
\emph{Step 1: one pair bounds the Rayleigh quotient.} Let $F=(h,g)\in\M$, $q:=h-|\Om|^{-1}\int_\Om h$, $v:=v_q$ and $s:=(\mathcal Sq,q)=|v|_1^2=(q,\Div v)$. As
$\int_\Om\Div v=0$, \eqref{dual0} gives $s=(h,\Div v)=-(g,\curl v)$, so $s^2\le\|q\|^2\|\Div v\|^2$ and $s^2\le\|g\|^2\|\curl v\|^2$. If $s>0$, then $q,g\ne0$;
dividing, adding and using Lemma~\ref{A:lem:divcurl} gives $s^2(\|q\|^{-2}+\|g\|^{-2})\le|v|_1^2=s$. Hence, also if $s=0$,
\begin{equation}\label{fv:pair}
(\mathcal Sq,q)\le\frac{\|q\|^2\,\|g\|^2}{\|q\|^2+\|g\|^2}\qquad(q\ne0\text{ or }g\ne0).
\end{equation}

\emph{Step 2: $\kappa\le1-2\beta^2$.} Let $F\in\Mz\setminus\{0\}$, so $q=h$. If $h\ne0$, then $\beta^2\|h\|^2\le(\mathcal Sh,h)$ (Lemma~\ref{A:lem:S}) and \eqref{fv:pair}
give $\beta^2\le\|g\|^2/(\|h\|^2+\|g\|^2)$, i.e.\ the quotient \eqref{kappa} of $F$ is at most $1-2\beta^2$. If $h=0$, it is $-1\le1-2\beta^2$, as $\beta^2\le1$.

\emph{Step 3: $b:=\beta^2<1$.} The pair $(x_1-\overline{x_1},(0,0,x_2))$, with $\overline{x_1}$ the mean over $\Om$, lies in $\Mz$ because $\curl(0,0,x_2)=e_1$, and
$h\ne0$; Step~2 gives $b\le\|g\|^2/(\|h\|^2+\|g\|^2)<1$.

\emph{Step 4: $\kappa\ge1-2b$.} As $b\in\sigma(\mathcal S)$, the spectral theorem gives $q_n\in L^2_0(\Om)$, $\|q_n\|=1$, $\|(\mathcal S-bI)q_n\|\to0$. Put
$h_n:=\mathcal Sq_n-q_n\in L^2_0(\Om)$ and $g_n:=\curl v_{q_n}$. By Definition~\ref{A:def:S}, $\Delta v_{q_n}=\nabla q_n$ in $\mathcal D'(\Om)$, so
$\curl g_n=\nabla\Div v_{q_n}-\Delta v_{q_n}=\nabla h_n$ and $\Div g_n=0$: $(h_n,g_n)\in\Mz$. Since $(\mathcal Sq_n,q_n)\to b$ and $\|\mathcal Sq_n\|^2\to b^2$,
Lemmas~\ref{A:lem:S} and~\ref{A:lem:divcurl} give
\begin{align*}
\|h_n\|^2&=\|\mathcal Sq_n\|^2-2(\mathcal Sq_n,q_n)+1\to(1-b)^2,\\
\|g_n\|^2&=|v_{q_n}|_1^2-\|\Div v_{q_n}\|^2=(\mathcal Sq_n,q_n)-\|\mathcal Sq_n\|^2\to b(1-b).
\end{align*}
The limits have sum $1-b>0$ (Step~3) and difference $(1-b)(1-2b)$, so the quotients of $(h_n,g_n)$ tend to $1-2b$.
\end{proof}

\subsection{Extruded pairs}
Extruded pairs reduce the Moisil--Teodorescu system to the Cauchy--Riemann equations \cite[Sect.~2]{Zsuppan2016} and carry planar bounds for $\kappa$ to prisms. For a bounded planar domain $\omega\subset\R^2$ let
\begin{equation}\label{A:kappa2}
\kappa(\omega):=\sup\Big\{\frac{|\int_\omega f^2|}{\|f\|^2_{L^2(\omega)}}:\ f\in L^2(\omega)\text{ holomorphic},\ \int_\omega f=0,\ f\ne0\Big\}.
\end{equation}
\begin{lemma}[extrusion]\label{A:lem:extr}
Let $\omega\subset\R^2$ be a bounded domain, $I\subset\R$ a bounded open interval, and $f=u+iv\in L^2(\omega)$ holomorphic with $\int_\omega f=0$, $f\ne0$. Then
\begin{equation}\label{A:extr}
F_f:=\big(u(x_1,x_2),\,(0,0,v(x_1,x_2))\big)\in\Mz(\omega\times I)\setminus\{0\},
\end{equation}
and
\[\max_{\phi\in\R}\frac{\|\Re(e^{i\phi}f)\|^2-\|\Im(e^{i\phi}f)\|^2}{\|\Re(e^{i\phi}f)\|^2+\|\Im(e^{i\phi}f)\|^2}=\frac{|\int_\omega f^2|}{\|f\|^2_{L^2(\omega)}},\]
where the left-hand side is the quotient of \eqref{kappa} for $F_{e^{i\phi}f}$ on $\omega\times I$.
\end{lemma}
\begin{proof}
By the Cauchy--Riemann equations, $\curl(0,0,v)=\nabla u$ and $\Div(0,0,v)=0$. Moreover $u^2-v^2=\Re f^2$ and $u^2+v^2=|f|^2$, and the factor $|I|$ cancels. The
details are in \docD{2}.
\end{proof}
\begin{corollary}\label{A:cor:extr}
$\kappa(\omega\times I)\ge\kappa(\omega)$. In particular $\kappa((-1,1)^3)\ge\kappa((-1,1)^2)\ge2/\pi$.
\end{corollary}
\begin{proof}
The first claim follows from Lemma~\ref{A:lem:extr}. The bound $\kappa((-1,1)^2)\ge2/\pi$ is equivalent to
$\beta^2\le\frac12-\frac1\pi$ for the square \cite[Thm.~3.3]{CCDL}, by the planar Friedrichs--Velte identity \cite[Thm.~2.1]{CostabelDauge} (see also \cite{companion}).
\end{proof}

\section{A multiplier identity}\label{sec:mult}
The integrand $h^2-|g|^2$ in \eqref{kappa} has no sign. Lemma~\ref{lem:mult} replaces its integral by $\int F^{\mathsf T}W(m)F$, with a multiplier $m$ that is free in $\Om$
and vanishes on $\partial\Om$. The proof uses only a conservation law; it needs no regularity of $\partial\Om$ and no traces of $F$.
For $F=(h,g)\in\M$ define the tensor
\begin{equation}\label{T}T_{ij}:=g_ig_j+\tfrac12\big(h^2-|g|^2\big)\delta_{ij}-h\,\varepsilon_{ijk}g_k .\end{equation}
\begin{lemma}[{cf.\ \cite[(2.9)]{Payne2007}}]\label{lem:T} For smooth $(h,g)$, with $E:=\nabla h-\curl g$,
\[\partial_jT_{ij}=g_i\Div g-(E\times g)_i+h\,E_i .\]
In particular $\Div T=0$ for $F\in\M$.\end{lemma}
\begin{proof}
Differentiating the three terms of $T_{ij}$ gives $\partial_j(g_ig_j)=g_j\partial_jg_i+g_i\Div g$, $\partial_i\big(\tfrac12(h^2-|g|^2)\big)=h\partial_ih-g_j\partial_ig_j$
and $-\partial_j(h\varepsilon_{ijk}g_k)=-(\nabla h\times g)_i-h(\curl g)_i$. Since $g_j\partial_jg_i-g_j\partial_ig_j=(\curl g\times g)_i$, the sum is the claim.
For $F\in\M$ the components are smooth (Lemma~\ref{A:lem:reg}), $E=0$ and $\Div g=0$ (details in \docD{2}).
\end{proof}
The law $\Div T=0$ is the identity $\nabla(h^2-|g|^2)=2\curl(hg)-2\Div(g\otimes g)$ used by Payne \cite[(2.9)]{Payne2007} for star-shaped domains (see also
\cite[Sect.~3.2, (33), (38)]{Zsuppan2020}). For a real $3\times3$ matrix $A$ let $\sym A:=\frac12(A+A^{\mathsf T})$ and $\operatorname{ax}(A)_k:=\varepsilon_{kij}A_{ji}$, so that $\operatorname{ax}(Dm)=\curl m$. We define the symmetric
matrix, with $F$ ordered as $(h,g_1,g_2,g_3)$,
\begin{equation}\label{W}
W(A):=\begin{pmatrix}1-\frac12\tr A&-\frac12\operatorname{ax}(A)^{\mathsf T}\\[2pt]-\frac12\operatorname{ax}(A)&(-1+\frac12\tr A)I_3-\sym A\end{pmatrix},\qquad W(m):=W(Dm),
\end{equation}
so that $W(m)$ has the entries $1-\frac12\Div m$, $-\frac12\curl m$ and $(-1+\frac12\Div m)I_3-\sym Dm$.
\begin{lemma}\label{lem:mult}
Let $m$ be Lipschitz on $\overline\Om$ with $m=0$ on $\partial\Om$. Then $\displaystyle\int_\Om\big(h^2-|g|^2\big)=\int_\Om F^{\mathsf T}W(m)F$ for all $F\in\M$.
\end{lemma}
\begin{proof}
\emph{Step 1: the contraction.} With $T:Dm:=T_{ij}\partial_jm_i$,
\begin{align*}T:Dm&=g_ig_j\partial_jm_i+\tfrac12(h^2-|g|^2)\Div m-hg_k\,\varepsilon_{ijk}\partial_jm_i\\
&=g^{\mathsf T}Dm\,g+\tfrac12(h^2-|g|^2)\Div m+h\,g\cdot\curl m,\end{align*}
because $\varepsilon_{ijk}\partial_jm_i=-\varepsilon_{kji}\partial_jm_i=-(\curl m)_k$. Since $g^{\mathsf T}Dm\,g=g^{\mathsf T}\sym Dm\,g$,
\begin{equation}\label{hgT}\begin{aligned}h^2-|g|^2-T:Dm&=h^2(1-\tfrac12\Div m)+|g|^2(-1+\tfrac12\Div m)\\
&\quad-g^{\mathsf T}\sym Dm\,g-h\,g\cdot\curl m=F^{\mathsf T}W(m)F .\end{aligned}\end{equation}

\emph{Step 2: compactly supported fields.} Let $\varphi$ be Lipschitz with compact support in $\Om$. As $T$ is smooth and $\Div T=0$ (Lemma~\ref{lem:T}), integration by parts on a
smooth open set $\Om'\Subset\Om$ containing $\supp\varphi$ gives $\int_\Om T:D\varphi=-\int_{\Om'}(\Div T)\cdot\varphi=0$.

\emph{Step 3: truncation.} For $t>0$ let $m_{t,i}:=\operatorname{sign}(m_i)(|m_i|-t)_+$, $i=1,2,3$. As $m$ is continuous on $\overline\Om$ and vanishes on $\partial\Om$, the set $\{|m_i|\ge t\}$ is a
compact subset of $\Om$; hence $m_t$ is Lipschitz with compact support, and $\int T:Dm_t=0$ by Step 2. By the chain rule for Lipschitz functions,
$Dm_{t,i}=Dm_i\mathbf 1_{\{|m_i|>t\}}$ a.e., which converges a.e.\ to $Dm_i\mathbf 1_{\{m_i\ne0\}}=Dm_i$ as $t\to0$; the last equality holds because
$Dm_i=0$ a.e.\ on $\{m_i=0\}$ \cite[Lemmas~7.6 and~7.7]{GilbargTrudinger}. Since $|T:Dm_t|\le|T|\,\|Dm\|_{L^\infty}$ and $|T|\le c|F|^2\in L^1$, dominated convergence gives $\int T:Dm=0$.

\emph{Step 4.} Integrating \eqref{hgT} and using $\int T:Dm=0$ gives the claim.
\end{proof}
If $m$ is continuous and piecewise linear on a tetrahedral mesh, $W(m)$ is constant on each tetrahedron $\tau$:
\begin{equation}\label{A:Wtau}
W_\tau:=W\big(Dm|_\tau\big).
\end{equation}
Near faces and edges a multiplier achieves the following.
\begin{proposition}\label{prop:local}
(a) Near a planar face with inward unit normal $\mathbf n$, the field $m:=2\dist(\cdot,\text{face})\,\mathbf n$ gives a matrix $W(m)$ with eigenvalues
$0,0,0,-2$; hence $F^{\mathsf T}W(m)F\le0$.
(b) Near an edge along the $x_3$-axis let $m=(m_1,m_2,0)$ with $m_1,m_2$ depending only on $(x_1,x_2)$, and let $w:=1-\partial_{\bar z}(m_1-im_2)$ be the planar weight of \cite{companion}. Then $W(m)$ splits into the
$(h,g_3)$-block, with eigenvalues $\pm|w|$, and the $(g_1,g_2)$-block, with eigenvalues $-1\pm\frac12(\lambda_1-\lambda_2)$, where $\lambda_1\ge\lambda_2$ are the eigenvalues of
$\sym D(m_1,m_2)$. In particular $\lambda_{\max}(W(m))\ge|w|$.
\end{proposition}
\begin{proof}(a) After a rigid motion the face is $\{x_3=0\}$ with $\mathbf n=e_3$, and $m=(0,0,2x_3)$. Then $Dm=2e_3e_3^{\mathsf T}$, $\Div m=2$, $\curl m=0$, so $1-\frac12\Div m=0$ and the
$g$-block is $0\cdot I-2e_3e_3^{\mathsf T}$: $W(m)=\operatorname{diag}(0,0,0,-2)$, and $F^{\mathsf T}W(m)F=-2g_3^2\le0$.

(b) Let $A:=D(m_1,m_2)$ ($2\times2$), $t:=\operatorname{tr}A$, $c:=\partial_1m_2-\partial_2m_1$. Since $m$ does not depend on $x_3$ and $m_3=0$, $Dm=\begin{pmatrix}A&0\\0&0\end{pmatrix}$, $\Div m=t$ and
$\curl m=c\,e_3$. Hence $W(m)$ has no entries coupling $(h,g_3)$ with $(g_1,g_2)$. The $(h,g_3)$-block is
$\begin{pmatrix}1-\frac t2&-\frac c2\\-\frac c2&-1+\frac t2\end{pmatrix}$; it has trace $0$ and determinant $-(1-\frac t2)^2-\frac{c^2}4$, so its eigenvalues are
$\pm|1-\frac t2+\frac{ic}2|$. With $\partial_{\bar z}=\frac12(\partial_1+i\partial_2)$ one has $\partial_{\bar z}(m_1-im_2)=\frac12(t-ic)$, so these eigenvalues are $\pm|w|$. The
$(g_1,g_2)$-block is $(-1+\frac t2)I-\sym A$; with $t=\lambda_1+\lambda_2$ its eigenvalues are $-1+\frac{\lambda_1+\lambda_2}2-\lambda_i=-1\mp\frac12(\lambda_1-\lambda_2)$.\end{proof}

\section{The certificate}\label{sec:cert}
The certificate proves $\kappa(\Om)\le\kappa_0$ by one matrix test (Theorem~\ref{thm:rel}). Near $\partial\Om$ the multiplier of Lemma~\ref{lem:mult} gives a pointwise
inequality; in the interior, $h$ is controlled by finitely many Taylor coefficients on balls.
\subsection{Matrix weight}
Fix $\kappa_0\in(0,1)$, a multiplier $m$ as in Lemma~\ref{lem:mult}, and a bounded measurable $\Phi\ge0$ on $\Om$ (the bump sum \eqref{A:Phi} below). The matrix weight is
\begin{equation}\label{Cmat}C(x):=\kappa_0I-W(m)(x)+\Phi(x)\,e_0e_0^{\mathsf T}.\end{equation}
\begin{lemma}[reformulation]\label{A:lem:reform}
For every $F\in\M$,
\begin{equation}\label{goal3}\int_\Om(h^2-|g|^2)-\kappa_0\int_\Om|F|^2=\int_\Om\Phi\,h^2-\int_\Om F^{\mathsf T}CF .\end{equation}
Consequently $\kappa(\Om)\le\kappa_0$ if and only if $\int_\Om\Phi h^2\le\int_\Om F^{\mathsf T}CF$ for all $F\in\Mz$.
\end{lemma}
\begin{proof}
Integrate $F^{\mathsf T}CF=\kappa_0|F|^2-F^{\mathsf T}W(m)F+\Phi h^2$ (all terms are integrable, as $W(m)$ and $\Phi$ are bounded) and use Lemma~\ref{lem:mult}.
By \eqref{kappa}, $\kappa(\Om)\le\kappa_0$ if and only if the left-hand side of \eqref{goal3} is $\le0$ on $\Mz$.
\end{proof}
Consequently the target is $\int\Phi h^2\le\int F^{\mathsf T}CF$ on $\Mz$. Adding $\Phi e_0e_0^{\mathsf T}$ charges only $h$. Where $Dm=0$, $C$ is positive definite if
$\Phi>1-\kappa_0$; where $\Phi=0$, positivity must come from $\kappa_0I-W(m)$. Remarks on $\tr W(A)$ and on the roles of $\Phi$ and $m$ are in \docD{3}.

\subsection{Ball forms}\label{A:sec:balls}
The term $\int\Phi h^2$ is estimated through the harmonic expansion of $h$ on interior balls. Finitely many Taylor coefficients enter the Gram matrix; the tail is a small
multiple of $\|h\|^2$ on a larger ball (Lemma~\ref{lem:ball}). Let $B_k:=B(x_k,r_k)$, $k=1,\dots,N_B$, be balls with $B(x_k,R_k)\subset\Om$ and $r_k<R_k$, let $\alpha_k\ge0$ be weights, and
\begin{equation}\label{A:Phi}
\phi_k(x):=\Big(1-\frac{|x-x_k|^2}{r_k^2}\Big)_+,\qquad \theta_k:=\frac{r_k}{R_k}<1,\qquad \Phi:=\sum_{k=1}^{N_B}\alpha_k\phi_k .
\end{equation}
For $F=(h,g)\in\M$, the harmonic $h$ has on $B(x_k,\dist(x_k,\partial\Om))\supset B(x_k,R_k)$ the expansion $h=\sum_{d\ge0}h_d$ into homogeneous harmonic polynomials
$h_d$ of degree $d$ in $x-x_k$, converging uniformly on compact subsets \cite[Cor.~5.34]{ABR}. The space $\mathcal Y_d(x_k)$ of homogeneous harmonic polynomials of degree $d$ in $x-x_k$
has dimension $2d+1$ \cite[Prop.~5.8]{ABR}. For $d\le p_k$ let $Y_{k,d,l}$, $l=1,\dots,2d+1$, be a basis of $\mathcal Y_d(x_k)$ that is orthonormal for $\int_\Om\phi_k(\cdot)(\cdot)$, and define the Taylor functional $\ell_{k,d,l}(F)$ as the coefficient of $Y_{k,d,l}$ in $h_d$, so that
\begin{equation}\label{A:ell}
h_d=\sum_{l=1}^{2d+1}\ell_{k,d,l}(F)\,Y_{k,d,l}\qquad(d\le p_k).
\end{equation}
\begin{lemma}\label{lem:ball}
(a) $\int\phi_kh^2=\sum_d\int\phi_kh_d^2$.
(b) For $d\ge0$,
\begin{gather*}\int\phi_kh_d^2=\frac{2\theta_k^{2d+3}}{2d+5}\|h_d\|^2_{L^2(B(x_k,R_k))};\quad\text{consequently}\\
\sum_{d>p}\int\phi_kh_d^2\le\frac{2\theta_k^{2p+5}}{2p+7}\|h\|^2_{L^2(B(x_k,R_k))}.\end{gather*}
(c) $\sum_{d\le p_k}\int\phi_kh_d^2=\sum_{d\le p_k}\sum_{l=1}^{2d+1}\ell_{k,d,l}(F)^2$.
\end{lemma}
\begin{proof}
Homogeneous harmonic polynomials of different degrees are orthogonal in $L^2$ of every sphere centred at $x_k$
\cite[Prop.~5.9]{ABR}, hence for every radial weight. By uniform convergence on $\overline{B(x_k,\rho)}$, $\rho<R$,
$\int_{B(x_k,\rho)}h^2=\sum_d\int_{B(x_k,\rho)}h_d^2$; letting $\rho\uparrow R$ (monotone convergence on both sides) gives $\|h\|^2_{L^2(B(x_k,R))}=\sum_d\|h_d\|^2_{L^2(B(x_k,R))}$.
The same argument with the radial weight $\phi_k$, supported in $\overline{B(x_k,r)}\subset B(x_k,R)$, gives (a). For (b), write $h_d(x_k+sy)=s^dY(y)$, $y\in S^2$. In polar coordinates
\begin{gather*}\int\phi_kh_d^2=\int_0^r\Big(1-\frac{s^2}{r^2}\Big)s^{2d+2}ds\int_{S^2}Y^2=\frac{2r^{2d+3}}{(2d+3)(2d+5)}\int_{S^2}Y^2,\\
\|h_d\|^2_{L^2(B(x_k,R))}=\frac{R^{2d+3}}{2d+3}\int_{S^2}Y^2,\end{gather*}
so the ratio is $2\theta^{2d+3}/(2d+5)$. It is nonincreasing in $d$, so $\sum_{d>p}\int\phi_kh_d^2\le\frac{2\theta^{2p+5}}{2p+7}\sum_{d>p}\|h_d\|^2_{L^2(B(x_k,R))}\le
\frac{2\theta^{2p+5}}{2p+7}\|h\|^2_{L^2(B(x_k,R))}$. Part (c) follows from \eqref{A:ell} and the orthonormality of $Y_{k,d,1},\dots,Y_{k,d,2d+1}$.
\end{proof}

\subsection{Representers}
The proof of Theorem~\ref{thm:rel} bounds the functionals $\ell_\iota$ by the Cauchy--Schwarz inequality in $L^2(\Om)^4$, so each must be the inner product with a fixed
field. One radial weight $\psi_k$ gives the coefficient: harmonic polynomials of different degrees are orthogonal on every sphere about $x_k$, and within one degree every
radial weight gives the same angular inner product up to a positive factor. Fields orthogonal to $\Mz$ (the corrector term) may be subtracted to make $G$ small.
\begin{definition}[representers]\label{A:def:rep}
Let $\iota=(k,d,l)$ with $d\le p_k$, and let $\psi_k\ge0$ be a bounded radial weight about $x_k$, not identically zero, supported in $\overline{B(x_k,R_k')}$ with $0<R_k'<\dist(x_k,\partial\Om)$. Put
\begin{equation}\label{A:bump}
b_\iota:=\mathcal N_\iota^{-1}\,\psi_k\,Y_\iota\,e_0,\qquad \mathcal N_\iota:=\int_\Om\psi_kY_\iota^2 ,
\end{equation}
where $Y_\iota:=Y_{k,d,l}$, and for $\xi=(\varphi,\vartheta)\in H^1_0(\Om)^3\times H^1_0(\Om)$
\begin{equation}\label{A:Dstar}
\mathcal D^*\xi:=(\Div\varphi,\ \curl\varphi+\nabla\vartheta).
\end{equation}
A representer of $\ell_\iota$ is any field
\begin{equation}\label{rep3}r_\iota:=b_\iota-\mathcal D^*\xi_\iota-\gamma_\iota e_0,\qquad\xi_\iota\in H^1_0(\Om)^3\times H^1_0(\Om),\ \gamma_\iota\in\R .\end{equation}
The field $\xi_\iota$ is called the corrector.
\end{definition}
\begin{lemma}\label{A:lem:rep}
(a) $\int_\Om F\cdot b_\iota=\ell_\iota(F)$ for $F\in\M$. (b) $\int_\Om F\cdot\mathcal D^*\xi=0$ for $F\in\M$ and every $\xi$. (c) $\int_\Om F\cdot e_0=0$ for $F\in\Mz$.
Hence $\int_\Om F\cdot r_\iota=\ell_\iota(F)$ for all $F\in\Mz$, for every corrector $\xi_\iota$ and every $\gamma_\iota$.
\end{lemma}
\begin{proof}
(a) On $\supp\psi_k\subset\overline{B(x_k,R_k')}$ the expansion of $h$ converges uniformly and may be integrated term by term; only $h_d$ contributes (proof of
Lemma~\ref{lem:ball}). For $Y,Y'\in\mathcal Y_d(x_k)$, polar coordinates give $\int\psi_kYY'=c_\psi\int_{S^2}YY'$ and $\int\phi_kYY'=c_\phi\int_{S^2}YY'$ with
$c_\psi,c_\phi>0$, so the basis $Y_{k,d,l}$ is orthogonal for $\psi_k$. With \eqref{A:ell},
\[\int_\Om F\cdot b_\iota=\mathcal N_\iota^{-1}\int_\Om\psi_k\,h\,Y_\iota=\mathcal N_\iota^{-1}\sum_{l'}\ell_{k,d,l'}(F)\int_\Om\psi_kY_{k,d,l'}Y_\iota=\ell_\iota(F).\]
(b) $\int F\cdot\mathcal D^*\xi=(h,\Div\varphi)+(g,\curl\varphi)+(g,\nabla\vartheta)$. The first two terms cancel by \eqref{dual0}, and $(g,\nabla\vartheta)=0$ because $\Div g=0$ in $\mathcal D'(\Om)$ and
$C_c^\infty(\Om)$ is dense in $H^1_0(\Om)$. (c) $\int F\cdot e_0=\int h=0$.
\end{proof}
\begin{remark}[choices in the verifier]\label{A:rem:repchoice}
The proof allows any corrector. The verifier uses $\psi_k(x)=(1-|x-x_k|^2/R'^2)_+^{\mu}$ with $\mu=2$ and $R'=\frac{63}{64}R_k$, and as corrector a finite element
approximation of the minimiser of $\int r_\iota^{\mathsf T}(C')^{-1}r_\iota$, computed in floating point. A poor corrector cannot produce a false certificate; it only makes $G$ larger.
\end{remark}

\subsection{The criterion}
The certificate data collect all choices. Theorem~\ref{thm:rel} needs only the pointwise inequalities \eqref{A:data} and the eigenvalue test.
\begin{definition}[certificate data]\label{A:def:data}
The data of a certificate are: $\kappa_0\in(0,1)$; a multiplier $m$ as in Lemma~\ref{lem:mult}; balls, weights and degrees $(x_k,r_k,R_k,\alpha_k,p_k)$, $k=1,\dots,N_B$, as in
Section~\ref{A:sec:balls}, with $\Phi$ from \eqref{A:Phi} and $C$ from \eqref{Cmat}; the index set
$\mathcal I:=\{(k,d,l):1\le k\le N_B,\ 0\le d\le p_k,\ 1\le l\le2d+1\}$ with $N:=\#\mathcal I=\sum_k(p_k+1)^2$ and weights $\alpha_\iota:=\alpha_k$ for $\iota=(k,d,l)$; representers $r_\iota$,
$\iota\in\mathcal I$, of the form \eqref{rep3}; a function $\rho_{\rm up}$ and a measurable symmetric matrix field $C'$ with
\begin{equation}\label{A:data}
\rho_{\rm up}\ge\sum_{k=1}^{N_B}\alpha_k\frac{2\theta_k^{2p_k+5}}{2p_k+7}\mathbf 1_{B(x_k,R_k)},\qquad 0\prec C'\preceq C-\rho_{\rm up}e_0e_0^{\mathsf T}\quad\text{a.e.\ in }\Om,
\end{equation}
such that the integrals in \eqref{G} below are finite. The Gram matrix of the certificate is
\begin{equation}\label{G}G_{\iota\iota'}:=\sqrt{\alpha_\iota\alpha_{\iota'}}\int_\Om r_\iota^{\mathsf T}(C')^{-1}r_{\iota'},\qquad\iota,\iota'\in\mathcal I .\end{equation}
\end{definition}
\begin{theorem}[reliability]\label{thm:rel}
Let the data of Definition~\ref{A:def:data} be given. If $\lambda_{\max}(G)\le1$ (or $\lambda_{\max}(G_{\rm up})\le1$ for some $G_{\rm up}\succeq G$), then $\kappa(\Om)\le\kappa_0$, i.e.\
$\beta(\Om)\ge\sqrt{(1-\kappa_0)/2}$.
\end{theorem}
\begin{proof}
Let $F\in\Mz$ and $u\in\R^N$, $u_\iota:=\sqrt{\alpha_\iota}\,\ell_\iota(F)=\int F\cdot\sqrt{\alpha_\iota}r_\iota$ (Lemma~\ref{A:lem:rep}). For $y\in\R^N$ put $r_y:=\sum_\iota y_\iota\sqrt{\alpha_\iota}r_\iota$. By the Cauchy--Schwarz inequality with weight $C'$,
\begin{align*}(y^{\mathsf T}u)^2&=\Big(\int(C')^{1/2}F\cdot(C')^{-1/2}r_y\Big)^2\\&\le\int F^{\mathsf T}C'F\int r_y^{\mathsf T}(C')^{-1}r_y=\int F^{\mathsf T}C'F\cdot y^{\mathsf T}Gy ;\end{align*}
with $y=u$ this gives $|u|^2\le\lambda_{\max}(G)\int F^{\mathsf T}C'F\le\int F^{\mathsf T}C'F$. By Lemma~\ref{lem:ball} and \eqref{A:data},
\begin{align*}\int\Phi h^2&\le\sum_\iota\alpha_\iota\ell_\iota(F)^2+\int\rho_{\rm up}h^2=|u|^2+\int\rho_{\rm up}h^2\\
&\le\int F^{\mathsf T}(C'+\rho_{\rm up}e_0e_0^{\mathsf T})F\le\int F^{\mathsf T}CF,\end{align*}
and Lemma~\ref{A:lem:reform} gives $\|h\|^2-\|g\|^2\le\kappa_0\|F\|^2$. The same argument applies with $G_{\rm up}$, since $y^{\mathsf T}Gy\le y^{\mathsf T}G_{\rm up}y$. The last statement is Theorem~\ref{thm:fv}.
\end{proof}
\begin{remark}\label{A:rem:clow}
A simple choice is $C':=C_{\rm low}-\rho_{\rm up}e_0e_0^{\mathsf T}$ with a piecewise constant symmetric matrix field $C_{\rm low}$ satisfying
\begin{equation}\label{A:Clow}
0\prec C_{\rm low}\preceq C\quad\text{a.e.\ in }\Om ;
\end{equation}
the designs of Section~\ref{sec:lower} use it, the exact verification the weight of Section~\ref{sec:verif}.
\end{remark}

\subsection{Symmetry reduction}
Symmetries of $\Om$ split $\Mz$ into orthogonal classes. The test can then be done class by class, with smaller Gram matrices and integrals over a fundamental region;
Proposition~\ref{oth:prop:classes} shows that this reduction is exact. For $R\in O(3)$ let
\begin{equation}\label{A:TR}
\Lambda_R:=\operatorname{diag}(\det R,R)\in O(4),\quad (\mathcal T_RF)(x):=\Lambda_R^{\mathsf T}F(Rx)=\big(\det R\,h(Rx),\,R^{\mathsf T}g(Rx)\big).
\end{equation}
\begin{lemma}[symmetry]\label{A:lem:sym}
Let $R\in O(3)$ with $R\Om=\Om$.
(a) $\mathcal T_R$ is unitary on $L^2(\Om)^4$, maps $\M$ onto $\M$ and $\Mz$ onto $\Mz$, and $\|(\mathcal T_RF)_h\|=\|h\|$, $\|(\mathcal T_RF)_g\|=\|g\|$, where the subscripts denote the
first and the last three components.
(b) If $m(Rx)=Rm(x)$, then
\[W(m)(Rx)=\Lambda_RW(m)(x)\Lambda_R^{\mathsf T},\qquad\int_\Om(\mathcal T_RF)^{\mathsf T}W(m)(\mathcal T_RF)=\int_\Om F^{\mathsf T}W(m)F.\]
\end{lemma}
\begin{proof}
(a) $x\mapsto Rx$ is a measure-preserving bijection of $\Om$ and $\Lambda_R$ is orthogonal, so $\mathcal T_R$ is unitary with inverse $\mathcal T_{R^{-1}}$. For
$F\in\M$ (smooth by Lemma~\ref{A:lem:reg}) put $\tilde h(x):=\det R\,h(Rx)$, $\tilde g(x):=R^{\mathsf T}g(Rx)$. The chain rule and the identity
\[\operatorname{ax}(R^{\mathsf T}AR)=\det R\,R^{\mathsf T}\operatorname{ax}(A)\qquad(R\in O(3))\]
give $\nabla\tilde h(x)=\det R\,R^{\mathsf T}(\nabla h)(Rx)=\curl\tilde g(x)$ and $\Div\tilde g(x)=(\Div g)(Rx)=0$; moreover $\int_\Om\tilde h=\det R\int_\Om h$.
(b) $Dm(Rx)=R\,Dm(x)R^{\mathsf T}$, so $\tr Dm$ is invariant, $\sym Dm(Rx)=R\,\sym Dm(x)R^{\mathsf T}$ and $\curl m(Rx)=\det R\,R\curl m(x)$; the blocks of
\eqref{W} give $W(m)(Rx)=\Lambda_RW(m)(x)\Lambda_R^{\mathsf T}$, and the change of variables $x\mapsto Rx$ gives the second claim. A detailed version is in
\docD{2}.
\end{proof}
\begin{lemma}[similarities]\label{lem:simil}
Let $\Psi(x):=rRx+b$ with $r>0$, $R\in SO(3)$ and $b\in\R^3$. Then $(h,g)\mapsto(h\circ\Psi^{-1},Rg\circ\Psi^{-1})$ maps $\M(\Om)$ onto $\M(\Psi\Om)$ and $\Mz(\Om)$
onto $\Mz(\Psi\Om)$, and multiplies $\|h\|^2$ and $\|g\|^2$ by $r^3$. The map $q\mapsto q\circ\Psi^{-1}$ is $r^{3/2}$ times a unitary map from $L^2_0(\Om)$ onto
$L^2_0(\Psi\Om)$, and $\mathcal S_{\Psi\Om}(q\circ\Psi^{-1})=(\mathcal S_\Om q)\circ\Psi^{-1}$. Hence $\kappa$, $\beta$ and the spectrum of $\mathcal S$, with its essential
part, are invariant under $\Psi$.
\end{lemma}
\begin{proof}
This is the change of variables $y=\Psi(x)$, with Jacobian $r^3$, in \eqref{kappa} and in Definition~\ref{A:def:S}, where $v_{q\circ\Psi^{-1}}=r\,Rv_q\circ\Psi^{-1}$;
$\curl(Rg\circ\Psi^{-1})=r^{-1}R(\curl g)\circ\Psi^{-1}$ by the identity for $\operatorname{ax}$ in the proof of Lemma~\ref{A:lem:sym}. Details are in
\docD{2}.
\end{proof}

\begin{corollary}[class-wise criterion]\label{cor:class}
Let $X\subset\Mz$ be a closed subspace, and let the data of Definition~\ref{A:def:data} be given, except that each $r_\iota$ is only required to satisfy $\int F\cdot r_\iota=\ell_\iota(F)$ for $F\in X$.
If $\lambda_{\max}(G)\le1$, then $\|h\|^2-\|g\|^2\le\kappa_0\|F\|^2$ for all $F\in X$.
\end{corollary}
\begin{proof}
The proof of Theorem~\ref{thm:rel} uses the representers only through $u_\iota=\int F\cdot\sqrt{\alpha_\iota}r_\iota$ for the pair $F$ at hand, so it applies to $F\in X$.
\end{proof}
\begin{definition}[symmetry classes]\label{A:def:classes}
Let $\mathcal G$ be a finite group of orthogonal $3\times3$ matrices $R$ with $R\Om=\Om$ and $R^2=I$ for every $R\in\mathcal G$. Such a group is abelian, because
$RR'=(RR')^{-1}=R'^{-1}R^{-1}=R'R$, and every $R\in\mathcal G$ is symmetric. A \emph{character} is a homomorphism $\chi:\mathcal G\to\{\pm1\}$. Put
\begin{equation}\label{oth:eq:class}
\begin{gathered}
\Mz^\chi:=\{F\in\Mz:\ \mathcal T_RF=\chi(R)F\ \text{for all }R\in\mathcal G\},\qquad P_\chi:=\frac1{|\mathcal G|}\sum_{R\in\mathcal G}\chi(R)\,\mathcal T_R ,\\
\kappa_\chi:=\sup_{F\in\Mz^\chi\setminus\{0\}}\frac{\|h\|^2-\|g\|^2}{\|h\|^2+\|g\|^2},
\end{gathered}
\end{equation}
with $\kappa_\chi:=-\infty$ if $\Mz^\chi=\{0\}$. A \emph{fundamental region} is a closed polyhedral set $\Om_{\mathcal G}\subset\overline\Om$ such that
$\overline\Om=\bigcup_{R\in\mathcal G}R\,\Om_{\mathcal G}$ and $|R\,\Om_{\mathcal G}\cap R'\,\Om_{\mathcal G}|=0$ for $R\ne R'$. The \emph{orbit size} $\operatorname{mult}_k$ of a
ball $B_k$ is the number of distinct balls $RB_k$, $R\in\mathcal G$. A corrector $\xi$ has the \emph{parity of the class} $\chi$ if
\begin{equation}\label{A:parity}
\mathcal T_R\,\mathcal D^*\xi=\chi(R)\,\mathcal D^*\xi\qquad\text{for all }R\in\mathcal G .
\end{equation}
\end{definition}

\begin{proposition}[class reduction for groups of reflections]\label{oth:prop:classes}\label{A:prop:cubeclass}
Let $\mathcal G$ and $\Om_{\mathcal G}$ be as in Definition~\ref{A:def:classes}.
\begin{enumerate}
\item[(a)] $\Mz=\bigoplus_\chi\Mz^\chi$, and the decomposition is orthogonal both for $\int h\,h'$ and for $\int g\cdot g'$. Consequently $\kappa(\Om)=\max_\chi\kappa_\chi$.
\item[(b)] Let the data of Definition~\ref{A:def:data} be given, and assume: (i) $m(Rx)=Rm(x)$ for all $R\in\mathcal G$; (ii) the balls form full $\mathcal G$-orbits, with equal
weights $\alpha_k$ and equal $r_k,R_k,p_k$ within an orbit, so that $\Phi$ and $\rho_{\rm up}$ are $\mathcal G$-invariant; (iii) $C'$ is given on $\Om_{\mathcal G}$ and
extended by $C'(Rx):=\Lambda_RC'(x)\Lambda_R$; (iv) the representers are $r_\iota^\chi:=P_\chi r_\iota$ with $r_\iota$ of the form \eqref{rep3}. Then $\int F\cdot r_\iota^\chi=\ell_\iota(F)$ for
$F\in\Mz^\chi$; the inequalities $0\prec C'\preceq C-\rho_{\rm up}e_0e_0^{\mathsf T}$ hold almost everywhere in $\Om$ as soon as they hold on $\Om_{\mathcal G}$; the integrands
$(r_\iota^\chi)^{\mathsf T}(C')^{-1}r_{\iota'}^\chi$ are $\mathcal G$-invariant, so that $\int_\Om=|\mathcal G|\int_{\Om_{\mathcal G}}$; and
$\lambda_{\max}(G^\chi)\le1$ for the Gram matrix $G^\chi$ of \eqref{G} built from the $r_\iota^\chi$ implies $\kappa_\chi\le\kappa_0$.
\item[(c)] In (b) the Gram matrix may be formed from one ball per $\mathcal G$-orbit with the weight $\operatorname{mult}_k\alpha_k$; functionals that vanish identically on
$\Mz^\chi$ may be omitted, and two functionals that agree up to sign on $\Mz^\chi$ may be replaced by one of them with the sum of their weights.
\end{enumerate}
\end{proposition}
\begin{proof}
The argument is written for orthogonal $R$ with $R\Om=\Om$; it uses $R^2=I$ only where stated, and for involutions $\Lambda_R^{\mathsf T}=\Lambda_R$, $R^{\mathsf T}=R$.

(a) By Lemma~\ref{A:lem:sym}(a) and \eqref{A:TR}, each $\mathcal T_R$ is unitary, $\mathcal T_R^*=\mathcal T_{R^{-1}}$, $\mathcal T_A\mathcal T_B=\mathcal T_{BA}$, and
$\mathcal T_R$ acts on $h$ and $g$ separately and preserves $\Mz$. As $\chi(R^{-1})=\chi(R)$, $P_\chi^*=P_\chi$. Substituting $S\mapsto SR$ in the sum gives
$\mathcal T_RP_\chi=\chi(R)P_\chi$, hence $P_\chi^2=P_\chi$. So $P_\chi$ is an orthogonal projection acting on $h$ and $g$ separately, $P_\chi F=F$ on $\Mz^\chi$ and
$P_\chi\Mz=\Mz^\chi$. Now let $R^2=I$ for all $R$. Then the characters are all irreducible representations of the abelian group $\mathcal G$, and orthogonality of
characters gives $P_\chi P_{\chi'}=0$ for $\chi\ne\chi'$ and $\sum_\chi P_\chi=I$; this is the orthogonal decomposition, for both forms. If $F=\sum_\chi F_\chi$,
$F_\chi=(h_\chi,g_\chi)\in\Mz^\chi$, put $a_\chi:=\|h_\chi\|^2$, $b_\chi:=\|g_\chi\|^2$. Then the quotient of $F$ is
$\sum_\chi(a_\chi-b_\chi)/\sum_\chi(a_\chi+b_\chi)\le\max_\chi(a_\chi-b_\chi)/(a_\chi+b_\chi)$ (maximum over $a_\chi+b_\chi>0$), so $\kappa(\Om)\le\max_\chi\kappa_\chi$;
the converse holds as $\Mz^\chi\subset\Mz$.

(b) By Corollary~\ref{cor:class} with $X=\Mz^\chi$, three facts suffice. First, $\int F\cdot P_\chi r_\iota=\int P_\chi F\cdot r_\iota=\ell_\iota(F)$ for
$F\in\Mz^\chi$ by (a) and Lemma~\ref{A:lem:rep}. Second, by (i) and Lemma~\ref{A:lem:sym}(b), $W(m)(Rx)=\Lambda_RW(m)(x)\Lambda_R^{\mathsf T}$; by (ii) and
$\Lambda_Re_0e_0^{\mathsf T}\Lambda_R^{\mathsf T}=e_0e_0^{\mathsf T}$ the same holds for $C-\rho_{\rm up}e_0e_0^{\mathsf T}$, and by (iii) for $C'$. Congruence with
the orthogonal $\Lambda_R$ shows that $0\prec C'\preceq C-\rho_{\rm up}e_0e_0^{\mathsf T}$ holds on $\Om$ (up to the null set of the boundaries of the $R\,\Om_{\mathcal G}$)
once it holds on $\Om_{\mathcal G}$. Third, $\mathcal T_Rr_\iota^\chi=\chi(R)r_\iota^\chi$ by (a), i.e.\ $r_\iota^\chi(Rx)=\chi(R)\Lambda_Rr_\iota^\chi(x)$, and
$C'(Rx)^{-1}=\Lambda_RC'(x)^{-1}\Lambda_R^{\mathsf T}$. As $\chi(R)^2=1$ and $\Lambda_R^{\mathsf T}\Lambda_R=I$, $(r_\iota^\chi)^{\mathsf T}(C')^{-1}r_{\iota'}^\chi$
takes the same value at $x$ and $Rx$, and the tiling gives $\int_\Om=|\mathcal G|\int_{\Om_{\mathcal G}}$.

(c) Let $F\in\Mz^\chi$ and $c_R:=\chi(R)\det R$. Then $h(Rx)=c_Rh(x)$, so $h(Rx_k+z)=c_R\,h(x_k+R^{\mathsf T}z)$, and the homogeneous parts of the Taylor
expansions about $Rx_k$ and $x_k$ satisfy $h^{(Rx_k)}_d(z)=c_R\,h^{(x_k)}_d(R^{\mathsf T}z)$. The bumps are radial, so the integral of $\big(h^{(Rx_k)}_d\big)^2$ against
the bump of $RB_k$ equals $\int\phi_k\big(h^{(x_k)}_d\big)^2$ for every $d$. Summing over the distinct balls of an orbit and using Lemma~\ref{lem:ball}(c), the
orbit contributes $\operatorname{mult}_k\alpha_k\sum_{d\le p_k,\,l}\ell_{k,d,l}(F)^2$ to $\sum_\iota\alpha_\iota\ell_\iota(F)^2$, for any ball $B_k$ of the orbit; this sum
is the only way in which the functionals enter the proof of Theorem~\ref{thm:rel}. Omitting a functional with $\ell_\iota=0$ on $\Mz^\chi$ changes neither side,
and $w\,\ell(F)^2+w'\,\ell'(F)^2=(w+w')\,\ell(F)^2$ if $\ell'=\pm\ell$ on $\Mz^\chi$.
\end{proof}
If the corrector $\xi_\iota$ has the parity \eqref{A:parity} of the class, then $P_\chi\mathcal D^*\xi_\iota=\mathcal D^*\xi_\iota$, so only $b_\iota$ and $\gamma_\iota e_0$ have to
be projected. Symmetries that are not involutions are used through the following lemma.
\begin{lemma}[conjugate classes]\label{oth:lem:normal}
Let $\mathcal G$ be as in Definition~\ref{A:def:classes} and let $\Pi$ be an orthogonal matrix with $\Pi\Om=\Om$ and $\Pi\mathcal G\Pi^{-1}=\mathcal G$. For a character $\chi$ put
$\chi^\Pi(R):=\chi(\Pi R\Pi^{-1})$. Then $\mathcal T_\Pi$ maps $\Mz^\chi$ bijectively onto $\Mz^{\chi^\Pi}$ and preserves $\|h\|^2$ and $\|g\|^2$; hence $\kappa_\chi=\kappa_{\chi^\Pi}$.
\end{lemma}
\begin{proof}
Let $F\in\Mz^\chi$, $R\in\mathcal G$ and $R':=\Pi R\Pi^{-1}\in\mathcal G$, so that $\Pi R=R'\Pi$. By the composition rule $\mathcal T_A\mathcal T_B=\mathcal T_{BA}$,
$\mathcal T_R(\mathcal T_\Pi F)=\mathcal T_{\Pi R}F=\mathcal T_{R'\Pi}F=\mathcal T_\Pi(\mathcal T_{R'}F)=\chi(R')\,\mathcal T_\Pi F=\chi^\Pi(R)\,\mathcal T_\Pi F$, so
$\mathcal T_\Pi F\in\Mz^{\chi^\Pi}$. $\mathcal T_\Pi$ preserves $\Mz$ and both norms (Lemma~\ref{A:lem:sym}(a)), and its inverse $\mathcal T_{\Pi^{-1}}$ maps $\Mz^{\chi^\Pi}$
back into $\Mz^\chi$ by the same argument, since $(\chi^\Pi)^{\Pi^{-1}}=\chi$.
\end{proof}

\emph{The cube.} Let $\Om=(-1,1)^3$. For $\varsigma\in\{\pm1\}^3$ let $J_\varsigma:=\operatorname{diag}(\varsigma_1,\varsigma_2,\varsigma_3)$, and let $J_i$ be the reflection in $x_i$
($\varsigma_i=-1$, $\varsigma_j=1$ for $j\ne i$). The eight reflections form a group $\mathcal G_{\rm C}$ as in Definition~\ref{A:def:classes}, with fundamental region the octant
$[0,1]^3$. Its characters are $\chi_s(J_\varsigma):=\prod_i\varsigma_i^{s_i}$, $s\in\{0,1\}^3$, and we write
\begin{gather}
\Mz^s:=\Mz^{\chi_s},\qquad \kappa_s:=\kappa_{\chi_s},\label{A:classes}\\
P_s:=P_{\chi_s}=\frac18\sum_{\varsigma\in\{\pm1\}^3}\chi_s(J_\varsigma)\,\mathcal T_{J_\varsigma}.\label{A:Ps}
\end{gather}
\begin{corollary}[symmetry classes of the cube]\label{A:lem:classes}
Let $\Om=(-1,1)^3$.
(a) $F\in\Mz^s$ if and only if $F\in\Mz$ and
\[h(J_ix)=-(-1)^{s_i}h(x),\qquad g_j(J_ix)=(-1)^{s_i+\delta_{ij}}g_j(x),\qquad i,j=1,2,3.\]
(b) A coordinate permutation maps $\Mz^s$ onto $\Mz^{s'}$, where $s'$ is the correspondingly permuted index, isometrically for both $\|h\|^2$ and $\|g\|^2$. Hence $\kappa_s=\kappa_{s'}$, and
\[\kappa(\Om)=\max\{\kappa_{000},\kappa_{100},\kappa_{110},\kappa_{111}\}.\]
\end{corollary}
\begin{proof}
(a) $\mathcal T_{J_i}F(x)=(-h(J_ix),J_ig(J_ix))$ and $\chi_s(J_i)=(-1)^{s_i}$; the eight reflections are generated by $J_1,J_2,J_3$.
(b) A permutation matrix $\Pi$ maps the cube onto itself, and $\Pi J_\varsigma\Pi^{-1}$ is the reflection with correspondingly permuted signs; hence
$\Pi\mathcal G_{\rm C}\Pi^{-1}=\mathcal G_{\rm C}$ and $\chi_s^\Pi=\chi_{s'}$, and Lemma~\ref{oth:lem:normal} gives the first claim and $\kappa_s=\kappa_{s'}$. Permutations map
indices onto indices with the same number of ones, so the four listed classes represent all eight, and Proposition~\ref{oth:prop:classes}(a) gives the formula for $\kappa(\Om)$.
\end{proof}
Only class $111$ contains functions $h$ with nonzero mean, so $\gamma_\iota e_0$ in \eqref{rep3} is needed only there; a ball centred on a symmetry plane
contributes only the harmonic polynomials of the parity of the class.
\begin{remark}[checks in the verifier]\label{A:rem:verifcheck}
For the cube, the hypotheses (i)--(iv) of Proposition~\ref{oth:prop:classes}(b) with $\mathcal G=\mathcal G_{\rm C}$ are checked as follows: (i) exactly at the mesh vertices
($m_i=0$ on $\{x_i=0\}$, and $m=0$ on $\{x_i=1\}$); (ii) holds by construction, because the ball list is generated as full reflection orbits with one weight per orbit; (iii) is how
$C'$ is assembled; in (iv) the projection coefficients are exact multiples of $\frac18$, and the correctors have the parity \eqref{A:parity} of the class, imposed by fixing their
nodal values on the symmetry planes. The fixed values are checked to vanish.
\end{remark}

\section{Upper bounds for the inf-sup constant}\label{sec:upper}
An upper bound for $\beta$ needs one admissible pair with a large quotient. For $F=(h,g)$ with $h\ne0$ or $g\ne0$ we write
\begin{equation}\label{B-eq:R}
\mathcal R(h,g):=\frac{\|h\|^2-\|g\|^2}{\|h\|^2+\|g\|^2}
\end{equation}
for its quotient. By \eqref{kappa} every $F\in\Mz\setminus\{0\}$ gives $\kappa(\Om)\ge\mathcal R(F)$ and, by Theorem~\ref{thm:fv}, $\beta(\Om)\le\sqrt{(1-\mathcal R(F))/2}$.
Trial pairs are constructed in $\M$ and then corrected by their mean.
\begin{lemma}[mean correction]\label{B-lem:mean}
Let $\Om$ be bounded, $F=(h,g)\in\M$ and
\begin{equation}\label{B-eq:hbar}
\bar h:=|\Om|^{-1}\int_\Om h .
\end{equation}
Then $(h-\bar h,g)\in\Mz$ and $\|h-\bar h\|^2=\|h\|^2-|\Om|\,\bar h^2$.
\end{lemma}
\begin{proof}
Constant pairs $(c,0)$ belong to $\M$, because $\nabla c=0=\curl0$ and $\Div0=0$, and $\M$ is a linear space; $\int_\Om(h-\bar h)=0$ by \eqref{B-eq:hbar}. Expanding,
$\|h-\bar h\|^2=\|h\|^2-2\bar h\int_\Om h+|\Om|\bar h^2=\|h\|^2-|\Om|\bar h^2$.
\end{proof}

\subsection{Polynomial pairs}
For $k\ge0$ let $\mathcal H_{\le k}$ be the space of harmonic polynomials on $\R^3$ of degree at most $k$, and for $d\ge0$ let
\begin{equation}\label{B-eq:Md}
\M_d:=\{F=(h,g)\in\M:\ h,g_1,g_2,g_3\ \text{are polynomials of degree}\le d\}.
\end{equation}
For $h\in\mathcal H_{\le d}$ define the Poincar\'e-gauge field (the classical Poincar\'e homotopy operator, see e.g.\ \cite{CostabelMcIntosh})
\begin{equation}\label{B-eq:gauge}
g_h(x):=\int_0^1t\,(\nabla h)(tx)\times x\,dt .
\end{equation}
\begin{lemma}[polynomial pairs]\label{B-lem:poly}
Let $d\ge0$.
(a) For $h\in\mathcal H_{\le d}$, $g_h$ is a polynomial field of degree at most $d$ with $\curl g_h=\nabla h$ and $\Div g_h=0$.
(b) $\M_d=\{(h,\,g_h+\nabla\vartheta):\ h\in\mathcal H_{\le d},\ \vartheta\in\mathcal H_{\le d+1}\}$, and $\dim\M_d=2(d+1)(d+2)$.
(c) The functional $F\mapsto\int_\Om h$ does not vanish on $\M_d$; hence $\dim(\M_d\cap\Mz)=2(d+1)(d+2)-1$.
(d) For $F=(h,g_1,g_2,g_3)$ put $\mathcal AF:=(\partial_1g_2-\partial_2g_1,\ \partial_2h+\partial_1g_3,\ \partial_2g_3-\partial_1h,\ -\partial_1g_1-\partial_2g_2)$. The trace map
$F\mapsto F|_{x_3=0}$ is a bijection from $\M_d$ onto the $4$-tuples of polynomials of degree at most $d$ in $(x_1,x_2)$, with inverse, the \emph{Cauchy--Kovalevskaya
extension},
\[F_0\mapsto\sum_{k=0}^d\frac{x_3^k}{k!}\,\mathcal A^kF_0 .\]
\end{lemma}
\begin{proof}
(a) and (b) are the polynomial case of the Poincar\'e homotopy formulas \cite{CostabelMcIntosh}, with $\dim\mathcal H_{\le k}=(k+1)^2$ \cite[Prop.~5.8]{ABR};
(c) holds because $(1,0)\in\M_d$; (d) is the Cauchy--Kovalevskaya expansion of the system $\partial_3F=\mathcal AF$, which is equivalent to $\nabla h=\curl g$,
$\Div g=0$, and $\mathcal A$ lowers the degree. A detailed proof is in \docD{2}.
\end{proof}
Already degree seven gives a pair that can be checked by hand.
\begin{proposition}[an explicit pair]\label{prop:explicit7}
Let $D:=(-1,1)^3$, $\langle f\rangle_D:=|D|^{-1}\int_Df$, $p(w):=-4w+4w^3+2w^5-w^7$, $u+iv:=p(x+it)$ for real $x,t$, and
\[h:=u(x,y)+u(x,z),\qquad g:=\big(0,-v(x,z),v(x,y)\big).\]
Then $(h,g)\in\Mz(D)$, $\langle h^2\rangle_D=\frac{13064512}{135135}$, $\langle|g|^2\rangle_D=\frac{2897056}{135135}$, and
\[\beta(D)^2\le\frac{\langle|g|^2\rangle_D}{\langle h^2\rangle_D+\langle|g|^2\rangle_D}=\frac{90533}{498799}<\frac{57}{314}<\frac12-\frac1\pi .\]
\end{proposition}
\begin{proof}
By the Cauchy--Riemann equations $u_x=v_t$, $u_t=-v_x$,
\[\curl g=\big(v_t(x,y)+v_t(x,z),\,-v_x(x,y),\,-v_x(x,z)\big)=\nabla h,\]
and $\Div g=\partial_y(-v(x,z))+\partial_z v(x,y)=0$. Since $p$ is odd with real coefficients, $u$ is odd in $x$; so $\int_Dh=0$ and $(h,g)\in\Mz(D)$.
The averages follow by expansion into monomials: $\langle x^ay^bz^c\rangle_D$ is $0$ if an exponent is odd and $1/((a+1)(b+1)(c+1))$ otherwise.
By \eqref{kappa} and Theorem~\ref{thm:fv}, $\beta^2=\frac12(1-\kappa)\le\frac12(1-\mathcal R(h,g))$, which is the quotient of averages. Finally $90533\cdot314<57\cdot498799$, and $\pi>\frac{157}{50}$ gives $\frac12-\frac1\pi>\frac12-\frac{50}{157}=\frac{57}{314}$.
\end{proof}
This pair places $\beta(D)^2=\min\sigma(\mathcal S)$ below the edge interval (Theorem~\ref{B-thm:edge}) without tables; for $\sigma_{\rm ess}(\mathcal S)$, Section~\ref{sec:ess} needs a sequence of pairs.

\emph{Higher degrees.} The proof needs only an admissible pair and its exact quotient. The reported computation finds the pair numerically: in each symmetry class the Rayleigh--Ritz vector for \eqref{B-eq:R} on $\M_d\cap\Mz$ is computed in floating point and rounded to rational
coefficients. The quotient of the mean-corrected pair (Lemma~\ref{B-lem:mean}) is then an exact integral of polynomials over boxes.
On the cube the largest values occur in the class $s=011$, which a coordinate permutation maps onto the class $110$ of Section~\ref{sec:cert}. Table~\ref{tab:poly}
lists the cube value for $d=48$, together with the bounds for the box, the L-shaped prism and the Fichera domain; the latter are proved in
Propositions~\ref{oth:prop:boxup}, \ref{oth:prop:lpup} and~\ref{oth:prop:fichup}. The values for lower degrees are listed in
\docD{4}.
\begin{table}[tbp]
\caption{Lower bounds for $\kappa$ from polynomial trial pairs of degree $d$ (exact quotients) and the resulting upper bounds for $\beta$. $^\ast$: degree $36$
enriched by pairs singular at the re-entrant edge (exponents from $\sin(\lambda\omega)=\pm\lambda\sin\omega/\kappa^\sharp$, $\omega=\frac{3\pi}2$, $\kappa^\sharp=0.913$,
Section~\ref{oth:ssec:lp}; polynomial in $z$ up to degree $8$); Gram entries enclosed in ball arithmetic.}\label{tab:poly}
\begin{tabular}{llll}
\toprule
domain & $d$ & $\kappa\ge$ & $\beta\le$\\
\midrule
cube $(-1,1)^3$ & 48 & $0.743070318916$ & $0.358420$\\
box $(0,2)\times(0,1)^2$ & 64 & $0.856921056806$ & $0.267469$\\
L-prism $\Om_{\rm L}=\omega_{\rm L}\times(0,1)$ \eqref{oth:eq:domains} & 24 & $0.904417176994$ & $0.218613$\\
 & $36^\ast$ & $0.913334979008$ & $0.208165$\\
Fichera domain $(-1,1)^3\setminus[0,1]^3$ & 20 & $0.802453899040$ & $0.314282$\\
\bottomrule
\end{tabular}
\end{table}
\begin{proposition}[polynomial trial pairs]\label{B-prop:poly}
For each row of Table~\ref{tab:poly} there is a pair $F\in\Mz(\Om)$ on the domain $\Om$ of that row with $\mathcal R(F)$ at least the listed value; consequently $\kappa(\Om)$ and $\beta(\Om)$ satisfy the listed bounds.
\end{proposition}
\begin{proof}[Proof (computer-assisted)]
For the rows without $^\ast$, let $F$ be the element of $\M_d$ with the rational coefficients of records R1 and R2 of \cite{lbb3drepo}, corrected by its mean (Lemma~\ref{B-lem:mean}). It is admissible by Lemma~\ref{B-lem:poly}, and its quotient is evaluated in exact rational arithmetic; for the cube, admissibility and quotient are also checked in Lean (Section~\ref{sec:verif}).
For the row $36^\ast$, admissibility of the edge-singular pairs is exact and the Gram entries, hence the quotient, are enclosed in ball arithmetic. The bounds for $\beta$ follow from Theorem~\ref{thm:fv}, rounded up.
\end{proof}
The polynomial quotients increase slowly with $d$ (\docD{4}). Section~\ref{sec:vertex} suggests that the vertices are responsible; we therefore build pairs singular at a vertex.

\subsection{Vertex-singular pairs}\label{B-sec:vpairs}
These pairs are used twice: one exponent close to $-\frac32$ gives the upper bound for the cube (Theorem~\ref{thm:upper}); a sequence of exponents tending to $-\frac32$ gives the bound on the essential spectrum (Theorem~\ref{thm:ess}). Throughout,
\begin{equation}\label{B-eq:GQ}
\Gamma:=(0,\infty)^3,\qquad Q:=(0,1)^3,\qquad \Sigma_i:=\{y\in[0,1]^3:\ y_i=1\}\quad(i=1,2,3);
\end{equation}
the $\Sigma_i$ are the three faces of $Q$ that do not contain the vertex $0$, with surface measure $dS$. By Lemma~\ref{lem:simil}, $\kappa((-1,1)^3)=\kappa(Q)$, and we work on $Q$ with the vertex at the origin.

\begin{lemma}[face formula]\label{lem:face}\label{oth:lem:planar}
Let $n\in\{2,3\}$, let $E_1,\dots,E_I\subset\R^n$ be relatively open flat pieces, $E_i$ lying in a hyperplane at distance $\delta_i>0$ from $0$, such that the cones
$K_i:=\{ty:\ 0<t<1,\ y\in E_i\}$ are pairwise disjoint, and let $K$ coincide with $\bigcup_iK_i$ up to a null set. Let $f$ be measurable on $\{ty:\ t>0,\ y\in\bigcup_iE_i\}$ and
positively homogeneous of degree $d\in\mathbb C$ with $\Re d>-n$, i.e.\ $f(tx)=t^df(x)$ for $t>0$, and let $\sum_i\int_{E_i}|f|\,dS<\infty$. Then $f\in L^1(rK)$ and
\[\int_{rK}f\,dx=\frac{r^{d+n}}{d+n}\sum_{i=1}^I\delta_i\int_{E_i}f\,dS\qquad(r>0).\]
This applies to (a) $K=Q$, $n=3$, $E_i$ the three far faces $\Sigma_i$ of \eqref{B-eq:GQ} (relative interiors), $\delta_i=1$; and (b) the L-shaped cross-section
\begin{gather*}\omega_{\rm L}:=(-1,1)^2\setminus[0,1)^2,\quad n=2,\quad E_1:=\{-1\}\times(-1,1),\quad E_2:=(-1,1)\times\{-1\},\\
E_3:=(-1,0)\times\{1\},\quad E_4:=\{1\}\times(-1,0),\end{gather*}
the four edges of $\omega_{\rm L}$ that do not pass through $0$, $\delta_i=1$.
\end{lemma}
\begin{proof}
The map $(t,y)\mapsto ty$ is a bijection from $(0,r)\times E_i$ onto $rK_i$ with Jacobian $\delta_it^{n-1}$. By Tonelli's theorem applied to $|f|$, which is homogeneous of degree $\Re d>-n$,
and by Fubini's theorem,
\begin{gather*}\int_{rK_i}|f|=\delta_i\int_{E_i}|f|\,dS\int_0^rt^{\Re d+n-1}dt<\infty,\\ \int_{rK_i}f=\delta_i\int_{E_i}f(y)\,dS(y)\int_0^rt^{d+n-1}dt;\end{gather*}
sum over $i$.
In (a), $K_i$ is the pyramid $\{x\in Q:\ x_i>x_j\ (j\ne i)\}$, and these pyramids cover $Q$ up to a null set. In (b), every point of $\omega_{\rm L}$ has $x_1<0$ or $x_2<0$, so
$\omega_{\rm L}$ is star-shaped with respect to $0$, and up to the four diagonal rays through its corners it is the disjoint union of the triangles $K_1,\dots,K_4$.
\end{proof}

\begin{lemma}[isotropic vertex pairs]\label{B-lem:vpairs}
Let $(\mathbf p,\mathbf q,\mathbf n)$ be an orthonormal frame of $\R^3$ with $\mathbf n=\mathbf p\times\mathbf q$, and put $\mathbf a:=\mathbf p+i\mathbf q\in\mathbb C^3$.
(a) $\mathbf a\cdot\mathbf a=0$, $\mathbf a\times\mathbf a=0$, $\mathbf a\cdot\mathbf n=0$ and $\mathbf a\times\mathbf n=i\mathbf a$.
(b) On the half-space $\{x:\mathbf p\cdot x>0\}$ let $\nu\in\R$ and
\begin{equation}\label{B-eq:zeta}
\zeta^\nu(x):=(\mathbf a\cdot x)^\nu
\end{equation}
(principal branch; $\Re(\mathbf a\cdot x)=\mathbf p\cdot x>0$). For every constant $w=(w_0,\vec w)\in\mathbb C\times\mathbb C^3$,
\[\nabla(\zeta^\nu w_0)-\curl(\zeta^\nu\vec w)=\nu(\mathbf a\cdot x)^{\nu-1}(w_0\mathbf a-\mathbf a\times\vec w),\qquad\Div(\zeta^\nu\vec w)=\nu(\mathbf a\cdot x)^{\nu-1}\,\mathbf a\cdot\vec w .\]
(c) The complex fields (first entry $h$, last three $g$)
\begin{equation}\label{B-eq:Z}
Z_1:=\zeta^\nu\,(0,\mathbf a),\qquad Z_2:=\zeta^\nu\,(i,\mathbf n)
\end{equation}
satisfy $\nabla h=\curl g$ and $\Div g=0$ on $\{\mathbf p\cdot x>0\}$, and so do their real and imaginary parts. Moreover $\zeta^\nu$ is harmonic there.
\end{lemma}
\begin{proof}
(a) $\mathbf a\cdot\mathbf a=|\mathbf p|^2-|\mathbf q|^2+2i\,\mathbf p\cdot\mathbf q=0$; $\mathbf a\times\mathbf a=0$ for every vector; $\mathbf a\cdot\mathbf n=0$ because $\mathbf n\perp\mathbf p,\mathbf q$; and
$\mathbf a\times\mathbf n=\mathbf p\times(\mathbf p\times\mathbf q)+i\,\mathbf q\times(\mathbf p\times\mathbf q)=-\mathbf q+i\mathbf p=i\mathbf a$.
(b) The principal power is holomorphic on $\{\Re z>0\}$, so $\nabla\zeta^\nu=\nu(\mathbf a\cdot x)^{\nu-1}\mathbf a$; for constant $\vec w$, $\curl(\zeta^\nu\vec w)=\nabla\zeta^\nu\times\vec w$ and $\Div(\zeta^\nu\vec w)=\nabla\zeta^\nu\cdot\vec w$.
(c) By (a), $w_0\mathbf a-\mathbf a\times\vec w=0$ and $\mathbf a\cdot\vec w=0$ for $w=(0,\mathbf a)$ and for $w=(i,\mathbf n)$. The operators $\nabla$, $\curl$, $\Div$ are real, so real and imaginary parts are again solutions. Finally $\Delta\zeta^\nu=\nu(\nu-1)(\mathbf a\cdot x)^{\nu-2}\,\mathbf a\cdot\mathbf a=0$.
\end{proof}

\begin{definition}[vertex-singular family]\label{B-def:family}
Fix $J\in\N$ and, for $j=1,\dots,J$, rational orthonormal frames $(\mathbf p_j,\mathbf q_j,\mathbf n_j)$ with $\mathbf n_j=\mathbf p_j\times\mathbf q_j$, such that all components of every $\mathbf p_j$ are positive. Put
\begin{equation}\label{B-eq:frames}
\mathbf a_j:=\mathbf p_j+i\mathbf q_j,\qquad p_{\min}:=\min_{j,i}(\mathbf p_j)_i>0,\qquad \zeta_j^\nu(x):=(\mathbf a_j\cdot x)^\nu ,
\end{equation}
\begin{equation}\label{B-eq:w}
w_{j,0}:=(0,\mathbf a_j),\quad w_{j,1}:=-i\,w_{j,0},\quad w_{j,2}:=(i,\mathbf n_j),\quad w_{j,3}:=-i\,w_{j,2}\qquad(\in\mathbb C^4).
\end{equation}
For a real vector $c\in\R^{4J}$ and $\nu\in[-\frac32,-1]$ let
\begin{equation}\label{B-eq:Fnu}
F_\nu=(h_\nu,g_\nu):=\sum_{j=1}^J\sum_{l=0}^3c_{4(j-1)+l+1}\,\Re\big(\zeta_j^\nu\,w_{j,l}\big)\qquad\text{on }\Gamma,
\end{equation}
and define the face sums and, if $A_h(-\frac32)>0$, the \emph{limit quotient}
\begin{equation}\label{B-eq:Ahg}
A_h(\nu):=\sum_{i=1}^3\int_{\Sigma_i}h_\nu^2\,dS,\qquad A_g(\nu):=\sum_{i=1}^3\int_{\Sigma_i}|g_\nu|^2\,dS,
\end{equation}
\begin{equation}\label{B-eq:R0}
\mathcal R_0(c):=\frac{A_h(-\frac32)-A_g(-\frac32)}{A_h(-\frac32)+A_g(-\frac32)} .
\end{equation}
\end{definition}
Thus $F_\nu=\Re\sum_j(\gamma_jZ_{1,j}+\delta_jZ_{2,j})$, with the fields \eqref{B-eq:Z} of the $j$-th frame and $\gamma_j:=c_{4j-3}-ic_{4j-2}$, $\delta_j:=c_{4j-1}-ic_{4j}$, Gaussian-rational for rational $c$.

\begin{lemma}[properties of the family]\label{B-lem:family}
Let the frames and $c$ be as in Definition~\ref{B-def:family}.
(a) For $\nu\in\R$, each $\zeta_j^\nu$ is defined and real-analytic on a neighbourhood of $\overline\Gamma\setminus\{0\}$ and positively homogeneous of degree $\nu$; for $x\in\overline\Gamma\setminus\{0\}$, $\Re(\mathbf a_j\cdot x)\ge p_{\min}(x_1+x_2+x_3)\ge p_{\min}|x|>0$, and
\begin{equation}\label{sbound}|\zeta_j^\nu(x)|=|\mathbf a_j\cdot x|^\nu\le(p_{\min}|x|)^\nu\qquad(\nu\le0).\end{equation}
(b) For $-\frac32\le\nu\le-1$ and $x\in\overline\Gamma\setminus\{0\}$,
\begin{equation}\label{mubound}|F_\nu(x)|\le\sqrt2\,\|c\|_1\,(p_{\min}|x|)^\nu\le\Theta(x):=\sqrt2\,\|c\|_1\max\{1,(p_{\min}|x|)^{-3/2}\}.\end{equation}
(c) For $\nu\in(-\frac32,-1]$ and every bounded open set $\Om'\subset\Gamma$, $F_\nu\in\M(\Om')$, and the mean-corrected pair $(h_\nu-\bar h_\nu,g_\nu)$ of Lemma~\ref{B-lem:mean} belongs to $\Mz(\Om')$.
(d) $A_h$ and $A_g$ are continuous on $[-\frac32,-1]$.
(e) For $\nu=-\frac32+\epsilon$, $0<\epsilon\le\frac12$, and $r>0$,
\begin{equation}\label{B-eq:rQ}
\int_{rQ}h_\nu^2=\frac{r^{2\epsilon}}{2\epsilon}A_h(\nu),\qquad \int_{rQ}|g_\nu|^2=\frac{r^{2\epsilon}}{2\epsilon}A_g(\nu),\qquad \int_Qh_\nu=\frac1{\nu+3}\sum_{i=1}^3\int_{\Sigma_i}h_\nu\,dS .
\end{equation}
(f) If $f$ is continuous on $\overline\Gamma\setminus\{0\}$ and positively homogeneous of degree $-3$, then $\sum_i\int_{\Sigma_i}f\,dS$ equals
$\int_{S^2\cap\Gamma}f\,d\sigma$. In particular $A_h(-\frac32)$, $A_g(-\frac32)$ and $\mathcal R_0(c)$ are integrals over the spherical triangle $S^2\cap\Gamma$; they depend only on $c$ and the frames, not on a domain.
(g) On $Q$, for $\nu=-\frac32+\epsilon$, $0<\epsilon\le\frac12$, the mean-corrected pair of (c) has the quotient
\begin{equation}\label{B-eq:Reps}
\mathcal R(h_\nu-\bar h_\nu,g_\nu)=\frac{A_h(\nu)-2\epsilon\,\bar h_\nu^2-A_g(\nu)}{A_h(\nu)-2\epsilon\,\bar h_\nu^2+A_g(\nu)},\qquad \bar h_\nu=\int_Qh_\nu=\frac1{\nu+3}\sum_{i=1}^3\int_{\Sigma_i}h_\nu\,dS ,
\end{equation}
computed from face integrals alone.
\end{lemma}
\begin{proof}
(a) The inequality holds since $x_i\ge0$, so $\{\mathbf p_j\cdot x>0\}\supset\overline\Gamma\setminus\{0\}$; the principal power is analytic on $\{\Re z>0\}$, and $(\mathbf a_j\cdot tx)^\nu=t^\nu(\mathbf a_j\cdot x)^\nu$ for $t>0$. Finally $|\mathbf a_j\cdot x|\ge\Re(\mathbf a_j\cdot x)$ and $\nu\le0$ give \eqref{sbound}.
(b) $|w_{j,l}|=\sqrt2$, so $|\Re(\zeta_j^\nu w_{j,l})|\le\sqrt2\,|\zeta_j^\nu|$, and \eqref{sbound} gives the first inequality. For $t:=p_{\min}|x|$ and $\nu\in[-\frac32,-1]$, $t^\nu\le1$ if $t\ge1$ and $t^\nu\le t^{-3/2}$ if $t<1$.
(c) By Lemma~\ref{B-lem:vpairs}(c) every term of \eqref{B-eq:Fnu} satisfies $\nabla h=\curl g$, $\Div g=0$ pointwise on $\Gamma$, hence in the sense of distributions. By (b), $|F_\nu|^2\le2\|c\|_1^2(p_{\min}|x|)^{2\nu}$ with $2\nu>-3$, which is integrable on bounded sets; so $F_\nu\in\M(\Om')$. The last statement is Lemma~\ref{B-lem:mean}.
(d) On the compact set $[-\frac32,-1]\times\Sigma_i$ the map $(\nu,y)\mapsto\zeta_j^\nu(y)=\exp(\nu\operatorname{Log}(\mathbf a_j\cdot y))$ is continuous, because $\Re(\mathbf a_j\cdot y)\ge p_{\min}>0$ there; hence the integrands of \eqref{B-eq:Ahg} are uniformly continuous, and $A_h,A_g$ are continuous.
(e) Apply Lemma~\ref{lem:face}(a) to $h_\nu^2$ and $|g_\nu|^2$ (homogeneous of degree $2\nu=-3+2\epsilon>-3$) and to $h_\nu$ (degree $\nu>-3$, $r=1$); the face integrals are finite because the integrands are continuous on the compact faces.
(f) The change of variables in the proof of Lemma~\ref{lem:face}, applied on the shell $Q\setminus e^{-1}\overline Q$, gives $\int_{Q\setminus e^{-1}\overline Q}f\,dx=\sum_i\int_{\Sigma_i}f\,dS\int_{1/e}^1t^{-1}dt=\sum_i\int_{\Sigma_i}f\,dS$. In polar coordinates, in each direction $\hat y\in S^2\cap\Gamma$ the shell is the radial interval $(\varrho_Q(\hat y)/e,\varrho_Q(\hat y))$, with $\varrho_Q$ the radial function of $Q$, and $f(s\hat y)s^2=s^{-1}f(\hat y)$; so the same integral equals $\int_{S^2\cap\Gamma}f(\hat y)\,d\sigma(\hat y)$.
(g) By Lemma~\ref{B-lem:mean} with $|Q|=1$, $\|h_\nu-\bar h_\nu\|^2=\|h_\nu\|^2-\bar h_\nu^2$; by (e) with $r=1$, $\|h_\nu\|^2=A_h(\nu)/(2\epsilon)$, $\|g_\nu\|^2=A_g(\nu)/(2\epsilon)$,
and $\bar h_\nu$ is given by the third identity of \eqref{B-eq:rQ}. Multiplying the numerator and the denominator of \eqref{B-eq:R} by $2\epsilon$ gives \eqref{B-eq:Reps}.
\end{proof}

\subsection{Upper bound for the cube}
The proof needs one pair of the family \eqref{B-eq:Fnu} with an enclosed quotient. The trial spaces are built from rational frames: $K$ \emph{axis frames}, with $\mathbf p$ close to $(1,1,1)/\sqrt3$ and $\mathbf q$ rotated about $\mathbf p$ by the
angles $\pi k/K$, and \emph{corner families}, in which $\mathbf n$ approaches a corner of $S^2\cap\Gamma$ geometrically from outside; the rounding to rational frames is
described in Section~\ref{B-sec:faceint}. Table~\ref{tab:vtx} lists the two trial spaces used below; their construction, all nine trial spaces V1--V9 and their
enclosures are in \docD{4}.
\begin{table}[tbp]
\caption{Vertex-singular trial spaces used in Theorem~\ref{thm:upper} and Lemma~\ref{lem:esscomp}: number $K$ of axis frames, frames per corner family, number $J$
of frames ($4J$ real coefficients), and the enclosed quotient (rounded down).}\label{tab:vtx}\label{tab:ess}
\small
\begin{tabular}{@{}lrlrll@{}}
\toprule
space & $K$ & corner families & $J$ & enclosed quotient & used in\\
\midrule
V6 & 9 & $3\times18$ & 63 & $\mathcal R\ge0.767204561906$ \eqref{B-eq:Reps}, $\epsilon=10^{-12}$ & Theorem~\ref{thm:upper}\\
V9 & 9 & $3\times26$ & 87 & $\mathcal R_0(c)\ge0.767204562046$ \eqref{B-eq:R0} & Lemma~\ref{lem:esscomp}\\
\bottomrule
\end{tabular}
\end{table}
\begin{theorem}\label{thm:upper}
$\kappa((-1,1)^3)\ge0.767204561906$; hence $\beta((-1,1)^3)\le0.3411711$, and $\beta^2<\frac12-\frac1\pi$.
\end{theorem}
\begin{proof}[Proof (computer-assisted)]
Work on $Q$, so that $\kappa((-1,1)^3)=\kappa(Q)$. Take the frames of V6, $\nu=-\frac32+10^{-12}$, and the rational coefficient vector $c$ of \recR{3}; it is the floating-point Ritz vector rounded to Gaussian rationals.
By Lemma~\ref{B-lem:family}(c), the mean-corrected pair $(h_\nu-\bar h_\nu,g_\nu)$ lies in $\Mz(Q)$. Its quotient \eqref{B-eq:Reps} is enclosed, by the procedure of Section~\ref{B-sec:faceint}, in
\[[0.76720456190628513928961686\pm9.2\cdot10^{-30}].\]
Hence $\kappa(Q)\ge0.767204561906$. By Theorem~\ref{thm:fv}, $\beta^2\le(1-0.767204561906)/2=0.116397719047$, so that $\beta\le0.34117109938<0.3411711$ and $\beta^2<\frac12-\frac1\pi=0.18169\ldots$.
\end{proof}
A comparison with the upper bounds inherited from planar cross-sections is given in \docD{3}.

\section{Lower bound for the cube}\label{sec:lower}
The lower bound for $\beta$ needs one design $(m,\Phi)$ and correctors for which the verifier certifies $\kappa_0$ in all four classes. The octant $[0,1]^3$ is meshed by radial layers and by cones from its vertex $(1,1,1)$ over a triangulation of the spherical triangle with alternating
diagonals and $n$ fan sectors around each edge direction; the mesh is mirrored to the cube. The multiplier $m$ is continuous, piecewise linear, equivariant under
the reflections and permutations, and zero on $\partial\Om$ (on the octant: $m=0$ on $\{x_i=1\}$, $m_i=0$ on $\{x_i=0\}$). It is computed together with the ball
weights $\alpha_k$ by a semidefinite program that maximises the smallest eigenvalue of $C_{\rm low}$ \eqref{A:Clow} at a design value of $\kappa_0$. The balls cover
the interior, with radii decreasing towards the boundary. The correctors are continuous piecewise quadratic on the conical mesh, with the parity of the class
and zero on $\partial\Om$. Numbers of functionals are listed in \docD{3}.

The \emph{limit of positive definiteness} of a design $(m,\Phi)$ is
\begin{equation}\label{B-eq:kPD}
\kappa_{\rm PD}(m,\Phi):=\inf\big\{\kappa_0\in\R:\ \kappa_0I-W(m)+\Phi\,e_0e_0^{\mathsf T}\succ0\ \text{a.e.\ in }\Om\big\};
\end{equation}
in the floating-point design computations it is evaluated with $C_{\rm low}$ in place of $C$. Where $\Phi=0$, positive definiteness requires
$\kappa_0>\lambda_{\max}(W(m))$, so near the vertices and edges $\kappa_{\rm PD}$ depends on the multiplier alone. For every design used here $\Phi$ vanishes near the vertices, so
$\kappa_{\rm PD}\ge L_*=0.767204562046$ by Corollary~\ref{B-cor:obstr}.

Table~\ref{tab:certlow} lists the certified values; the best one uses $32$ sectors and $22$ balls, and the binding class is always $110$.
The $12$-ball design passes the floating-point criterion of Appendix~\ref{app:design} only from $\kappa_0=0.76851$ on, yet it certifies $\kappa_0=0.7681$. The reason is that the verifier
(piecewise quadratic correctors, chord bound of Lemma~\ref{B-lem:Cprime}) is sharper than the floating-point evaluation ($Q_1$ correctors, cellwise weight).
Moreover $0.7681$ lies above the exact positive-definiteness floor $\kappa_{\rm PD}^{\rm ex}=0.768012$ \eqref{oth:eq:kpd} of this design.
\begin{table}[tbp]
\caption{Cube: certified $\kappa_0$ (all four classes), bounds $\lambda_{\rm up}\ge\lambda_{\max}(G)$ per class, numbers of balls and of pieces
(Section~\ref{sec:verif}), and the lower bound for $\beta$; correctors with $392\,985$ ($n=16$) and $773\,625$ ($n=32$) degrees of freedom per class; records R6 ($n=16$) and R5 ($n=32$) of \cite{lbb3drepo}. Further
values: \docD{3}.}\label{tab:certlow}
\small
\begin{tabular}{rlrrllll@{\ }l}
\toprule
$n$ & $\kappa_0$ & balls & pieces & \multicolumn{4}{c}{$\lambda_{\rm up}$: $000$, $100$, $110$, $111$} & $\beta\ge$\\
\midrule
16 & $7681/10000$ & 12 & $3\,653\,046$ & $0.5798$ & $0.7119$ & $0.8887$ & $0.6597$ & $0.340514$\\
32 & $96/125$ & 22 & $4\,484\,238$ & $0.6916$ & $0.7791$ & $0.9185$ & $0.6464$ & $0.340587$\\
\bottomrule
\end{tabular}
\end{table}
\begin{theorem}\label{thm:lower}
$\kappa((-1,1)^3)\le\frac{96}{125}$; hence $\beta((-1,1)^3)\ge\sqrt{0.116}>0.340587$.
\end{theorem}
\begin{proof}[Proof (computer-assisted)]
Apply Proposition~\ref{oth:prop:classes}(b),(c) with $\mathcal G=\mathcal G_{\rm C}$ in the classes $000$, $100$, $110$, $111$, with $\kappa_0=\frac{96}{125}$ and the data of \recR{5}: $32$ sectors, $22$ balls, piecewise quadratic correctors; the hypotheses (i)--(iv) are checked as in Remark~\ref{A:rem:verifcheck}.
The verifier of Section~\ref{sec:verif} returns $\lambda_{\rm up}=0.6916$, $0.7791$, $0.9185$, $0.6464$ for the classes $000$, $100$, $110$, $111$ (Table~\ref{tab:certlow}). These are upper bounds for $\lambda_{\max}(G_{\rm up})$ with $G_{\rm up}\succeq G^{\chi_s}$, and all are below $1$; hence $\kappa_s\le\frac{96}{125}$ for these classes, and $\kappa\le\frac{96}{125}$ by Corollary~\ref{A:lem:classes}(b).
By Theorem~\ref{thm:fv}, $\beta^2\ge\frac12\big(1-\frac{96}{125}\big)=\frac{29}{250}=0.116$, and $0.340587^2=0.115999504569<\frac{29}{250}$ exactly.
\end{proof}
\begin{corollary}[enclosure for the cube]\label{B-cor:cube}
$0.767204561906\le\kappa((-1,1)^3)\le\frac{96}{125}=0.768$ and $0.340587\le\beta((-1,1)^3)\le0.3411711$; the relative width of the enclosure of $\beta$ is $1.7\cdot10^{-3}$.
\end{corollary}
This follows from Theorems~\ref{thm:upper} and~\ref{thm:lower}; $(0.3411711-0.340587)/0.3411711=1.712\ldots\cdot10^{-3}$.
The remaining gap in $\kappa$ is analysed in \docD{3}.
\begin{remark}[published finite element values]\label{B-rem:fe}
The lowest discrete Cosserat eigenvalue on the cube in the finite element computations of \cite[Sect.~5.2, Fig.~8]{CCDL} (about $0.11$) lies below our certified bound
$\beta^2\ge0.116$; this is possible because discrete Cosserat eigenvalues of conforming Stokes pairs are not one-sided bounds for $\beta^2$. The box is discussed in
Remark~\ref{oth:rem:boxfe}.
\end{remark}
Validation runs and computing times are reported in \docD{8}.

\section{Essential spectrum generated by the vertices}\label{sec:ess}
This section shows that a convex rectangular vertex produces essential spectrum below the edge interval. Throughout, $\Om\subset\R^3$ is a bounded domain and $\mathcal S=\mathcal S_\Om$ is the Cosserat operator of Definition~\ref{A:def:S}; by Lemma~\ref{A:lem:S} it is
bounded and self-adjoint on $L^2_0(\Om)$, $0\le\mathcal S\le I$, and $\beta(\Om)^2=\min\sigma(\mathcal S)$.
\begin{definition}[essential spectrum]\label{B-def:ess}
$\sigma_{\rm ess}(\mathcal S)$ is the set of $\lambda\in\sigma(\mathcal S)$ that are not isolated eigenvalues of finite multiplicity, and
\begin{equation}\label{B-eq:lame}
\lambda_{\rm e}(\Om):=\inf\sigma_{\rm ess}(\mathcal S_\Om),\qquad \kess(\Om):=1-2\lambda_{\rm e}(\Om).
\end{equation}
\end{definition}
\begin{theorem}[{\cite[Thm.~4.1]{CCDL}}]\label{B-thm:edge}
Let $\Om$ coincide, near a point of a straight edge, with a wedge of opening $\omega\in(0,2\pi)$. Then
$[\frac12-\frac{|\sin\omega|}{2\omega},\frac12+\frac{|\sin\omega|}{2\omega}]\subset\sigma_{\rm ess}(\mathcal S_\Om)$.
\end{theorem}
The theorem is stated in \cite{CCDL} with $\sin\omega$; for re-entrant edges the absolute value is intended, as the characteristic equation carries $\pm$.
For the right angle the lower end is $\frac12-\frac1\pi=0.18169\ldots$. On the cube $\min\sigma(\mathcal S)=\beta^2$ lies below it (Proposition~\ref{prop:explicit7}),
so it is an isolated eigenvalue of finite multiplicity or a point of $\sigma_{\rm ess}(\mathcal S)$. A single quotient bounds only $\min\sigma(\mathcal S)$.
A bound for $\inf\sigma_{\rm ess}(\mathcal S)$ needs a sequence of normalized pressures $q_k\rightharpoonup0$ (Lemma~\ref{lem:weak}), because a finite-rank spectral projection maps such a sequence to zero in norm.
The vertex family supplies it: for $\nu_k\to-\frac32$ the normalized pressures concentrate at the vertex, hence converge weakly to zero, while their quotients tend to $\mathcal R_0(c)$ (Theorem~\ref{thm:ess}).

Theorem~\ref{thm:ess} proves, for convex rectangular vertices, the corner contribution conjectured in \cite[Sect.~5.2]{CCDL}, where corner conditions are formulated through the Mellin symbol on $\Re\lambda=-\frac12$ \cite[Sect.~4.3, (4.8)]{CCDL}. Its only computed ingredient is one quotient of two angular integrals (Lemma~\ref{lem:esscomp}).

\begin{lemma}[{Persson-type criterion; see e.g.\ \cite[Prop.~2.1]{HKW2013}}]\label{lem:weak}
Let $(q_k)\subset L^2_0(\Om)$ with $\|q_k\|=1$ and $q_k\rightharpoonup0$ weakly. Then $\sigma_{\rm ess}(\mathcal S)\ne\emptyset$ and
$\lambda_{\rm e}(\Om)\le\liminf_{k\to\infty}(\mathcal Sq_k,q_k)$.
\end{lemma}
Appendix~\ref{B-app:weak} gives the proof for $\mathcal S$.

\begin{lemma}[{pairwise form of the first half of the proof of \cite[Thm.~2.1]{CostabelForms}}]\label{lem:pairS}
Let $F=(h,g)\in\M$, let $\bar h$ be as in \eqref{B-eq:hbar} and $q:=h-\bar h\in L^2_0(\Om)$. Then
\[(\mathcal Sq,q)\le\frac{\|q\|^2\,\|g\|^2}{\|q\|^2+\|g\|^2}\qquad(\text{if }q\ne0\text{ or }g\ne0),\]
that is, $(\mathcal Sq,q)/\|q\|^2\le\frac12\big(1-\mathcal R(q,g)\big)$ whenever $q\ne0$.
\end{lemma}
\begin{proof}
By Lemma~\ref{B-lem:mean}, $(q,g)\in\Mz$. The first step of the proof of Theorem~\ref{thm:fv} (Cauchy--Schwarz with Lemma~\ref{A:lem:divcurl}) gives the inequality for every pair in $\Mz$.
\end{proof}

\begin{theorem}\label{thm:ess}
Let $\Om\subset\R^3$ be a bounded domain such that, after a translation and a rotation, $r_0Q\subset\Om\subset\Gamma$ for some $r_0>0$ (that is,
$\Om$ has a convex rectangular vertex at $0$ and lies in the octant of that vertex). Let the frames and $c\in\R^{4J}$ be as in Definition~\ref{B-def:family}, with $A_h(-\frac32)>0$. Then
\[\lambda_{\rm e}(\Om)\le\frac{A_g(-\frac32)}{A_h(-\frac32)+A_g(-\frac32)}=\frac{1-\mathcal R_0(c)}2,\qquad\text{i.e.}\qquad\kess(\Om)\ge\mathcal R_0(c).\]
More precisely, for every sequence $\epsilon_k\downarrow0$ with $\epsilon_k\le\frac12$ there is $k_0$ such that for $k\ge k_0$ the functions
\[q_k:=\frac{h_{\nu_k}-\bar h_{\nu_k}}{\|h_{\nu_k}-\bar h_{\nu_k}\|},\qquad\nu_k:=-\tfrac32+\epsilon_k,\]
are defined (i.e.\ $h_{\nu_k}\ne\bar h_{\nu_k}$), belong to $L^2_0(\Om)$ and satisfy $\|q_k\|=1$. As $k\to\infty$, $q_k\rightharpoonup0$, $\|q_k\|_{L^2(\Om\setminus B(0,\delta))}\to0$ for every $\delta>0$, and
$\limsup_k(\mathcal Sq_k,q_k)\le\frac12(1-\mathcal R_0(c))$.
\end{theorem}
\begin{proof}
By Lemma~\ref{lem:simil} with $r=1$, translations and rotations leave the problem invariant; we may therefore assume
$r_0Q\subset\Om\subset\Gamma$.

\emph{Step 1: admissibility.} Let $\nu\in(-\frac32,-1]$. By Lemma~\ref{B-lem:family}(c), $F_\nu\in\M(\Om)$. Moreover
$\Theta\in L^1(\Om)$, and $\Theta\in L^\infty(\Om\setminus B(0,\delta))$ for each $\delta>0$.

\emph{Step 2: asymptotics.} Let $\nu=-\frac32+\epsilon$, $0<\epsilon\le\frac12$. By \eqref{B-eq:rQ}, $\int_{r_0Q}h_\nu^2=r_0^{2\epsilon}A_h(\nu)/(2\epsilon)$ and $\int_{r_0Q}|g_\nu|^2=r_0^{2\epsilon}A_g(\nu)/(2\epsilon)$.
On $\Om\setminus r_0Q$ we have $|x|\ge\max_ix_i\ge r_0$. With
\[K_1:=|\Om|\,\|\Theta\|^2_{L^\infty(\Om\setminus B(0,r_0))},\qquad K_2:=|\Om|^{-1}\|\Theta\|_{L^1(\Om)},\qquad N_\nu:=\|h_\nu-\bar h_\nu\|,\]
which do not depend on $\nu$ except $N_\nu$, we get $0\le\int_{\Om\setminus r_0Q}|F_\nu|^2\le K_1$ and $|\bar h_\nu|\le K_2$, and $N_\nu^2=\|h_\nu\|^2-|\Om|\bar h_\nu^2$ by Lemma~\ref{B-lem:mean}.
By the continuity of $A_h,A_g$ (Lemma~\ref{B-lem:family}(d)) and $r_0^{2\epsilon}\to1$,
\begin{equation}\label{asym}2\epsilon N_\nu^2\to A_h(-\tfrac32),\qquad2\epsilon\|g_\nu\|^2\to A_g(-\tfrac32),\qquad
\frac{\|h_\nu\|^2}{N_\nu^2}\to1\qquad(\epsilon\to0+).\end{equation}
Since $A_h(-\frac32)>0$, $N_\nu\to\infty$. In particular $h_\nu-\bar h_\nu\ne0$ for small $\epsilon$, and $q_k$ is well defined for large $k$.

\emph{Step 3: concentration and weak convergence.} Let $0<\delta<r_0$. Since $|h_\nu-\bar h_\nu|\le\Theta+K_2$ and $\Theta$ is bounded on $\Om\setminus B(0,\delta)$, $\|q_k\|_{L^2(\Om\setminus B(0,\delta))}=O(N_{\nu_k}^{-1})\to0$.
For $\varphi\in L^2(\Om)$,
\[|(q_k,\varphi)|\le\frac{\|h_{\nu_k}\|}{N_{\nu_k}}\|\varphi\|_{L^2(\Om\cap B(0,\delta))}+\frac{\|\Theta\|_{L^2(\Om\setminus B(0,\delta))}\|\varphi\|+K_2\|\varphi\|_{L^1(\Om)}}{N_{\nu_k}},\]
so $\limsup_k|(q_k,\varphi)|\le\|\varphi\|_{L^2(\Om\cap B(0,\delta))}$, which tends to $0$ as $\delta\to0$. Hence $q_k\rightharpoonup0$.

\emph{Step 4: the quotient.} Lemma~\ref{lem:pairS} applied to $F_{\nu_k}\in\M$ gives
\[(\mathcal Sq_k,q_k)\le\frac{\|g_{\nu_k}\|^2}{N_{\nu_k}^2+\|g_{\nu_k}\|^2}=\frac{2\epsilon_k\|g_{\nu_k}\|^2}{2\epsilon_kN_{\nu_k}^2+2\epsilon_k\|g_{\nu_k}\|^2}
\longrightarrow\frac{A_g(-\frac32)}{A_h(-\frac32)+A_g(-\frac32)}\]
by \eqref{asym}; the denominator of the limit is positive. Lemma~\ref{lem:weak} completes the proof.
\end{proof}

Remarks on the mean correction, on other polyhedral vertices and on the Mellin symbol are in \docD{3}.

\begin{lemma}\label{lem:esscomp}
Take the $J=87$ rational frames of the trial space V9 of Table~\ref{tab:vtx}, for which $p_{\min}=1.71\cdot10^{-8}$, and the dyadic coefficient vector $c$
with denominators $2^{80}$ of \recR{7}, which also lists the frames. Then
\begin{align*}A_h(-\tfrac32)&\in[106.7229874230301992008083\pm4.01\cdot10^{-23}],\\ A_g(-\tfrac32)&\in[14.05871461089678231778197\pm4.97\cdot10^{-24}],\\ \mathcal R_0(c)&\in[0.7672045620461987022298030299\pm1.94\cdot10^{-29}] .\end{align*}
\end{lemma}
\begin{proof}[Proof (computer-assisted)]
\emph{Reduction to face matrices.} Write $h_\nu=\Re\sum_j\eta_j\zeta_j^\nu$ and $g_\nu=\Re\sum_j\zeta_j^\nu\vec v_j$, with Gaussian-rational $\eta_j,\vec v_j$ determined
exactly by $c$. Since $(\Re z)^2=\frac12\Re(z^2)+\frac12|z|^2$, $h_\nu^2=\tfrac12\Re\sum_{j,k}(\eta_j\eta_k\zeta_j\zeta_k+\eta_j\bar\eta_k\zeta_j\bar\zeta_k)$,
and similarly for $|g_\nu|^2$. Hence $A_h(-\frac32)=c^{\mathsf T}M_hc$ and $A_g(-\frac32)=c^{\mathsf T}M_gc$, with real $4J\times4J$ matrices $M_h,M_g$ that are exact rational combinations of the complex face sums
\begin{equation}\label{B-eq:Pjk}
\mathcal P_{jk}:=\sum_{i=1}^3\int_{\Sigma_i}\zeta_j^\nu\zeta_k^\nu\,dS,\qquad \tilde{\mathcal P}_{jk}:=\sum_{i=1}^3\int_{\Sigma_i}\zeta_j^\nu\overline{\zeta_k^\nu}\,dS\qquad(\nu=-\tfrac32).
\end{equation}
\emph{Enclosure.} The face sums \eqref{B-eq:Pjk} are enclosed by the procedure of Section~\ref{B-sec:faceint}, with $1767$ panels, $N_{\rm GL}=24$ nodes per direction, working precision $160$ bits and panel tolerance $10^{-22}$. The a priori quadrature error bound, rounded up and at most $1.9\cdot10^{-26}$ per entry, is added to the radius of every entry.
\emph{Assembly.} $M_h$ and $M_g$ are assembled in ball arithmetic from these enclosures with exact rational multipliers. The quadratic forms $c^{\mathsf T}M_hc$,
$c^{\mathsf T}M_gc$ and the quotient $\mathcal R_0(c)$ are evaluated in ball arithmetic with the exact dyadic $c$, and positivity of $A_h(-\frac32)$ is checked on the
enclosure. The frames and their identities (Lemma~\ref{B-lem:vpairs}(a)) are checked in rational arithmetic. The relation between the internal enclosures and the printed balls is discussed in Remark~\ref{B-rem:printed}.
\end{proof}
The vector $c$ is only proposed numerically, by a Rayleigh--Ritz computation for the pencil $(M_h-M_g,M_h+M_g)$ on the midpoints, refined by inverse
iteration at $256$ bits; it is not part of the proof. Every $c$ for which the enclosure succeeds gives a valid bound.

\begin{corollary}\label{cor:ess}
Let
\begin{equation}\label{B-eq:Lstar}
L_*:=0.767204562046 .
\end{equation}
For each of the domains
\[(-1,1)^3,\quad(0,2)\times(0,1)^2,\quad\big((-1,1)^2\setminus[0,1)^2\big)\times(0,1),\quad(-1,1)^3\setminus[0,1]^3,\]
we have $\kess(\Om)\ge L_*$, i.e.
\[\inf\sigma_{\rm ess}(\mathcal S_\Om)\le\tfrac12(1-L_*)\le0.116397718977<\tfrac12-\tfrac1\pi=0.181690\ldots\]
\end{corollary}
\begin{proof}
Each domain has a convex rectangular vertex with $r_0Q\subset\Om-b\subset\Gamma$: the cube with $b=(-1,-1,-1)$, $r_0=2$; the box with $b=0$, $r_0=1$; the
L-shaped prism with $b=(-1,-1,0)$, $r_0=1$, since $((0,2)^2\setminus[1,2)^2)\times(0,1)\supset Q$; the Fichera domain with $b=(-1,-1,-1)$, $r_0=1$, since
$(0,2)^3\setminus[1,2]^3\supset Q$.
Theorem~\ref{thm:ess} and Lemma~\ref{lem:esscomp} give $\kess\ge L_*$; the decimal comparison is exact.
\end{proof}
\begin{corollary}[spectral position]\label{B-cor:cubeess}\label{oth:prop:spec}
Let $\Om$ be one of the four domains of Corollary~\ref{cor:ess}.
\begin{enumerate}
\item[(a)] $[\frac12-\frac1\pi,\frac12+\frac1\pi]\subset\sigma_{\rm ess}(\mathcal S_\Om)$; for the L-shaped prism and the Fichera domain also
$[\frac12-\frac1{3\pi},\frac12+\frac1{3\pi}]\subset\sigma_{\rm ess}(\mathcal S_\Om)$.
\item[(b)] $\beta(\Om)^2\le\inf\sigma_{\rm ess}(\mathcal S_\Om)\le\frac12(1-L_*)\le0.116397718977$. In particular the bottom of the essential spectrum lies strictly below the edge
interval $[\frac12-\frac1\pi,\frac12+\frac1\pi]$ of Theorem~\ref{B-thm:edge}.
\item[(c)] For the cube, $0.116\le\beta((-1,1)^3)^2\le\inf\sigma_{\rm ess}(\mathcal S)\le0.116397718977$, and $\kappa((-1,1)^3)\ge\kess((-1,1)^3)\ge L_*$.
\item[(d)] For the box, the L-shaped prism and the Fichera domain, $\beta(\Om)^2$ is at most $0.071539471597$, $0.043332510496$ and $0.098773050480$,
respectively. In particular
$\beta(\Om)^2<\frac12(1-U_{\rm cert})\le\frac12(1-L)$, with $U_{\rm cert}=0.767206641725$ and the angular supremum $L$ of \eqref{B-eq:L}: $\beta^2$ lies strictly below the value
generated by the convex rectangular vertices.
\end{enumerate}
\end{corollary}
\begin{proof}
(a) At an interior point of an edge of opening $\omega$ the domain coincides locally with a wedge of opening $\omega$, which is the hypothesis of Theorem~\ref{B-thm:edge}. All four
domains have edges of opening $\frac\pi2$, and the L-shaped prism and the Fichera domain have edges of opening $\frac{3\pi}2$ (Table~\ref{oth:tab:geom}).
(b) The essential spectrum is closed and, by Lemma~\ref{lem:weak}, nonempty; hence $\inf\sigma_{\rm ess}(\mathcal S)\in\sigma_{\rm ess}(\mathcal S)\subset\sigma(\mathcal S)$, and
$\beta^2=\min\sigma(\mathcal S)\le\inf\sigma_{\rm ess}(\mathcal S)$ by Lemma~\ref{A:lem:S}. The upper bound is Corollary~\ref{cor:ess}, and $0.116397718977<\frac12-\frac1\pi$.
(c) $\beta^2\ge0.116$ is Theorem~\ref{thm:lower}. By Lemma~\ref{A:lem:S} and Theorem~\ref{thm:fv}, $\kappa=1-2\min\sigma(\mathcal S)\ge1-2\inf\sigma_{\rm ess}(\mathcal S)=\kess$,
and $\kess\ge L_*$ by Corollary~\ref{cor:ess}; $L_*$ exceeds the value $0.767204561906$ of Theorem~\ref{thm:upper} by $1.4\cdot10^{-10}$.
(d) The three numbers are $\frac12(1-\kappa)$ for the lower bounds of $\kappa$ of Proposition~\ref{B-prop:poly} (Table~\ref{tab:poly}), by Theorem~\ref{thm:fv}. By Proposition~\ref{prop:UL}, whose proof does not use this
corollary, $L\le U_*\le U_{\rm cert}$, so $\frac12(1-L)\ge\frac12(1-U_{\rm cert})=0.1163966791375>0.098773050480$.
\end{proof}
The sequences of Theorem~\ref{thm:ess} concentrate at a vertex, so the bound of Corollary~\ref{cor:ess} comes from the vertices. Whether $\beta^2$ is an isolated eigenvalue remains open (Remark~\ref{oth:rem:notproved}; for the cube, Remark~\ref{B-rem:conj}).

\section{The vertex of the cube}\label{sec:vertex}
This section compares two constants of the octant $\Gamma$ \eqref{B-eq:GQ}: the best pointwise bound $U_*$ of homogeneous multipliers and the best angular quotient $L$ of vertex-singular pairs. The comparison shows why the criterion cannot certify $\kappa_0<L_*$ (Corollary~\ref{B-cor:obstr}) and how close verified multipliers come to $L_*$.
The matrix \eqref{W} involves $m$ only through $Dm$, since $\Div m=\tr Dm$ and $\curl m$ and $\sym Dm$ are linear in $Dm$; we write $W(Dm(x))$ for the value of $W(m)$ at $x$.
\begin{definition}[vertex constants]\label{B-def:vertex}
Let $\mathcal L_1(\Gamma)$ be the set of Lipschitz fields $m$ on $\overline\Gamma$ that are positively homogeneous of degree one and vanish on $\partial\Gamma$. Put
\begin{equation}\label{B-eq:Ustar}
U(m):=\operatorname*{ess\,sup}_{x\in\Gamma}\lambda_{\max}\big(W(Dm(x))\big),\qquad U_*:=\inf_{m\in\mathcal L_1(\Gamma)}U(m),
\end{equation}
and, calling $(J,(\mathbf p_j,\mathbf q_j,\mathbf n_j)_{j\le J},c)$ \emph{admissible} if the frames and $c\in\R^{4J}$ are as in Definition~\ref{B-def:family} and $A_h(-\frac32)>0$,
\begin{equation}\label{B-eq:L}
L:=\sup\big\{\mathcal R_0(c):\ (J,(\mathbf p_j,\mathbf q_j,\mathbf n_j)_{j\le J},c)\ \text{admissible}\big\}.
\end{equation}
\end{definition}
$U(m)<\infty$ because $Dm$ is bounded. Near a vertex a certificate cannot pass below $L\le U_*$ (Corollary~\ref{B-cor:obstr}, Proposition~\ref{prop:UL}); we do not prove that $U_*$ equals $L$ or a spectral constant of the vertex.
\begin{proposition}\label{B-prop:L}
(a) $L\ge L_*$. (b) $\kess(\Om)\ge L$ for every domain $\Om$ as in Theorem~\ref{thm:ess}, in particular for the four domains of Corollary~\ref{cor:ess}.
\end{proposition}
\begin{proof}
(a) Lemma~\ref{lem:esscomp} exhibits one admissible family with $\mathcal R_0(c)\ge L_*$. (b) Theorem~\ref{thm:ess} gives $\kess(\Om)\ge\mathcal R_0(c)$ for every admissible $c$; take the supremum.
\end{proof}

Minimising $U(m)$ over continuous multipliers that are piecewise linear on a conical mesh is a semidefinite program with one $4\times4$ inequality per cone; its
floating-point optima are listed in \docD{5}.
\begin{theorem}\label{thm:vexact}
For $U_{\rm cert}=0.767254086675$, $0.767218014123$, $0.767206641725$ there is a continuous field $m$ on the octant $\Gamma$, homogeneous of degree one, piecewise linear on a conical mesh with
$49\,344$, $184\,704$, $1\,660\,032$ cones and zero on $\partial\Gamma$, such that $W(Dm(x))\prec U_{\rm cert}I$ for almost every $x\in\Gamma$. In particular $U_*\le U_{\rm cert}$.
\end{theorem}
\begin{proof}[Proof (computer-assisted)]
For $n=64,128,384$ (\recR{8}) the floating-point multipliers were rounded to rational data: rays and nodal values on the grid $2^{-60}$, with the nodal values symmetrised exactly under
the permutations of the coordinates. In exact arithmetic it was checked that (i) the rounded cones form a conforming conical mesh of $\Gamma$; (ii) $m=0$ on $\partial\Gamma$; (iii) for every cone
$\tau$, the integer matrix obtained by scaling $U_{\rm cert}I-W(Dm|_\tau)$ by a positive integer has positive leading principal minors.
By (i), the cone-wise linear interpolant $m$ of the nodal values is continuous, homogeneous of degree one and Lipschitz, with $Dm$ constant on each cone. By (iii) and Sylvester's criterion, $U_{\rm cert}I-W(Dm)\succ0$ on every cone, hence almost everywhere. Since $m\in\mathcal L_1(\Gamma)$ by (ii), $U_*\le U(m)\le U_{\rm cert}$.
\end{proof}
Theorem~\ref{thm:vexact} is a statement about $W(Dm)$ only. On the unbounded cone the field $m$ does not vanish at infinity and the pairs of interest are not in $L^2(\Gamma)$, so Lemma~\ref{lem:mult} cannot be applied on $\Gamma$ directly. The following local lemma replaces it.
\begin{lemma}[local multiplier bound]\label{B-lem:local}
Let $\rho>0$, let $m$ be Lipschitz on $\rho\overline Q$ with $m=0$ on the three faces of $\rho Q$ through $0$, and let $W(Dm(x))\preceq U\,I$ for almost every $x\in\frac\rho2Q$. Then $U\ge\mathcal R_0(c)$ for all frames and $c$ as in Definition~\ref{B-def:family} with $A_h(-\frac32)>0$. Consequently $U\ge L\ge L_*$.
\end{lemma}
\begin{proof}
Let $\eta(x):=\min\{1,2(1-\max_ix_i/\rho)\}_+$ on $\rho\overline Q$ and $\tilde m:=\eta m$. Then $\tilde m$ is Lipschitz on $\rho\overline Q$ and vanishes on $\partial(\rho Q)$ (on the faces through $0$ because
$m=0$, on the others because $\eta=0$). Let $\nu=-\frac32+\epsilon\in(-\frac32,-1]$. By Lemma~\ref{B-lem:family}(c), $F_\nu\in\M(\rho Q)$, so Lemma~\ref{lem:mult} on $\rho Q$ gives
$\int_{\rho Q}(h_\nu^2-|g_\nu|^2)=\int_{\rho Q}F_\nu^{\mathsf T}W(D\tilde m)F_\nu$.

On $\frac\rho2Q$ we have $\eta=1$ and $W(D\tilde m)=W(Dm)\preceq UI$. On $\rho Q\setminus\frac\rho2Q$, $D\tilde m=\eta Dm+m\otimes\nabla\eta$ is
bounded, and $|x|\ge\max_ix_i\ge\frac\rho2$, so $|F_\nu|\le\Theta$ is bounded there uniformly in $\nu$ by \eqref{mubound}. Hence there is $K_3$, independent of $\nu$, with
\[\int_{\rho Q}(h_\nu^2-|g_\nu|^2)\le U\int_{\frac\rho2Q}|F_\nu|^2+K_3,\qquad \int_{\rho Q}(h_\nu^2-|g_\nu|^2)\ge\int_{\frac\rho2Q}(h_\nu^2-|g_\nu|^2)-K_3 .\]
By \eqref{B-eq:rQ} with $r=\frac\rho2$, multiplying by $2\epsilon(\rho/2)^{-2\epsilon}$,
$A_h(\nu)-A_g(\nu)\le U\big(A_h(\nu)+A_g(\nu)\big)+4\epsilon(\rho/2)^{-2\epsilon}K_3$. Letting $\epsilon\to0$, the continuity of $A_h,A_g$ (Lemma~\ref{B-lem:family}(d)) gives
$A_h(-\frac32)-A_g(-\frac32)\le U\big(A_h(-\frac32)+A_g(-\frac32)\big)$, i.e.\ $U\ge\mathcal R_0(c)$. Taking the supremum over $c$ and using Proposition~\ref{B-prop:L}(a) gives $U\ge L\ge L_*$.
\end{proof}
\begin{proposition}\label{prop:UL}
Let $m$ be Lipschitz on $\overline\Gamma$, positively homogeneous of degree one, $m=0$ on $\partial\Gamma$, and $W(Dm(x))\preceq U\,I$ for almost every $x\in\Gamma$. Then $U\ge\mathcal R_0(c)$ for every family of
frames and every $c$ with $A_h(-\frac32)>0$ as in Definition~\ref{B-def:family}. Consequently $U_*\ge L\ge L_*=0.767204562046$, and $L_*\le L\le U_*\le U_{\rm cert}=0.767206641725$ with the finest certificate of Theorem~\ref{thm:vexact}.
\end{proposition}
\begin{proof}
Apply Lemma~\ref{B-lem:local} with $\rho=1$ to the restriction of $m$ to $\overline Q$; the faces of $Q$ through $0$ lie in $\partial\Gamma$. This gives $U\ge\mathcal R_0(c)$ and $U\ge L$. Every $m\in\mathcal L_1(\Gamma)$ satisfies the hypothesis with $U=U(m)$, so $U_*\ge L$; $L\ge L_*$ is Proposition~\ref{B-prop:L}(a), and $U_*\le U_{\rm cert}$ is Theorem~\ref{thm:vexact}.
\end{proof}
\begin{corollary}[vertex obstruction]\label{B-cor:obstr}
Let $\Om$ be as in Theorem~\ref{thm:ess} (after a rigid motion, $r_0Q\subset\Om\subset\Gamma$).
(a) Let $m$ be Lipschitz on $\overline\Om$ with $m=0$ on $\partial\Om$ and $W(Dm)\preceq U'I$ almost everywhere in $\Om\cap B(0,\delta)$ for some $\delta>0$. Then $U'\ge L\ge L_*$.
(b) Let $m$ be as in (a), $\Phi\ge0$ with $\Phi=0$ on $\Om\cap B(0,\delta)$, and $\kappa_0I-W(m)+\Phi e_0e_0^{\mathsf T}\succ0$ almost everywhere in $\Om\cap B(0,\delta)$. Then $\kappa_0\ge L_*$.
In particular $\kappa_{\rm PD}(m,\Phi)\ge L_*$ in \eqref{B-eq:kPD} whenever $\Phi$ vanishes near the vertex, as the function $\Phi$ of Section~\ref{sec:cert} does. Moreover, if the data of Theorem~\ref{thm:rel} on $\Om$ satisfy $0\prec C'\preceq C-\rho_{\rm up}e_0e_0^{\mathsf T}$ almost everywhere, then $\kappa_0\ge L_*$: the criterion cannot certify any $\kappa_0<L_*$.
\end{corollary}
\begin{proof}
(a) Let $\rho:=\min\{r_0,\delta/2\}$. Then $\rho Q\subset\Om$ and $\frac\rho2Q\subset B(0,\frac{\sqrt3}2\rho)\subset B(0,\delta)$. The faces of $\rho Q$ through $0$ lie in $\partial\Gamma\cap\overline\Om\subset\partial\Om$, because $\Om\subset\Gamma$ is open; so $m=0$ there. Lemma~\ref{B-lem:local} applied to $m|_{\rho\overline Q}$ gives $U'\ge L\ge L_*$.
(b) Where $\Phi=0$, the hypothesis says $W(m)\prec\kappa_0I$, so (a) applies with $U'=\kappa_0$. The function $\Phi$ of Section~\ref{sec:cert} vanishes near every vertex, because it is a finite sum of bumps whose support balls have closures in $\Om$. Finally $\rho_{\rm up}\ge0$, so $0\prec C'\preceq C-\rho_{\rm up}e_0e_0^{\mathsf T}$ implies $C\succ0$ almost everywhere, and the previous statement applies.
\end{proof}
Corollary~\ref{B-cor:obstr} explains the limits of positive definiteness in Table~\ref{tab:lower} and in \docD{5}, all of which exceed $L_*$, and the refusals of the verifier reported in Appendix~\ref{app:design}.

The gap $U_{\rm cert}-L_*$ is $4.95\cdot10^{-5}$, $1.35\cdot10^{-5}$ and $2.08\cdot10^{-6}$ for $n=64,128,384$; observed to decrease at about second order in the number $n$ of fan sectors per edge
(\docD{5}). Each value compares one verified multiplier bound with one enclosed angular trial value. By Proposition~\ref{prop:UL} it bounds $U_*-L$, but it does not show $U_*=L$.

\begin{remark}[conjecture]\label{B-rem:conj}
Proved: $\kappa((-1,1)^3)\ge\kess\ge L$ (Corollary~\ref{B-cor:cubeess}(c), Proposition~\ref{B-prop:L}(b)), and near a vertex the criterion of Theorem~\ref{thm:rel} cannot pass below $L_*$ (Corollary~\ref{B-cor:obstr}).
Computed: $U_{\rm cert}-L_*\le2.1\cdot10^{-6}$. Observed: two discretisations of the angular problem without error analysis give $L\approx0.7672045620$
(\docD{5}). This suggests
$\kappa((-1,1)^3)=L\approx0.7672046$, i.e.\ $\beta\approx0.341171$, and $\beta^2=\inf\sigma_{\rm ess}(\mathcal S)$. Two questions remain open: whether $U_*=L$, and whether the cube has an isolated eigenvalue of $\mathcal S$ below $\inf\sigma_{\rm ess}(\mathcal S)$. The local constants of the other tangent cones are in Table~\ref{oth:tab:local}.
\end{remark}

\setlength{\emergencystretch}{1.5em}%
\section{Other polyhedra}\label{sec:other}
This section proves two-sided bounds for $\beta$ on three further rectangular polyhedra, i.e.\ polyhedra with axis-parallel faces,
\begin{equation}\label{oth:eq:domains}
\Om_{\rm B}:=(0,2)\times(0,1)^2,\qquad \Om_{\rm L}:=\omega_{\rm L}\times(0,1),\qquad \Om_{\rm F}:=(-1,1)^3\setminus[0,1]^3,
\end{equation}
with $\omega_{\rm L}=(-1,1)^2\setminus[0,1)^2$ as in Lemma~\ref{lem:face}(b): the \emph{box}, the \emph{L-shaped prism} and the \emph{Fichera domain}. For each domain we need
a reduction to finitely many symmetry classes, one trial pair (Section~\ref{sec:upper}), and $\kappa_\chi\le\kappa_0$ in each class needed (Sections~\ref{sec:cert},
\ref{sec:verif}). The method uses the cube only through its symmetry group and the exact description of one octant; these are replaced. By Lemma~\ref{lem:simil} we may
use translated copies; the coordinates are fixed at the
beginning of each subsection. Section~\ref{oth:ssec:cert} proves the three lower bounds by one argument (Proposition~\ref{oth:prop:cert}); Table~\ref{oth:tab:summary}
collects the results.

\subsection{Common tools}\label{oth:ssec:common}
Theorem~\ref{thm:fv} holds on every bounded domain. Lemma~\ref{lem:mult} holds for every multiplier that is Lipschitz on
$\overline\Om$ and vanishes on $\partial\Om$. Lemma~\ref{lem:ball} and the representers \eqref{rep3} only need $B(x_k,R_k)\subset\Om$. The proof of Theorem~\ref{thm:rel} is
domain independent, and Proposition~\ref{oth:prop:classes} holds for every group of reflections. The domain-dependent checks are collected in
Proposition~\ref{oth:prop:cert}. The new point is the Lipschitz bound on a non-convex domain, where a segment between points of $\overline\Om$ may leave $\overline\Om$.
\begin{lemma}[Lipschitz multipliers on non-convex domains]\label{oth:lem:lip}
Let $\Om\subset\R^3$ be a bounded open set, let $\mathfrak T$ be a finite set of closed tetrahedra with $\bigcup\mathfrak T=\overline\Om$, and let $m:\overline\Om\to\R^3$ be continuous,
affine on each $\tau\in\mathfrak T$, and $m=0$ on $\partial\Om$. Put $K_m:=\max_{\tau\in\mathfrak T}|Dm|_\tau|$ (operator norm). Then $|m(x)-m(y)|\le K_m|x-y|$ for all
$x,y\in\overline\Om$; in particular $m$ satisfies the hypothesis of Lemma~\ref{lem:mult}, also when $\Om$ is not convex.
\end{lemma}
The elementary proof is given in Appendix~\ref{A:app:lip}.

\emph{Upper bounds.} Each proof needs one pair in $\Mz$, exactly admissible, and an exact value or an enclosure of its quotient \eqref{kappa}. The polynomial parts
are exact Cauchy--Kovalevskaya extensions of rational traces (Lemma~\ref{B-lem:poly}(d)); both equations are verified as polynomial identities in exact arithmetic. The traces
are proposed by floating-point Rayleigh--Ritz computations and rounded; only the exact pair enters a proof.

\emph{Floors of a design.} Two numbers explain the limits of a design (Remark~\ref{oth:rem:limit}); they enter no proof. For fixed design
data and fixed refinement parameters of the verifier, let $B_P(\kappa_0)$ be the exact piecewise matrices \eqref{B-eq:BP} on the pieces $P$. Let $\mathfrak T_\partial$ be
the set of tetrahedra $\tau$ of the multiplier mesh on which the cellwise lower bound of $\Phi-\rho_{\rm up}$, computed in floating point by the design code, is not
positive. We call
\begin{equation}\label{oth:eq:kpd}
\kappa_{\rm PD}^{\rm ex}:=\inf\{\kappa_0:\ B_P(\kappa_0)\succ0\ \text{for every piece }P\},\qquad \kappa_\partial:=\max_{\tau\in\mathfrak T_\partial}\lambda_{\max}(W_\tau)
\end{equation}
the \emph{exact positive-definiteness floor} and the \emph{boundary floor} of the design. For $\kappa_0\le\kappa_{\rm PD}^{\rm ex}$ the verifier refuses before any Gram integral is
formed, whatever the correctors. At fixed pieces the smallest eigenvalue of $B_P(\kappa_0)$ moves exactly with $\kappa_0$. The adaptive pieces depend on $\kappa_0$, so we
bracket $\kappa_{\rm PD}^{\rm ex}$ by builds at several values of $\kappa_0$.
\begin{lemma}[floors]\label{oth:lem:floors}
(a) $\kappa_{\rm PD}(m,\Phi)\le\kappa_{\rm PD}^{\rm ex}$, with $\kappa_{\rm PD}(m,\Phi)$ from \eqref{B-eq:kPD}.
(b) If a piece $P$ of a tetrahedron $\tau$ has $\phi_P^{\min}\le0$ in Definition~\ref{B-def:BP}, then $B_P(\kappa_0)\succ0$ implies $\kappa_0>\lambda_{\max}(W_\tau)$.
(c) Let $\Om$ be as in Theorem~\ref{thm:ess}, after a rigid motion, with its vertex at $0$. If $W(m)\preceq\kappa_\partial I$ almost everywhere in $\Om\cap B(0,\delta)$ for some
$\delta>0$, then $\kappa_\partial\ge L_*$.
\end{lemma}
\begin{proof}
(a) If $B_P(\kappa_0)\succ0$ for every piece, then $0\prec C'\preceq C-\rho_{\rm up}e_0e_0^{\mathsf T}\preceq C$ almost everywhere by Lemma~\ref{B-lem:Cprime} and
$\rho_{\rm up}\ge0$, i.e.\ $\kappa_0I-W(m)+\Phi e_0e_0^{\mathsf T}\succ0$ almost everywhere; so $\kappa_{\rm PD}(m,\Phi)\le\kappa_0$.
(b) $B_P=\kappa_0'I-W_{{\rm up},P}+\phi_P^{\min}e_0e_0^{\mathsf T}\preceq\kappa_0I-W_\tau$, because $\kappa_0'\le\kappa_0$, $W_{{\rm up},P}\succeq W_\tau$ and $\phi_P^{\min}\le0$.
(c) Corollary~\ref{B-cor:obstr}(a) with $U'=\kappa_\partial$.
\end{proof}
Up to the difference between the cellwise bound of the design code and $\phi_P^{\min}$, (b) gives $\kappa_\partial\le\kappa_{\rm PD}^{\rm ex}$. Since $\Phi$ vanishes near every
vertex and $\rho_{\rm up}\ge0$, the tetrahedra near a vertex belong to $\mathfrak T_\partial$. Hence (c) applies, and $\kappa_\partial\ge L_*$ for every design with a convex
rectangular vertex; for the cube, $\kappa_{\rm PD}(m,\Phi)\ge L_*$ is Corollary~\ref{B-cor:obstr}.

\emph{Local geometry.} Table~\ref{oth:tab:geom} lists the faces, edges and vertices of the four domains by type, and Table~\ref{oth:tab:local} the local constants of the corresponding tangent cones
on the scale of $\kappa$. An edge interval $[\frac12-\frac{|\sin\omega|}{2\omega},\frac12+\frac{|\sin\omega|}{2\omega}]$ of the Cosserat spectrum corresponds to the value
$|\sin\omega|/\omega=1-2\big(\frac12-\frac{|\sin\omega|}{2\omega}\big)$. The \emph{L-cone} is the tangent cone at an end point of a re-entrant edge, a sector of opening
$\frac{3\pi}2$ times a half-line; the \emph{Fichera cone} is $\R^3$ minus a closed octant.
\begin{table}[t]
\caption{Faces, edges by opening, vertices by tangent cone, and order of the symmetry group.}\label{oth:tab:geom}
\small
\setlength{\tabcolsep}{4.5pt}
\begin{tabular}{@{}lrrrrrrr@{}}
\toprule
 & & \multicolumn{2}{c}{edges} & \multicolumn{3}{c}{vertices} & \\
\cmidrule(lr){3-4}\cmidrule(lr){5-7}
domain & faces & $\frac\pi2$ & $\frac{3\pi}2$ & octant & L-cone & Fichera cone & $|\text{group}|$\\
\midrule
cube $(-1,1)^3$ & 6 & 12 & 0 & 8 & 0 & 0 & 48\\
box $\Om_{\rm B}$ & 6 & 12 & 0 & 8 & 0 & 0 & 16\\
L-prism $\Om_{\rm L}$ & 8 & 17 & 1 & 10 & 2 & 0 & 4\\
Fichera $\Om_{\rm F}$ & 9 & 18 & 3 & 10 & 3 & 1 & 6\\
\bottomrule
\end{tabular}
\end{table}
\begin{table}[t]
\caption{Local constants of the tangent cones on the scale of $\kappa$. Edges: value $|\sin\omega|/\omega$ of the edge interval of \cite[Thm.~4.1]{CCDL}. Vertices: pointwise multiplier
bound $U$ of optimised conical multipliers (Section~\ref{sec:vertex}). ``fp'': floating-point semidefinite programs, not certified.}\label{oth:tab:local}
\small
\begin{tabular}{@{}lp{0.30\textwidth}p{0.42\textwidth}@{}}
\toprule
tangent cone & value & status\\
\midrule
half-space & $0$ & $\lambda_{\max}(W)=0$ for $m=2\dist(\cdot,\text{face})\,n$ (Proposition~\ref{prop:local}(a))\\
wedge $\frac\pi2$ & $2/\pi=0.636620\ldots$ & edge interval proved \cite[Thm.~4.1]{CCDL}\\
wedge $\frac{3\pi}2$ & $2/(3\pi)=0.212207\ldots$ & edge interval proved \cite[Thm.~4.1]{CCDL}\\
octant & $[0.767204562046,$ $0.767206641725]$ & proved (Theorem~\ref{thm:vexact}, Proposition~\ref{prop:UL})\\
L-cone & $U=0.63668370$ ($n=64$) & fp; equals the $n=64$ fan value $1/(n\sin\frac\pi{2n})=0.63668369$ of a $\frac\pi2$ edge to $4.9\cdot10^{-9}$ (limit $2/\pi$)\\
Fichera cone & $U\in[0.2123537,$ $0.2123676]$ ($n\ge64$) & fp; $1.5$--$1.6\cdot10^{-4}$ above $2/(3\pi)$\\
\bottomrule
\end{tabular}
\end{table}
\begin{remark}[what is not proved]\label{oth:rem:notproved}
Corollary~\ref{B-cor:cubeess} leaves two alternatives for each of the three domains: $\beta^2$ is an isolated eigenvalue of $\mathcal S$ of finite multiplicity, or
$\beta^2=\inf\sigma_{\rm ess}(\mathcal S)$. The second would require essential spectrum at or below $0.0716$ (box), $0.0434$ (L-prism) or $0.0988$ (Fichera domain), strictly
below the vertex value $\frac12(1-L)$ and below the edge intervals of Corollary~\ref{B-cor:cubeess}(a). Nothing in this paper produces spectrum there, but we do not exclude
it. Excluding it needs a \emph{lower} bound for $\inf\sigma_{\rm ess}(\mathcal S)$, for instance a localisation principle through the local constants of
Table~\ref{oth:tab:local}, which we do not prove; whether eigenvalues of the edge symbol add further essential spectrum is open \cite[Rem.~4.2]{CCDL}. The certified bounds
give $\inf\sigma_{\rm ess}(\mathcal S)\ge\beta^2\ge\frac{13}{200}$, $\frac3{80}$ and $\FichBetaSqLow$ (Corollaries~\ref{oth:cor:box}, \ref{oth:cor:lp} and~\ref{oth:cor:fich}).
The fast convergence of the trial values in the binding classes is consistent with an isolated eigenvalue but is not a proof. For the cube see Remark~\ref{B-rem:conj}.
\end{remark}

\subsection{The box}\label{oth:ssec:box}
The box needs no new tool: its group is generated by reflections and one swap, so Proposition~\ref{oth:prop:classes} applies directly.
\emph{Coordinates and geometry.} We work with the translate $\Om_{\rm B}-(1,\frac12,\frac12)=(-1,1)\times(-\frac12,\frac12)^2$, with $x=x_1$ the long direction. Its
symmetry group has order $16$. It is generated by the eight reflections $J_\varsigma:=\operatorname{diag}(\varsigma_1,\varsigma_2,\varsigma_3)$, $\varsigma_i=\pm1$, and the swap
$\Pi_{23}(x_1,x_2,x_3):=(x_1,x_3,x_2)$. The reflections form a group $\mathcal G_{\rm B}$ as in Definition~\ref{A:def:classes}, with fundamental region the octant
$\Om_{\mathcal G}=[0,1]\times[0,\frac12]^2$. The octant has three mirror faces $\{x_i=0\}$ and three outer faces $\{x_1=1\}$, $\{x_2=\frac12\}$, $\{x_3=\frac12\}$. It contains
one vertex $(1,\frac12,\frac12)$ of the box and three outer edges, all of opening $\frac\pi2$: $\{x_1=1,x_2=\frac12\}$ and $\{x_1=1,x_3=\frac12\}$ of length $\frac12$ in the
octant, and $\{x_2=x_3=\frac12\}$ of length $1$. The edges of the octant on a mirror plane are not edges of the box.
\emph{Symmetry classes.} We label the characters of $\mathcal G_{\rm B}$ as for the cube, $\chi_s(J_\varsigma)=\prod_i\varsigma_i^{s_i}$ with $s\in\{0,1\}^3$, and write
$\kappa_s:=\kappa_{\chi_s}$.
\begin{corollary}[classes of the box]\label{oth:cor:boxclasses}
\begin{equation}\label{oth:eq:boxclasses}
\kappa(\Om_{\rm B})=\max_{s\in\mathcal S_{\rm B}}\kappa_s,\qquad \mathcal S_{\rm B}:=\{000,\,100,\,010,\,110,\,011,\,111\}.
\end{equation}
In class $s$, $h(J_ix)=-(-1)^{s_i}h(x)$ for the reflection $J_i$ in $x_i$; the constants belong to class $111$ only.
\end{corollary}
\begin{proof}
Proposition~\ref{oth:prop:classes}(a) gives $\kappa(\Om_{\rm B})=\max_{s\in\{0,1\}^3}\kappa_s$. The swap $\Pi_{23}$ maps the centred box onto itself, because its second and third
side lengths agree, and $\Pi_{23}J_\varsigma\Pi_{23}^{-1}=\operatorname{diag}(\varsigma_1,\varsigma_3,\varsigma_2)\in\mathcal G_{\rm B}$; hence $\chi_s^{\Pi_{23}}=\chi_{(s_1,s_3,s_2)}$, and
Lemma~\ref{oth:lem:normal} gives $\kappa_{001}=\kappa_{010}$ and $\kappa_{101}=\kappa_{110}$. The parity follows from $\mathcal T_{J_i}F=\chi_s(J_i)F$ and $\det J_i=-1$. A constant
$h$ is even in all variables, i.e.\ $s=111$.
\end{proof}
No further reduction is available: a permutation involving $x_1$ does not map the box onto itself.

\emph{Upper bound.} Table~\ref{tab:poly} lists the box with polynomial pairs of degree $64$. The largest values occur in class $011$, in which $h$ is odd in the long
direction and even in $x_2,x_3$.
\begin{proposition}[upper bound for the box]\label{oth:prop:boxup}
$\kappa(\Om_{\rm B})\ge\kappa_{011}\ge0.856921056806$; hence $\beta(\Om_{\rm B})^2\le0.071539471597$ and $\beta(\Om_{\rm B})\le0.267469$.
\end{proposition}
\begin{proof}[Proof (computer-assisted)]
The trial pair (\recR{2}) is a polynomial pair of total degree $64$ in class $011$ with rational coefficients, constructed and checked as in Section~\ref{oth:ssec:common}
(extension in $x_3$ of a trace on $\{x_3=0\}$); $\int h=0$ because $h$ is odd in $x_1$. Its quotient \eqref{kappa} is a rational number, computed by exact integration of polynomials over the box; rounded down to twelve
decimals it is $0.856921056806$. Theorem~\ref{thm:fv} gives $\beta^2\le\frac12(1-0.856921056806)=0.071539471597$ and $\beta\le0.2674686367$.
\end{proof}
The trace is proposed on the space of polynomial pairs of degree at most $64$ in class $011$ (dimension $1056$). The exact bounds for lower degrees and the floating-point
Ritz values of the other classes are listed in \docD{6}.

\begin{remark}[published finite element values for the box]\label{oth:rem:boxfe}
For the box, the plotted lowest discrete Cosserat eigenvalue of \cite[Sect.~5.2, Fig.~8]{CCDL} (about $0.07$) is consistent with Proposition~\ref{oth:prop:boxup}, and the
quasimode bound \cite[Lemma~5.1, (5.7)]{CCDL} gives $\beta^2\le0.0775$; Proposition~\ref{oth:prop:boxup} and Theorem~\ref{oth:thm:box} give $0.065\le\beta^2\le0.071539471597$.
\end{remark}

\subsection{The L-shaped prism}\label{oth:ssec:lp}
On this non-convex domain polynomial pairs converge slowly, so we add pairs singular at the re-entrant edge.
\emph{Coordinates and geometry.} We work with the translate $\Om_{\rm L}-(0,0,\frac12)=\omega_{\rm L}\times(-\frac12,\frac12)$. The cross-section $\omega_{\rm L}$ has five
convex corners and the re-entrant corner $0$. Accordingly $\Om_{\rm L}$ has the re-entrant edge $\{x_1=x_2=0\}$ of opening $\frac{3\pi}2$, $17$ edges of opening $\frac\pi2$,
two L-cone vertices $(0,0,\pm\frac12)$ and ten octant vertices (Table~\ref{oth:tab:geom}). The symmetry group is $\mathcal G_{\rm L}:=\{I,Z,\Pi_{12},\Pi_{12}Z\}$ with
\begin{equation}\label{oth:eq:GL}
Z(x_1,x_2,x_3):=(x_1,x_2,-x_3),\qquad \Pi_{12}(x_1,x_2,x_3):=(x_2,x_1,x_3),
\end{equation}
a group of reflections as in Definition~\ref{A:def:classes}; there is no reflection in $x_1$ or in $x_2$. A fundamental region is
\begin{equation}\label{oth:eq:FL}
\Om_{\mathcal G}:=\{x\in\overline{\omega_{\rm L}}\times[0,\tfrac12]:\ x_2\ge x_1\}=\big([-1,0]\times[0,1]\ \cup\ \{-1\le x_1\le x_2\le0\}\big)\times[0,\tfrac12],
\end{equation}
with volume $|\Om_{\mathcal G}|=\frac34$. Its boundary consists of pieces of $\partial\Om$ (the faces $x_1=-1$, $x_2=1$, the notch face $x_1=0$, $0\le x_2\le1$, and the top $x_3=\frac12$) and of the mirror
pieces $\{x_3=0\}$ and $\{x_1=x_2\}$. The re-entrant edge lies in the mirror plane $\{x_1=x_2\}$.

\emph{Symmetry classes.} A character of $\mathcal G_{\rm L}$ is determined by $\chi(Z)=(-1)^{s_3}$ and $\chi(\Pi_{12})=(-1)^e$, $s_3,e\in\{0,1\}$. We label it by the string
$00s_3e$. The two leading zeros record that there is no reflection in $x_1$, $x_2$; the labels agree with those used in \recR{10}.
\begin{corollary}[classes of the L-shaped prism]\label{oth:cor:lpclasses}
\begin{equation}\label{oth:eq:lpclasses}
\kappa(\Om_{\rm L})=\max_{\chi\in\mathcal S_{\rm L}}\kappa_\chi,\qquad \mathcal S_{\rm L}:=\{0000,\,0001,\,0010,\,0011\}.
\end{equation}
In class $00s_3e$, $h(Zx)=-(-1)^{s_3}h(x)$ and $h(\Pi_{12}x)=-(-1)^eh(x)$: $h$ is odd in $x_3$ if $s_3=0$ and even if $s_3=1$, odd under $x_1\leftrightarrow x_2$ if $e=0$ and even if
$e=1$. The constants belong to class $0011$ only.
\end{corollary}
\begin{proof}
Proposition~\ref{oth:prop:classes}(a) with $\mathcal G=\mathcal G_{\rm L}$; the parities follow from $\mathcal T_RF=\chi(R)F$ and $\det Z=\det\Pi_{12}=-1$.
\end{proof}
The trial computations below use only the reflection $Z$: their class ``$h$ even in $x_3$'' is $\Mz^{0010}\oplus\Mz^{0011}$, and ``$h$ odd in $x_3$'' is
$\Mz^{0000}\oplus\Mz^{0001}$.

\emph{Upper bound.} Polynomial pairs converge slowly on $\Om_{\rm L}$; their bounds are listed in \docD{6}. We add pairs that are singular at the re-entrant edge, generated
from harmonic potentials.
\begin{lemma}[pairs from harmonic potentials]\label{oth:lem:pot}
Let $\mathcal U\subset\R^3$ be open and $u$ a complex-valued harmonic function on $\mathcal U$. Then
\begin{equation}\label{oth:eq:Pj}
\mathcal Q_0(u):=(0,\nabla u),\qquad \mathcal Q_j(u):=(\partial_ju,\ \nabla u\times e_j),\quad j=1,2,3,
\end{equation}
satisfy $\nabla h=\curl g$ and $\Div g=0$ in $\mathcal U$, and so do their real and imaginary parts.
\end{lemma}
\begin{proof}
For $\mathcal Q_0(u)$: $\curl\nabla u=0=\nabla0$ and $\Div\nabla u=\Delta u=0$. For $\mathcal Q_j(u)$, the identities $\curl(A\times B)=A\Div B-B\Div A+(B\cdot\nabla)A-(A\cdot\nabla)B$ and
$\Div(A\times B)=B\cdot\curl A-A\cdot\curl B$ with $A=\nabla u$ and the constant $B=e_j$ give $\curl(\nabla u\times e_j)=-e_j\Delta u+\partial_j\nabla u=\nabla(\partial_ju)$
and $\Div(\nabla u\times e_j)=e_j\cdot\curl\nabla u=0$. The operators are real, so real and imaginary parts are admissible.
\end{proof}
Let $\zeta_{\rm L}:=\varpi\,(x_1+ix_2)$ with $\varpi:=e^{-5\pi i/4}$; then $\omega_{\rm L}=\{x:-\frac{3\pi}4<\arg\zeta_{\rm L}<\frac{3\pi}4\}\cap(-1,1)^2$, and the principal branch of $\zeta_{\rm L}^\nu$ is analytic
except on the ray $\arg(x_1+ix_2)=\frac\pi4$, which lies in the removed quadrant. For $\lambda\in\mathbb C$ with $\Re\lambda>0$ and integers $p,k\ge0$ put
\begin{equation}\label{oth:eq:upot}
\begin{gathered}
u_{\lambda,p,k}:=\sum_{n=0}^{\lfloor k/2\rfloor}c_n\,\zeta_{\rm L}^{\lambda+p+n}\,\bar\zeta_{\rm L}^{\,n}\,x_3^{k-2n},\\
c_0:=1,\qquad c_{n+1}:=-\frac{(k-2n)(k-2n-1)}{4(\lambda+p+n+1)(n+1)}\,c_n .
\end{gathered}
\end{equation}
\begin{lemma}\label{oth:lem:upot}
$u_{\lambda,p,k}$ is harmonic on $\R^3$ minus the half-plane $\{\arg(x_1+ix_2)=\frac\pi4\}$ and the $x_3$-axis; in particular on a neighbourhood of
$\overline{\Om_{\rm L}}$ minus the re-entrant edge. The components of the pairs \eqref{oth:eq:Pj} built from $u_{\lambda,p,k}$ are of the form $\Re(\zeta_{\rm L}^{\lambda-1}\Upsilon)$ or
$\Im(\zeta_{\rm L}^{\lambda-1}\Upsilon)$ with polynomials $\Upsilon$ in $x_1,x_2,x_3,\bar\zeta_{\rm L}$ and $\zeta_{\rm L}$, and they belong to $L^2(\Om_{\rm L})$.
\end{lemma}
\begin{proof}
Write $\zeta:=\zeta_{\rm L}$ in this proof. The map $x_1+ix_2\mapsto\zeta$ is a rotation, so $\Delta=4\partial_\zeta\partial_{\bar\zeta}+\partial_3^2$ in the variables $(\zeta,\bar\zeta,x_3)$. For the term
$t_n:=c_n\zeta^{\lambda+p+n}\bar\zeta^{\,n}x_3^{k-2n}$, $4\partial_\zeta\partial_{\bar\zeta}t_n=4c_n(\lambda+p+n)n\,\zeta^{\lambda+p+n-1}\bar\zeta^{\,n-1}x_3^{k-2n}$ and
$\partial_3^2t_n=c_n(k-2n)(k-2n-1)\zeta^{\lambda+p+n}\bar\zeta^{\,n}x_3^{k-2n-2}$. In $\Delta\sum_nt_n$ the coefficient of $\zeta^{\lambda+p+n}\bar\zeta^{\,n}x_3^{k-2n-2}$ is
$4c_{n+1}(\lambda+p+n+1)(n+1)+c_n(k-2n)(k-2n-1)$, which vanishes by the recursion; the first term of $4\partial_\zeta\partial_{\bar\zeta}t_0$ vanishes because of the factor $n=0$,
and the last second derivative vanishes because $k-2n\le1$ for $n=\lfloor k/2\rfloor$. Every derivative of $u_{\lambda,p,k}$ is a sum of terms
$\zeta^{\lambda+p+n-1}\times$(polynomial), and $p+n\ge0$. Since $|\zeta^{\lambda-1}|=|\zeta|^{\Re\lambda-1}e^{-\Im\lambda\arg\zeta}$ with $\Re\lambda-1>-1$ and $\arg\zeta$ bounded, these
functions are square integrable near the edge.
\end{proof}
\begin{proposition}[upper bound for the L-shaped prism]\label{oth:prop:lpup}
$\kappa(\Om_{\rm L})\ge0.913334979008$; hence $\beta(\Om_{\rm L})^2\le0.043332510496$ and $\beta(\Om_{\rm L})\le0.208165$.
\end{proposition}
\begin{proof}[Proof (computer-assisted)]
\emph{The pair} (\recR{2}). $F=F_{\rm pol}+F_{\rm sing}$ with $h$ even in $x_3$. $F_{\rm pol}$ is an exactly admissible polynomial pair of degree $36$ with rational
coefficients, constructed as in Section~\ref{oth:ssec:common}. $F_{\rm sing}$ is a combination with dyadic rational coefficients of the
real and imaginary parts of $924$ pairs \eqref{oth:eq:Pj}, built in the rotated frame from the potentials \eqref{oth:eq:upot} with $p\le2$, $k\le10$ and the seven exponents
\[\lambda\in\{\tfrac{5341}{10000},\ \tfrac{9514}{10000}\pm\tfrac{790}{10000}i,\ \tfrac23,\ \tfrac43,\ \tfrac{16273}{10000}\pm\tfrac{2550}{10000}i\},\]
with $\mathcal Q_0,\mathcal Q_3$ for odd $k$ and $\mathcal Q_1,\mathcal Q_2$ (taken in the rotated frame) for even $k$, which makes $h$ even in $x_3$. By Lemmas~\ref{oth:lem:pot}, \ref{oth:lem:upot} and the invariance of the
system under rotations (Lemma~\ref{lem:simil}), $F\in\M$; in addition $\Delta u=0$ and both equations were checked term by term in Gaussian rational arithmetic. The mean of $h$
is subtracted, so the pair lies in $\Mz$.

\emph{The quotient.} Integrals of products of two polynomial components are exact rationals. Every other integral is a product of an integral in $x_3$ of a polynomial
(exact) and a moment $\int_{\omega_{\rm L}}x_1^ix_2^j\mathcal K$ with $\mathcal K=\zeta_{\rm L}^\gamma$ or $\mathcal K=\zeta_{\rm L}^a\bar\zeta_{\rm L}^{\,b}$ (principal branches), which is homogeneous of a complex degree $d$ with
$\Re d\ge2(0.5341-1)>-2$. Lemma~\ref{lem:face}(b) reduces it to integrals over $E_1,\dots,E_4$ of functions analytic in a neighbourhood of each edge (the branch
points are at distance at least $1$); these are enclosed by rigorous integration in ball arithmetic \cite{Arb} at $200$ bits, with exponents and $\sqrt2$ entering as
balls. The quotient \eqref{kappa} of $F$ is then evaluated in ball arithmetic; its enclosure has radius $7\cdot10^{-47}$, and its lower end, rounded down to twelve
decimals, is $0.913334979008$. Theorem~\ref{thm:fv} gives $\beta^2\le0.043332510496$ and $\beta\le0.2081646236$.
\end{proof}
The polynomial part is proposed on the space of polynomial pairs of degree at most $36$ with $h$ even in $x_3$ (dimension $1406$), enlarged by the singular pairs.
The exponents are a heuristic choice, explained in \docD{6}; like the trace, they affect only the quality of the bound, not its validity.

\subsection{The Fichera domain}\label{oth:ssec:fich}
The symmetry group consists of the six coordinate permutations and contains rotations of order three, so Proposition~\ref{oth:prop:classes} applies only to its
subgroups of order two. The main step is that two classes suffice (Proposition~\ref{oth:prop:fichclasses}); it rests on Lemma~\ref{oth:lem:Esector}, proved on paper.
\emph{Geometry.} $\Om_{\rm F}=(-1,1)^3\setminus[0,1]^3$ is the cube with one closed octant removed. Its faces are the three squares $\{x_i=-1\}$, the three L-shaped
faces $\{x_i=1\}\cap\partial\Om_{\rm F}$ and the three inner squares $\{x_i=0\}\cap[0,1]^3$. Its edges of opening $\frac{3\pi}2$ are the three edges
$\{x_i=x_j=0,\,0<x_k<1\}$ of the removed octant; its eighteen other edges have opening $\frac\pi2$. Its vertices are ten octant vertices, the three L-cone vertices
$(1,0,0)$, $(0,1,0)$, $(0,0,1)$, where a re-entrant edge ends on an outer face, and the Fichera vertex $0$ (Table~\ref{oth:tab:geom}; the full list is in
\docD{7}).

\emph{Symmetry.} No reflection in a coordinate plane maps $\Om_{\rm F}$ onto itself. Let
\begin{equation}\label{oth:eq:GF}
\mathcal G_{\rm F}:=\{\Pi\in O(3):\ \Pi\ \text{a permutation matrix}\},\qquad \mathcal H_{\rm F}:=\{I,\Pi_{12}\}\subset\mathcal G_{\rm F},
\end{equation}
with $\Pi_{12}$ as in \eqref{oth:eq:GL}, and $\Pi_{23}(x_1,x_2,x_3):=(x_1,x_3,x_2)$. For $\Pi\in\mathcal G_{\rm F}$, $\det\Pi$ is the sign of the permutation and $\tr\Pi$
the number of its fixed indices. The fundamental regions are
\begin{equation}\label{oth:eq:FF}
\Om_{\mathcal H_{\rm F}}:=\{x\in\overline{\Om_{\rm F}}:\ x_1\le x_2\},\qquad \Om_{\mathcal G_{\rm F}}:=\{x\in\overline{\Om_{\rm F}}:\ x_1\le x_2\le x_3\},
\end{equation}
of volumes $\frac72$ and $\frac76$; the images $\Pi\,\Om_{\mathcal G_{\rm F}}$, $\Pi\in\mathcal G_{\rm F}$, are the six chambers $\{x_{\pi(1)}\le x_{\pi(2)}\le x_{\pi(3)}\}$ of the walls
$\{x_i=x_j\}$ intersected with $\overline{\Om_{\rm F}}$. For $e\in\{0,1\}$ the classes are
\begin{equation}\label{oth:eq:fichclasses}
\begin{gathered}
\Mz^{{\rm Q}e}:=\{F\in\Mz:\ \mathcal T_{\Pi_{12}}F=(-1)^eF\},\\
\Mz^{{\rm S}e}:=\{F\in\Mz:\ \mathcal T_\Pi F=(\det\Pi)^eF\ \text{for all }\Pi\in\mathcal G_{\rm F}\},
\end{gathered}
\end{equation}
with suprema $\kappa_{{\rm Q}e}$, $\kappa_{{\rm S}e}$ defined as in \eqref{oth:eq:class}. $\Mz^{{\rm Q}e}$ is the class of the character $\chi(\Pi_{12})=(-1)^e$ of
$\mathcal H_{\rm F}$, and $\Mz^{{\rm S}e}$ that of the character $\chi=\det^e$ of $\mathcal G_{\rm F}$; the labels agree with those of the records. Since
$\mathcal T_\Pi F=\chi(\Pi)F$ means $h(\Pi x)=\chi(\Pi)\det\Pi\,h(x)$, $h$ is odd under $x_1\leftrightarrow x_2$ in ${\rm Q}0$ and even in ${\rm Q}1$; it is symmetric under
all permutations in ${\rm S}1$ and alternating in ${\rm S}0$. The constants belong to ${\rm Q}1$ and ${\rm S}1$ only. The certificate uses ${\rm Q}0$ (group
$\mathcal H_{\rm F}$, region $\Om_{\mathcal H_{\rm F}}$) and ${\rm S}1$ (group $\mathcal G_{\rm F}$, region $\Om_{\mathcal G_{\rm F}}$). $\mathcal H_{\rm F}$ is a group of
reflections, so Proposition~\ref{oth:prop:classes} applies to ${\rm Q}0$. For ${\rm S}1$ we need the following extension to groups that contain rotations.
\begin{lemma}[one-dimensional classes of a finite group]\label{oth:lem:genclass}
Let $\mathcal G$ be a finite group of orthogonal $3\times3$ matrices with $R\Om=\Om$ for all $R\in\mathcal G$, not necessarily with $R^2=I$, and let
$\chi:\mathcal G\to\{\pm1\}$ be a homomorphism. Define $\Mz^\chi$, $P_\chi$, $\kappa_\chi$ by \eqref{oth:eq:class}, and the fundamental region $\Om_{\mathcal G}$, the
orbit size and the parity \eqref{A:parity} as in Definition~\ref{A:def:classes}.
\begin{enumerate}
\item[(a)] $P_\chi$ is an orthogonal projection on $L^2(\Om)^4$ that acts on $h$ and $g$ separately, $\mathcal T_RP_\chi=\chi(R)P_\chi$ for $R\in\mathcal G$, and
$P_\chi\Mz=\Mz^\chi$.
\item[(b)] Proposition~\ref{oth:prop:classes}(b),(c) hold for $\mathcal G$ and $\chi$, with (iii) replaced by: $C'$ is given on $\Om_{\mathcal G}$ and extended by
$C'(Rx):=\Lambda_RC'(x)\Lambda_R^{\mathsf T}$.
\item[(c)] Let $\xi=(\varphi,\vartheta)$ be continuous and piecewise polynomial on $\Om_{\mathcal G}$, zero on $\partial\Om\cap\Om_{\mathcal G}$, and such that
$\varphi(Rx)=\chi(R)\det R\,R\varphi(x)$ and $\vartheta(Rx)=\chi(R)\vartheta(x)$ whenever $x$ and $Rx$ lie in $\Om_{\mathcal G}$. Then
\begin{equation}\label{oth:eq:corrext}
\varphi(Rx):=\chi(R)\det R\,R\,\varphi(x),\qquad \vartheta(Rx):=\chi(R)\,\vartheta(x),\qquad x\in\Om_{\mathcal G},\ R\in\mathcal G,
\end{equation}
defines an extension $\xi\in H^1_0(\Om)^3\times H^1_0(\Om)$ with the parity \eqref{A:parity} of the class $\chi$.
\end{enumerate}
\end{lemma}
The proof is routine and given in Appendix~\ref{A:app:genclass}. For $\Om_{\mathcal H_{\rm F}}$ and $\Om_{\mathcal G_{\rm F}}$ the hypothesis of (c) reduces to
conditions on the walls. If $x$ and $Rx$ both satisfy $x_1\le x_2$ (resp.\ $x_1\le x_2\le x_3$), then $Rx=x$. The stabiliser of $x$ is generated by $\Pi_{12}$ if
$x_1=x_2$ and, for $\mathcal G_{\rm F}$, by $\Pi_{23}$ if $x_2=x_3$. Since $R\mapsto\chi(R)\det R\,R$ is a homomorphism, it suffices that
$\varphi=-\chi(\Pi)\Pi\varphi$ and $\vartheta=\chi(\Pi)\vartheta$ on the wall fixed by $\Pi\in\{\Pi_{12},\Pi_{23}\}$. For ${\rm Q}0$ on $\{x_1=x_2\}$ this is
$\varphi_2=-\varphi_1$, $\varphi_3=0$; for ${\rm S}1$ it is $\varphi_1=\varphi_2$, $\vartheta=0$ on $\{x_1=x_2\}$ and $\varphi_2=\varphi_3$, $\vartheta=0$ on $\{x_2=x_3\}$.

\emph{Why two classes suffice.} Under $\mathcal G_{\rm F}$, $\Mz$ splits into two one-dimensional parts and a part of type $E$, the two-dimensional irreducible
representation of $S_3\cong\mathcal G_{\rm F}$. ${\rm Q}0$ and ${\rm Q}1$ each contain one half of the $E$-part; both halves have the same $\kappa$
(Lemma~\ref{oth:lem:Esector}). The characters are $\chi_{\rm triv}=1$, $\chi_{\rm sgn}=\det$ and
$\chi_E(\Pi)=\tr\Pi-1$ (the values $2$, $0$, $-1$ on the identity, the transpositions and the rotations of order three; their orthogonality relations
\cite[Sect.~2.3]{Serre} are checked directly). Put
\begin{equation}\label{oth:eq:isotypic}
\begin{gathered}
P_\rho:=\frac{d_\rho}6\sum_{\Pi\in\mathcal G_{\rm F}}\chi_\rho(\Pi)\,\mathcal T_\Pi,\qquad X_\rho:=P_\rho\Mz,\\
P_\pm:=\tfrac12\big(I\pm\mathcal T_{\Pi_{12}}\big),\qquad X_E^\pm:=P_\pm X_E ,
\end{gathered}
\end{equation}
for $\rho\in\{{\rm triv},{\rm sgn},E\}$ with $d_{\rm triv}=d_{\rm sgn}=1$ and $d_E=2$,
and let $\mathcal E_hF:=(h,0)$ and $\mathcal E_gF:=(0,g)$ be the component projections, so that $\|h\|^2=(F,\mathcal E_hF)$ and $\|g\|^2=(F,\mathcal E_gF)$.
\begin{lemma}[isotypic decomposition]\label{oth:lem:isotypic}
The $P_\rho$ are orthogonal projections on $L^2(\Om_{\rm F})^4$ with $P_\rho P_\sigma=0$ for $\rho\ne\sigma$ and $\sum_\rho P_\rho=I$; they commute with every
$\mathcal T_\Pi$, with $P_\pm$, $\mathcal E_h$ and $\mathcal E_g$, and map $\Mz$ into itself. Moreover $X_{\rm triv}=\Mz^{{\rm S}0}$, $X_{\rm sgn}=\Mz^{{\rm S}1}$,
\[\Mz=X_{\rm triv}\oplus X_{\rm sgn}\oplus X_E^+\oplus X_E^-,\qquad \Mz^{{\rm Q}0}=X_{\rm triv}\oplus X_E^+,\qquad \Mz^{{\rm Q}1}=X_{\rm sgn}\oplus X_E^-,\]
and the sums are orthogonal both for $\int h\,h'$ and for $\int g\cdot g'$.
\end{lemma}
\begin{proof}
$V_\Pi:=\mathcal T_{\Pi^{-1}}$ is a unitary representation of $\mathcal G_{\rm F}$ on $L^2(\Om_{\rm F})^4$ (Lemma~\ref{A:lem:sym}(a) and
$\mathcal T_A\mathcal T_B=\mathcal T_{BA}$), which preserves $\Mz$ and commutes with $\mathcal E_h$, $\mathcal E_g$; it extends to an algebra homomorphism
$a=\sum_\Pi a_\Pi\Pi\mapsto V(a):=\sum_\Pi a_\Pi V_\Pi$ of the group algebra $\mathbb Q[\mathcal G_{\rm F}]$, with $V(a)^*=V(a^*)$ for
$a^*:=\sum_\Pi a_\Pi\Pi^{-1}$. The characters are real and $\chi_\rho(\Pi^{-1})=\chi_\rho(\Pi)$, so $P_\rho=V(e_\rho)$ and $P_\pm=V(p_\pm)$ with
$e_\rho:=\frac{d_\rho}6\sum_\Pi\chi_\rho(\Pi)\Pi=e_\rho^*$ and $p_\pm:=\frac12(I\pm\Pi_{12})=p_\pm^*$. In $\mathbb Q[\mathcal G_{\rm F}]$ the $e_\rho$ are central,
$e_\rho^2=e_\rho$, $e_\rho e_\sigma=0$ for $\rho\ne\sigma$ and $\sum_\rho e_\rho=1$ (the projection formula \cite[Sect.~2.6, Thm.~8]{Serre} for the regular
representation, which is faithful on $\mathbb Q[\mathcal G_{\rm F}]$); moreover $p_-e_{\rm triv}=0$ and $p_+e_{\rm sgn}=0$, because $\Pi_{12}e_{\rm triv}=e_{\rm triv}$
and $\Pi_{12}e_{\rm sgn}=-e_{\rm sgn}$. Applying $V$ gives the first statement. By Lemma~\ref{oth:lem:genclass}(a), $P_{\rm triv}$, $P_{\rm sgn}$ and $P_\pm$ are the
projections $P_\chi$ of the classes ${\rm S}0$, ${\rm S}1$, ${\rm Q}0$, ${\rm Q}1$, so $X_{\rm triv}=\Mz^{{\rm S}0}$, $X_{\rm sgn}=\Mz^{{\rm S}1}$ and
$\Mz^{{\rm Q}0}=P_+\Mz=\sum_\rho P_+P_\rho\Mz=X_{\rm triv}\oplus X_E^+$, using $P_+P_{\rm triv}=P_{\rm triv}$ and $P_+P_{\rm sgn}=0$; likewise for ${\rm Q}1$. The
projections onto the summands, $P_{\rm triv}$, $P_{\rm sgn}$, $P_+P_E$, $P_-P_E$, are self-adjoint, mutually annihilating and commute with $\mathcal E_h$ and
$\mathcal E_g$; hence $(F,\mathcal E_hF')=0$ and $(F,\mathcal E_gF')=0$ for $F$, $F'$ in different summands.
\end{proof}
\begin{lemma}[the $E$-sector]\label{oth:lem:Esector}
$\kappa(X_E^+)=\kappa(X_E^-)$, where, for a subspace $X\subset\Mz$, $\kappa(X):=\sup_{F\in X\setminus\{0\}}(\|h\|^2-\|g\|^2)/(\|h\|^2+\|g\|^2)$ and
$\kappa(\{0\}):=-\infty$.
\end{lemma}
\begin{proof}
Let $c\in\mathcal G_{\rm F}$ be a rotation of order three. In $\mathbb Q[\mathcal G_{\rm F}]$,
\begin{equation}\label{oth:eq:s3id}
e_E\,p_-\,c^{-1}p_+\,c\,p_-=\tfrac34\,e_E\,p_-,\qquad e_E\,p_+\,c^{-1}p_-\,c\,p_+=\tfrac34\,e_E\,p_+ .
\end{equation}
Both identities are verified by direct multiplication in the six-dimensional algebra, in exact rational arithmetic. They express a fact about
$E$, realised by the permutation matrices on the plane $\{x_1+x_2+x_3=0\}$: there $\Pi_{12}$ is a reflection with orthonormal eigenvectors $u_\pm$, $p_\pm$ acts as
$u_\pm u_\pm^{\mathsf T}$, $c$ is a rotation by $\frac{2\pi}3$, and $p_\mp c^{-1}p_\pm cp_\mp$ acts as $(u_\pm\cdot cu_\mp)^2u_\mp u_\mp^{\mathsf T}=\sin^2\!\frac{2\pi}3\,u_\mp
u_\mp^{\mathsf T}$. Let $F\in X_E^-$, so $F=P_EP_-F$, and put $\mathcal JF:=\frac2{\sqrt3}P_+V_cF$. Since $P_E$ commutes with $V_c$ and
$P_+$, $\mathcal JF\in P_+P_E\Mz=X_E^+$. For $\mathcal E\in\{\mathcal E_h,\mathcal E_g\}$, which commutes with all $V_\Pi$, using $V_c^*=V_{c^{-1}}$, $F=P_EP_-F$
and \eqref{oth:eq:s3id},
\[(\mathcal JF,\mathcal E\mathcal JF)=\tfrac43\big(F,\mathcal E\,V(e_Ep_-c^{-1}p_+cp_-)F\big)=\tfrac43\cdot\tfrac34\,(F,\mathcal E P_EP_-F)=(F,\mathcal EF).\]
So $\mathcal J$ maps $X_E^-$ into $X_E^+$ and preserves $\|h\|^2$ and $\|g\|^2$; hence $\kappa(X_E^-)\le\kappa(X_E^+)$. The second identity of \eqref{oth:eq:s3id}
gives the reverse inequality in the same way.
\end{proof}
\begin{proposition}[classes of the Fichera domain]\label{oth:prop:fichclasses}
For the classes \eqref{oth:eq:fichclasses},
\[\kappa(\Om_{\rm F})=\max\{\kappa_{{\rm Q}0},\kappa_{{\rm S}1}\}=\max\{\kappa_{{\rm Q}1},\kappa_{{\rm S}0}\}.\]
\end{proposition}
\begin{proof}
By Lemma~\ref{oth:lem:isotypic} and the argument of the proof of Proposition~\ref{oth:prop:classes}(a) (a quotient over an orthogonal sum for both forms is at most the
largest quotient over the summands, and each summand is a subspace of $\Mz$),
\begin{gather*}
\kappa(\Om_{\rm F})=\max\{\kappa(X_{\rm triv}),\kappa(X_{\rm sgn}),\kappa(X_E^+),\kappa(X_E^-)\},\\
\kappa_{{\rm Q}0}=\max\{\kappa(X_{\rm triv}),\kappa(X_E^+)\},\qquad \kappa_{{\rm Q}1}=\max\{\kappa(X_{\rm sgn}),\kappa(X_E^-)\},\\
\kappa_{{\rm S}0}=\kappa(X_{\rm triv}),\qquad \kappa_{{\rm S}1}=\kappa(X_{\rm sgn}).
\end{gather*}
Lemma~\ref{oth:lem:Esector} gives both formulas.
\end{proof}
Without Lemma~\ref{oth:lem:Esector} a certificate would need both ${\rm Q}0$ and ${\rm Q}1$; the lemma replaces ${\rm Q}1$ by the smaller class ${\rm S}1$
($\FichNS$ instead of $\FichNQone$ functionals in the design of Table~\ref{oth:tab:design}).

\emph{Upper bound.}
\begin{proposition}[upper bound for the Fichera domain]\label{oth:prop:fichup}
$\kappa(\Om_{\rm F})\ge0.802453899040$; hence $\beta(\Om_{\rm F})^2\le0.098773050480$ and $\beta(\Om_{\rm F})\le0.314282$.
\end{proposition}
\begin{proof}[Proof (computer-assisted)]
The trial pair (\recR{1}) is an exactly admissible polynomial pair of degree $20$ with rational coefficients, constructed as in Section~\ref{oth:ssec:common}, without
symmetry reduction; the mean of $h$ is subtracted exactly. Its quotient is a rational number, computed by exact integration of polynomials over $\Om_{\rm F}$, a union of
seven unit cubes; rounded down to twelve decimals it is $0.802453899040$. Theorem~\ref{thm:fv} gives $\beta^2\le0.098773050480$ and $\beta\le0.3142818011$.
\end{proof}
The polynomial trial values have not converged at $d=20$; they are listed for lower degrees in \docD{6}.

\subsection{The certified lower bounds}\label{oth:ssec:cert}
Proposition~\ref{oth:prop:cert} states the conditions under which a bound $\lambda_{\rm up}$ gives $\kappa_\chi\le\kappa_0$; the paragraph after it says how the verifier checks
them. The domain enters only through the data. Table~\ref{oth:tab:design} in Appendix~\ref{app:design} lists the data of each record, and
Section~\ref{tab:certs} explains where their hashes are listed; the centres and weights of the balls are in the data files of the records.
\begin{proposition}[verified certificate]\label{oth:prop:cert}
Let $\Om$ be $\Om_{\rm B}$, $\Om_{\rm L}$ or $\Om_{\rm F}$, let $\mathcal G$ be one of its groups $\mathcal G_{\rm B}$, $\mathcal G_{\rm L}$, $\mathcal H_{\rm F}$,
$\mathcal G_{\rm F}$ with fundamental region $\Om_{\mathcal G}$, let $\chi$ be a character of $\mathcal G$ and $\kappa_0$ rational. Let the data of
Definition~\ref{A:def:data} consist of a mesh of $\Om_{\mathcal G}$ into closed tetrahedra, a continuous multiplier $m$ affine on each tetrahedron, balls with weights and
harmonic bases, and continuous piecewise quadratic correctors. Assume:
(E1) the tetrahedra tile $\Om_{\mathcal G}$;
(E2) $m=0$ at the mesh vertices on $\partial\Om$, and $m(x)=Rm(x)$ at the mesh vertices $x$ on the mirror wall of $\Om_{\mathcal G}$ fixed by $R\in\mathcal G$;
(E3) the balls form $\mathcal G$-orbits with equal $r_k,R_k,\alpha_k,p_k$, $B(x_k,R_k)\subset\Om$ and $r_k\le R_k$;
(E4) the harmonic bases have degree $\le p_k$ and are homogeneous, harmonic, orthogonal and of the parity of the class; the representers are projected with the exact
coefficients of $P_\chi$; the correctors vanish on $\partial\Om$ and have the parity of the class on the mirror walls (as after Lemma~\ref{oth:lem:genclass});
(E5) $B_P\succ0$ on all pieces.
Let $G^\chi$ be the Gram matrix of Proposition~\ref{oth:prop:classes}(b),(c) if $\mathcal G$ is a group of reflections, and of Lemma~\ref{oth:lem:genclass}(b) if
$\mathcal G=\mathcal G_{\rm F}$. If $\lambda_{\max}(G^\chi)\le\lambda_{\rm up}\le1$, then $\kappa_\chi\le\kappa_0$.
\end{proposition}
\begin{proof}
By (E1), (E2) the $\mathcal G$-equivariant extension of $m$ is continuous, affine on each tetrahedron $R\tau$ and zero on $\partial\Om$. By Lemma~\ref{oth:lem:lip} it is
Lipschitz on $\overline\Om$, also when $\Om$ is not convex, and it satisfies (i) of Proposition~\ref{oth:prop:classes}(b). (E3) gives (ii) and the containment needed by
Lemma~\ref{lem:ball} and \eqref{rep3}. By (E5) and Lemma~\ref{B-lem:Cprime}, $0\prec C'\preceq C-\rho_{\rm up}e_0e_0^{\mathsf T}$ on $\Om_{\mathcal G}$, and $C'$ is
extended as in (iii), resp.\ by $C'(Rx)=\Lambda_RC'(x)\Lambda_R^{\mathsf T}$. By (E4) the representers are projected onto the class, as in (iv), and by
Lemma~\ref{oth:lem:genclass}(c) the correctors, extended by \eqref{oth:eq:corrext}, lie in $H^1_0(\Om)^3\times H^1_0(\Om)$ and have the parity of the class.
Proposition~\ref{oth:prop:classes}(b),(c), resp.\ Lemma~\ref{oth:lem:genclass}(b), gives $\kappa_\chi\le\kappa_0$.
\end{proof}
\emph{How the conditions are checked.} The verifier checks (E1)--(E5) in integer and rational arithmetic. For the non-convex regions of $\Om_{\rm L}$ and $\Om_{\rm F}$ it
checks (E1) through six conditions: all vertices lie in $\Om_{\mathcal G}$, no tetrahedron is degenerate, the volumes sum to $|\Om_{\mathcal G}|$, each face belongs to at most
two tetrahedra, each boundary face lies in one closed convex boundary piece, and the two tetrahedra of an interior face lie on opposite sides of it. They imply the tiling:
the covering number is locally constant off the boundary faces, vanishes outside, and equals one by the volume identity. (E5) is checked by exact leading principal minors.
$G^\chi$ is accumulated on $\Om_{\mathcal G}$ with the factor $|\mathcal G|$ in interval arithmetic, with Grundmann--M\"oller rules of orders $s$ checked in rational
arithmetic to be exact to degree $2s+1$ and exact for the degree of every integrand. The interval sums therefore enclose the entries of $G^\chi$, and $\lambda_{\rm up}$ is the
bound of Section~\ref{sec:verif}. All three records use the representer weight $\mu=2$ with $R'=\frac{63}{64}R$; the containment tests of (E3) are listed in \docD{7}.
\begin{table}[t]
\caption{Certified bounds $\lambda_{\rm up}\ge\lambda_{\max}(G^\chi)$ per class (rounded up to nine decimals) and number $N$ of functionals. Symmetry of $h$ (o: odd,
e: even): box, parity in $(x_1,x_2,x_3)$; L-prism, class $00s_3e$, parity in $x_3$ and under $x_1\leftrightarrow x_2$; Fichera domain, group $\mathcal H_{\rm F}$ and
region $\Om_{\mathcal H_{\rm F}}$ for ${\rm Q}0$, group $\mathcal G_{\rm F}$ and region $\Om_{\mathcal G_{\rm F}}$ for ${\rm S}1$ \eqref{oth:eq:FF}.}\label{oth:tab:cert}\label{oth:tab:boxcert}\label{oth:tab:lpcert}\label{oth:tab:fichcert}
\small
\begin{tabular}{@{}llllr@{}}
\toprule
domain, $\kappa_0$ & class & symmetry of $h$ & $\lambda_{\rm up}$ & $N$\\
\midrule
box, $\frac{87}{100}$ & $000$ & (o,o,o) & $0.058763585$ & 244\\
 & $100$ & (e,o,o) & $0.058818125$ & 246\\
 & $010\cong001$ & (o,e,o) & $0.123335403$ & 265\\
 & $110\cong101$ & (e,e,o) & $0.121787006$ & 265\\
 & $011$ & (o,e,e) & $0.968926254$ & 286\\
 & $111$ & (e,e,e) & $0.260628396$ & 289\\
\midrule
L-prism, $\frac{37}{40}$ & $0000$ & (o; o) & $0.076857855$ & 208\\
 & $0001$ & (o; e) & $0.076664915$ & 212\\
 & $0010$ & (e; o) & $0.982951223$ & 252\\
 & $0011$ & (e; e) & $0.485093615$ & 256\\
\midrule
Fichera, $\FichKappaZero$ & ${\rm Q}0$ & odd under $x_1\leftrightarrow x_2$ & $\FichLamQ$ & $\FichNQ$\\
 & ${\rm S}1$ & symmetric & $\FichLamS$ & $\FichNS$\\
\bottomrule
\end{tabular}
\end{table}

\begin{theorem}[lower bound for the box]\label{oth:thm:box}
$\kappa_s\le\frac{87}{100}$ for every $s\in\mathcal S_{\rm B}$; hence $\kappa(\Om_{\rm B})\le\frac{87}{100}$ and $\beta(\Om_{\rm B})^2\ge\frac{13}{200}$.
\end{theorem}
\begin{theorem}[lower bound for the L-shaped prism]\label{oth:thm:lp}
$\kappa_\chi\le\frac{37}{40}$ for every $\chi\in\mathcal S_{\rm L}$; hence $\kappa(\Om_{\rm L})\le\frac{37}{40}$ and $\beta(\Om_{\rm L})^2\ge\frac3{80}$.
\end{theorem}
\begin{theorem}[lower bound for the Fichera domain]\label{oth:thm:fich}
$\kappa_{{\rm Q}0}\le\kappa_0^{\rm F}$ and $\kappa_{{\rm S}1}\le\kappa_0^{\rm F}$ for $\kappa_0^{\rm F}:=\FichKappaZero$; hence $\kappa(\Om_{\rm F})\le\kappa_0^{\rm F}$,
$\beta(\Om_{\rm F})^2\ge\frac12(1-\kappa_0^{\rm F})=\FichBetaSqLow$ and $\beta(\Om_{\rm F})\ge\FichBetaLow$.
\end{theorem}
\begin{proof}[Proof of Theorems~\ref{oth:thm:box}--\ref{oth:thm:fich} (computer-assisted)]
The records are \recR{9}, \recR{10} and \recR{11}. In each case the verifier checks (E1)--(E5) of Proposition~\ref{oth:prop:cert} and bounds $\lambda_{\max}(G^\chi)$ for
the classes of Table~\ref{oth:tab:cert}. All bounds are below $1$, the largest being $0.968926254$ (box, class $011$), $0.982951223$ (L-prism, class $0010$) and $\FichLamQ$
(Fichera domain, class $\FichBindClass$). Proposition~\ref{oth:prop:cert} gives $\kappa_\chi\le\kappa_0$ for these classes; Corollaries~\ref{oth:cor:boxclasses},
\ref{oth:cor:lpclasses} and Proposition~\ref{oth:prop:fichclasses} give the bounds for $\kappa(\Om)$, and Theorem~\ref{thm:fv} those for $\beta$. The domain-specific checks
are as follows.

\emph{Box:} $\mathcal G=\mathcal G_{\rm B}$ on the octant, with $R_k=\min_i(\frac{w_i}2-x_{k,i})=\dist(x_k,\partial\Om)$ exactly for the side lengths
$(w_1,w_2,w_3)=(2,1,1)$; rule orders up to $9$.
\emph{L-shaped prism:} the mesh, the multiplier and the balls on $\omega_{\rm L}\times[0,\frac12]$ are checked to be exactly invariant under $\Pi_{12}$ before they are
restricted to $\Om_{\mathcal G}$ \eqref{oth:eq:FL}. The coefficients of $P_\chi$ are exact multiples of $\frac14$, and the coinciding functionals of the two diagonal
balls are merged; rule orders up to $8$.
\emph{Fichera domain:} both classes are verified in one run. It reports a result only if the verified classes contain $\{{\rm Q}0,{\rm S}1\}$ or $\{{\rm Q}1,{\rm S}0\}$ and
all of them pass. The data are given on $\Om_{\rm F}$. Before they are restricted to a region \eqref{oth:eq:FF}, the vertex set, the tetrahedra and the multiplier are
checked to be invariant under all six permutations, and no tetrahedron crosses a plane $\{x_i=x_j\}$. One ball per orbit of the group of the class is kept, and the orbits
of the kept balls are checked to be exactly the balls of the data. The coefficients of $P_\chi$ are exact multiples of $\frac12$ for $\mathcal H_{\rm F}$ (class ${\rm Q}0$)
and of $\frac16$ for $\mathcal G_{\rm F}$ (class ${\rm S}1$; the projected representers are multiplied by $3$ and their weights divided by $9$, which leaves
$\alpha_\iota\ell_\iota(F)^2$ unchanged). $\FichMergedQ$ functionals are merged in ${\rm Q}0$ and $\FichMergedS$ in ${\rm S}1$; rule orders $s\in\FichRules$.
\end{proof}
\begin{corollary}[enclosure for the box]\label{oth:cor:box}
$0.856921056806\le\kappa(\Om_{\rm B})\le0.87$ and $0.254950\le\beta(\Om_{\rm B})\le0.267469$, a relative width $(0.267469-0.254950)/\allowbreak 0.267469<4.69\cdot10^{-2}$.
Moreover $\frac{13}{200}\le\beta(\Om_{\rm B})^2\le0.071539471597$ and $\frac{13}{200}\le\inf\sigma_{\rm ess}(\mathcal S_{\Om_{\rm B}})\le0.116397718977$.
\end{corollary}
\begin{proof}
Proposition~\ref{oth:prop:boxup} and Theorem~\ref{oth:thm:box}; $0.254950^2=0.0649995025\le\frac{13}{200}$ exactly. The last statement is
Corollary~\ref{B-cor:cubeess}(b) with $\beta^2\ge\frac{13}{200}$.
\end{proof}
\begin{corollary}[enclosure for the L-shaped prism]\label{oth:cor:lp}
$0.913334979008\le\kappa(\Om_{\rm L})\le0.925$ and $0.193649\le\beta(\Om_{\rm L})\le0.208165$, a relative width
$(0.208165-0.193649)/\allowbreak 0.208165<6.98\cdot10^{-2}$. Moreover $\frac3{80}\le\beta(\Om_{\rm L})^2\le0.043332510496$ and
$\frac3{80}\le\inf\sigma_{\rm ess}(\mathcal S_{\Om_{\rm L}})\le0.116397718977$.
\end{corollary}
\begin{proof}
Proposition~\ref{oth:prop:lpup} and Theorem~\ref{oth:thm:lp}; $0.193649^2=0.037499935201\le\frac3{80}$ exactly. The last statement is Corollary~\ref{B-cor:cubeess}(b)
with $\beta^2\ge\frac3{80}$.
\end{proof}
\begin{corollary}[enclosure for the Fichera domain]\label{oth:cor:fich}
$0.802453899040\le\kappa(\Om_{\rm F})\le\kappa_0^{\rm F}$ and $\FichBetaLow\le\beta(\Om_{\rm F})\le0.314282$, a relative width
$(0.314282-\FichBetaLow)/\allowbreak 0.314282<\FichWidth$. Moreover $\frac12(1-\kappa_0^{\rm F})\le\beta(\Om_{\rm F})^2\le0.098773050480$ and
$\frac12(1-\kappa_0^{\rm F})\le\inf\sigma_{\rm ess}(\mathcal S_{\Om_{\rm F}})\le0.116397718977$.
\end{corollary}
\begin{proof}
Proposition~\ref{oth:prop:fichup}, Theorem~\ref{oth:thm:fich} and Corollary~\ref{B-cor:cubeess}(b); exactly, $\FichBetaLow^2=\FichBetaLowSq\le\FichBetaSqLow$.
\end{proof}

\begin{remark}[what limits the certificates]\label{oth:rem:limit}\label{oth:rem:boxlimit}\label{oth:rem:lplimit}\label{oth:rem:fichlimit}\label{oth:rem:local}
For the three domains the trial values of $\kappa$ ($0.856921$, $0.913334$, $0.802453$; Table~\ref{tab:poly}) exceed every local constant of Table~\ref{oth:tab:local}. For
the cube, in contrast, $\kappa$ and the octant constant agree to about $2\cdot10^{-6}$ (Section~\ref{sec:vertex}). Near every edge and vertex a multiplier can therefore keep
$\lambda_{\max}(W)$ below the certified $\kappa_0$. The certificates are limited by global quantities: the gap between $\kappa_0$ and the top of the spectrum, the corrector
resolution and the piecewise weight. In particular, the excess of the Fichera-cone value over $2/(3\pi)$ does not limit a certificate near $\kappa\approx0.80$. The gaps
in $\kappa$ to the trial values are $0.013078943194$, $0.011665020992$ and $\FichGap$. For the box and the L-prism no $\kappa_0$ below $\kappa_{\rm PD}^{\rm ex}$
\eqref{oth:eq:kpd} ($>0.866$, resp.\ $\approx0.9143$) can be certified with these data, whatever the correctors. For the Fichera domain $B_P$ is positive definite down to
$\frac{83}{100}$ at least, and the limit is the margin $\FichMarginQ$ of $\lambda_{\rm up}$ in class ${\rm Q}0$. The boundary floors $\kappa_\partial=0.84982$, $0.89490$,
$\FichKpartial$ lie below the certified values. Lowering $\kappa_0$ therefore needs a finer design in the interior, not a different treatment of edges or vertices
(details in \docD{7}).
\end{remark}

\subsection{Summary}\label{oth:ssec:summary}\label{thm:other}
Table~\ref{oth:tab:summary} collects the bounds for the four domains. For the cube the lower bound is limited by the octant vertex, and $\beta^2$ lies within
$4\cdot10^{-4}$ of the vertex value $\frac12(1-L)$ (Section~\ref{sec:vertex}). For the three other domains $\beta^2$ lies strictly below that value
(Corollary~\ref{B-cor:cubeess}(d)). Their certificates are limited by the design (Remark~\ref{oth:rem:limit}), and their relative widths are larger: $4.69\cdot10^{-2}$,
$6.98\cdot10^{-2}$ and $\FichWidth$. The upper bound $0.116397718977$ for $\inf\sigma_{\rm ess}(\mathcal S)$ is the same for all four domains, because each of them has a
convex rectangular vertex (Corollary~\ref{cor:ess}).
\begin{table}[t]
\caption{Summary. Lower bounds for $\beta$ from the certified $\kappa_0$, upper bounds from trial pairs; width: relative width $(\beta_{\rm up}-\beta_{\rm low})/\beta_{\rm up}$,
rounded up to three digits; upper bound for $\beta^2$; lower bound $\beta_{\rm low}^2$ for $\inf\sigma_{\rm ess}(\mathcal S)$; binding class of the certificate.}\label{oth:tab:summary}
\scriptsize
\setlength{\tabcolsep}{2pt}
\begin{tabular}{@{}llllllll@{}}
\toprule
domain & $\kappa_0$ & $\beta\ge$ & $\beta\le$ & width & $\beta^2\le$ & $\inf\sigma_{\rm ess}\ge$ & class\\
\midrule
cube & $\frac{96}{125}$ & $0.340587$ & $0.3411711$ & $1.72\cdot10^{-3}$ & $0.116397718977$ & $\frac{29}{250}$ & $110$\\
box & $\frac{87}{100}$ & $0.254950$ & $0.267469$ & $4.69\cdot10^{-2}$ & $0.071539471597$ & $\frac{13}{200}$ & $011$\\
L-prism & $\frac{37}{40}$ & $0.193649$ & $0.208165$ & $6.98\cdot10^{-2}$ & $0.043332510496$ & $\frac3{80}$ & $0010$\\
Fichera & $\FichKappaZero$ & $\FichBetaLow$ & $0.314282$ & $\FichWidth$ & $0.098773050480$ & $\FichBetaSqLow$ & $\FichBindClass$\\
\bottomrule
\end{tabular}
\end{table}

\setlength{\emergencystretch}{0pt}\section{Verified computation and formal checks}\label{sec:verif}
This section shows how the hypotheses of Theorem~\ref{thm:rel} are verified. Floating point only proposes data: the multiplier, the weights, the correctors, the
refinement and the eigenvector matrix of the final test. Everything on the proof path is exact or enclosed.

\subsection{From pointwise enclosures to the matrix test}
The verifier works on a subdivision of the tetrahedra of the mesh into simplices, the \emph{pieces} $P$ (Table~\ref{tab:certlow}); on $P$, $W(m)=W_\tau$
\eqref{A:Wtau} for the tetrahedron $\tau\supseteq P$. The proof needs, on each piece, $B_P\succ0$ with $B_P$ of \eqref{B-eq:BP}, $\kappa_0'\le\kappa_0$,
$W_{{\rm up},P}\succeq W_\tau$, $0\le\phi_P$ and $\phi_P^{\min}+\phi_P\le\Phi-\rho_{\rm up}$ on $P$. These give $0\prec C'\preceq C-\rho_{\rm up}e_0e_0^{\mathsf T}$
(Lemma~\ref{B-lem:Cprime}). Integration must then give a \emph{quadratic-form} bound $G\preceq G_{\rm up}$, since Theorem~\ref{thm:rel} uses
$y^{\mathsf T}Gy\le y^{\mathsf T}G_{\rm up}y$ for all $y$. Entrywise bounds alone are not enough; Lemma~\ref{B-lem:gramerr} turns them into such a bound.

The computation stores $2^{50}B_P$ as an integer matrix and checks $B_P\succ0$ by exact leading principal minors. It integrates by cubature rules exact for the
degree of the integrand and accumulates in interval arithmetic. The data checks (tiling of the fundamental region, boundary and mirror conditions on the
multiplier and the correctors, containment of the balls, harmonic basis) are done in integer and rational arithmetic; for the box, the L-shaped prism and the
Fichera domain they are listed in the proofs of Theorems~\ref{oth:thm:box}, \ref{oth:thm:lp} and~\ref{oth:thm:fich}. Parts (a)--(c) of the next lemma are
standard (proofs in Appendix~\ref{B-app:verif}).
\begin{lemma}\label{lem:verif}
(a) (Sherman--Morrison \cite[Sect.~2.1]{GolubVanLoan}.) Let $B\succ0$ be a symmetric $4\times4$ matrix, $u:=B^{-1}e_0$, $b_{00}:=e_0^{\mathsf T}B^{-1}e_0$, and $0\le\phi\le\phi_{\max}$. Then
\[(B+\phi\, e_0e_0^{\mathsf T})^{-1}=B^{-1}-\frac{\phi}{1+\phi\, b_{00}}\,uu^{\mathsf T}\preceq B^{-1}-\lambda_{\rm SM}\phi\,uu^{\mathsf T},\qquad\lambda_{\rm SM}:=\frac1{1+\phi_{\max}b_{00}} .\]
(b) (Young.) For a symmetric $M\succeq0$, vectors $u,v$ and $t>0$: $(u+v)^{\mathsf T}M(u+v)\le(1+t)\,u^{\mathsf T}Mu+(1+1/t)\,v^{\mathsf T}Mv$.
(c) (Gershgorin after congruence; cf.\ \cite[Thm.~6.1.1]{HornJohnson}.) Let $A$ be symmetric positive semidefinite, $V$ invertible, $Y:=V^{\mathsf T}AV$ and $\Xi:=V^{\mathsf T}V$. If $d_i:=\Xi_{ii}-\sum_{j\ne i}|\Xi_{ij}|>0$ for all $i$, then
$\lambda_{\max}(A)\le\max_i\big(Y_{ii}+\sum_{j\ne i}|Y_{ij}|\big)/d_i$.
\end{lemma}
\begin{lemma}[entrywise Gram error]\label{B-lem:gramerr}
Let $\hat G,\hat A$ be symmetric $N\times N$ matrices with $|\hat G_{\iota\iota'}-\hat A_{\iota\iota'}|\le\varepsilon_G$ for all $\iota,\iota'$, let $d\in\R^N$ and $D_d:=\operatorname{diag}(d_\iota)$. Then
\begin{equation}\label{B-eq:Gup}
D_d\hat GD_d\preceq G_{\rm up}:=D_d\hat AD_d+\varepsilon_G|d|^2I .
\end{equation}
\end{lemma}
\begin{proof}
For $y\in\R^N$, $\big|y^{\mathsf T}D_d(\hat G-\hat A)D_dy\big|\le\varepsilon_G\big(\sum_\iota|d_\iota||y_\iota|\big)^2\le\varepsilon_G|d|^2|y|^2$ by the Cauchy--Schwarz inequality.
\end{proof}

\subsection{Pointwise weight}
On a piece, $C'$ is a constant matrix plus a polynomial multiple of $e_0e_0^{\mathsf T}$; so the Gram integrands are bounded by polynomials (Lemma~\ref{B-lem:Cprime}).
\begin{definition}[implemented weight]\label{B-def:BP}
On a piece $P$ let $\kappa_0'\le\kappa_0$, let $W_{{\rm up},P}\succeq W_\tau$ be symmetric, and let $\phi_P^{\min}\in\R$ and a polynomial $\phi_P\ge0$ on $P$ satisfy $\phi_P^{\min}+\phi_P\le\Phi-\rho_{\rm up}$ on $P$. Here $\phi_P$ is the polynomial part of $\Phi$ above its minimum on $P$. Put
\begin{equation}\label{B-eq:BP}
B_P:=\kappa_0'I-W_{{\rm up},P}+\phi_P^{\min}\,e_0e_0^{\mathsf T},\qquad C'|_P:=B_P+\phi_P\,e_0e_0^{\mathsf T}.
\end{equation}
\end{definition}
\begin{lemma}\label{B-lem:Cprime}
If $B_P\succ0$ for every piece, then $C'$ of \eqref{B-eq:BP} satisfies $0\prec C'\preceq C-\rho_{\rm up}e_0e_0^{\mathsf T}$ almost everywhere, as required in Theorem~\ref{thm:rel}. Moreover, on each piece,
\[(C')^{-1}\preceq B_P^{-1}-\lambda_P\,\phi_P\,u_Pu_P^{\mathsf T},\qquad u_P:=B_P^{-1}e_0,\qquad \lambda_P:=\frac1{1+\phi_{P,\max}\,e_0^{\mathsf T}B_P^{-1}e_0},\]
with $\phi_{P,\max}\ge\max_P\phi_P$, so the Gram integrands are bounded by polynomials on $P$.
\end{lemma}
\begin{proof}
$C'\succeq B_P\succ0$ because $\phi_P\ge0$. On $P$, $C-\rho_{\rm up}e_0e_0^{\mathsf T}-C'=(\kappa_0-\kappa_0')I+(W_{{\rm up},P}-W_\tau)+(\Phi-\rho_{\rm up}-\phi_P^{\min}-\phi_P)e_0e_0^{\mathsf T}\succeq0$. The bound for $(C')^{-1}$ is Lemma~\ref{lem:verif}(a) with $B=B_P$, $\phi=\phi_P$ and $\lambda_{\rm SM}=\lambda_P$.
\end{proof}

\subsection{Gram integrals}
The rules are Grundmann--M\"oller rules \cite{GrundmannMoeller}, with exactness checked in rational arithmetic; they are selected from the polynomial degree of
the data, and the verifier refuses if an integrand exceeds their exactness. The representers use the radial weight $(1-s^2/R'^2)_+^{2}$ of
Section~\ref{sec:cert}, with $s$ the distance to the ball centre and $R'=\frac{63}{64}R$.

\emph{Pieces cut by a representer sphere.} On a piece $P$ cut by a sphere $\{|x-x_k|=R_k'\}$ the weight $\psi_k$ of \eqref{A:bump} is not a polynomial. For each
of the $n_P$ cutting spheres (reflected copies included) the verifier replaces $\psi_k$ by a polynomial branch ($\psi_k$ if the centroid of $P$ lies inside the
sphere, $0$ otherwise). The difference $e_c$ is \emph{one} scalar function per sphere, common to all modes of that ball, since they all carry the factor $\psi_k$.
So the remainder of $r_y$ is $\big(\sum_ce_c(x)(v_c(x)\cdot y)\big)e_0$, with $v_c(x)$ the exact polynomial vector of the modes of the ball. A bound
$\bar e_c\ge|e_c|$ on $P$ follows from the squared distances of the vertices of $P$ to the centre.
\begin{lemma}[pieces cut by a representer sphere]\label{B-lem:cut}
Let $P$ be a piece with $B_P\succ0$ and $C'=B_P+\phi_Pe_0e_0^{\mathsf T}$ on $P$ as in \eqref{B-eq:BP}. For $y\in\R^N$ let the combination of representers on $P$ be
\[r_y=r_y^{\rm pol}+\Big(\sum_{c=1}^{n_P}e_c\,(v_c\cdot y)\Big)e_0 ,\]
with $r_y^{\rm pol}$ and $v_c:P\to\R^N$ polynomial and measurable $e_c$ with $|e_c|\le\bar e_c$ on $P$. Then for every $t>0$,
\[\int_Pr_y^{\mathsf T}(C')^{-1}r_y\le(1+t)\int_P(r_y^{\rm pol})^{\mathsf T}(C')^{-1}r_y^{\rm pol}+(1+t^{-1})\,b_{00,P}\,n_P\sum_{c=1}^{n_P}\bar e_c^{\,2}\int_P(v_c\cdot y)^2 ,\]
where $b_{00,P}:=e_0^{\mathsf T}B_P^{-1}e_0$.
\end{lemma}
\begin{proof}
Pointwise on $P$, Lemma~\ref{lem:verif}(b) with $M=(C')^{-1}$ gives the factor $1+t$ for $r_y^{\rm pol}$ and $1+t^{-1}$ for the remainder $\big(\sum_ce_c(v_c\cdot y)\big)e_0$.
Since $C'\succeq B_P\succ0$, $e_0^{\mathsf T}(C')^{-1}e_0\le b_{00,P}$, and by the Cauchy--Schwarz inequality $\big(\sum_ce_c(v_c\cdot y)\big)^2\le n_P\sum_c\bar e_c^{\,2}(v_c\cdot y)^2$.
Integrate over $P$.
\end{proof}
The bound keeps $v_c\cdot y$ signed, holds for every real $y$ and is a quadratic form in $y$ with polynomial integrands (Lemma~\ref{B-lem:Cprime}). The verifier multiplies the polynomial-branch contribution by $1+t$ ($t$ chosen in floating point) and adds, in interval arithmetic, the positive
semidefinite matrix $(1+t^{-1})\,b_{00,P}\,n_P\sum_c\bar e_c^{\,2}\int_Pv_cv_c^{\mathsf T}$. In every lower-bound verifier of the records this matrix enters
$\lambda_{\rm up}$; a scalar trace bound is only logged. Moduli $|r_j|\le E_j$ of single modes would not suffice: $(\sum_jE_jy_j)^2$ is no bound for $(\sum_jy_jr_j)^2$ ($0$ against $1$ for
$E=(1,1)$, $r=(1,0)$, $y=(1,-1)$), and one would need $(\sum_jE_j)\sum_jE_jy_j^2$.

\emph{Scaling and eigenvalue bound.} The integrals are accumulated for unscaled representers $\hat r_\iota$, with $\sqrt{\alpha_\iota}\,r_\iota=d_\iota\hat
r_\iota$; $d_\iota$ is built from $\alpha_\iota$, the orbit size of Proposition~\ref{oth:prop:classes}(c) and the normalisation of $Y_\iota$ and $\mathcal
N_\iota$, and is enclosed in ball arithmetic. So $G=D_d\hat GD_d$, with $\hat G$ the Gram matrix of the $\hat r_\iota$. For the computed $\hat A$ the
accumulation gives $\varepsilon_G\ge|\hat G_{\iota\iota'}-\hat A_{\iota\iota'}|$, and $G\preceq G_{\rm up}$ by Lemma~\ref{B-lem:gramerr}. Lemma~\ref{lem:verif}(c)
bounds $\lambda_{\max}(D_d\hat AD_d)$ in ball arithmetic \cite{Arb}, after congruence with a dyadic approximate eigenvector matrix $V$, and $\varepsilon_G|d|^2$ is
added. The result is $\lambda_{\rm up}\ge\lambda_{\max}(G_{\rm up})$ of Table~\ref{tab:certlow} and of the records (Section~\ref{tab:certs}).

\subsection{Face integrals of vertex-singular pairs}\label{B-sec:faceint}
The upper bounds need face integrals of the pairs of Definition~\ref{B-def:family} on $Q$. The integrands are homogeneous of degree $2\nu$ or $\nu$, larger than
$-3$; by \eqref{B-eq:rQ} every integral reduces to the far faces $\Sigma_i$, where they are real-analytic up to the boundary. The frames are exactly orthonormal
rational matrices: their columns $(\mathbf p,\mathbf q,\mathbf n)$ are those of $R(\omega)/|\omega|^2$, with $R(\omega)$ the quadratic rotation matrix of a
quaternion $\omega\ne0$ with dyadic entries; the identities of Lemma~\ref{B-lem:vpairs}(a) are checked in rational arithmetic. No frame vector equals the
irrational diagonal $(1,1,1)/\sqrt3$ (for the first frame of V9, $\mathbf p$ is within $8.5\cdot10^{-10}$\,rad of it). The frames, the rounded coefficients, the
panels and the quadrature error bounds are in \recR{3} (V6) and \recR{7} (V9).
The integrals are enclosed by tensor Gauss--Legendre rules with $N_{\rm GL}$ nodes per direction on dyadic squares (panels), with nodes and weights enclosed in
Arb \cite{Arb,FLINT}. On the Bernstein ellipse of a panel side, $\Re(\mathbf a_j\cdot x)\ge b_j>0$ with $b_j$ computed in Arb, hence $|\zeta_j|\le b_j^{-3/2}$ for
$\nu=-\frac32$, also for the conjugated factor; Lemma~\ref{B-lem:gl} bounds the error, which is added to the radius (panel tolerance $10^{-22}$). The rounded
coefficients (Gaussian rationals for V6, dyadic for V9) are exact data. The full procedure is in \docD{9}.
\begin{lemma}[Gauss--Legendre quadrature for analytic integrands]\label{B-lem:gl}
(a) Let $f$ be holomorphic and bounded by $M$ in the open region bounded by the Bernstein ellipse $\mathcal E_\varrho$ with foci $\pm1$, $\varrho>1$, and let $\mathcal Q_{N_{\rm GL}}$ be the $N_{\rm GL}$-point Gauss--Legendre rule on $[-1,1]$. Then
\[\Big|\int_{-1}^1f-\mathcal Q_{N_{\rm GL}}(f)\Big|\le\frac{8M\varrho^{1-2N_{\rm GL}}}{\varrho-1}.\]
(b) On a rectangle $[s_0,s_1]\times[t_0,t_1]$, let $f$ be continuous and, in each variable separately (the other one fixed and real in its interval), satisfy the hypothesis of (a), transplanted to the side in question, with bounds $M_s$ and $M_t$. Then the tensor rule has error at most $4(s_1-s_0)(t_1-t_0)(M_s+M_t)\varrho^{1-2N_{\rm GL}}/(\varrho-1)$.
\end{lemma}
The proof (Appendix~\ref{B-app:gl}) combines $|c_m|\le2M\varrho^{-m}$ for the Chebyshev coefficients \cite[Thm.~8.1]{TrefethenATAP}, the exactness of
$\mathcal Q_{N_{\rm GL}}$ to degree $2N_{\rm GL}-1$ and its positive weights; \cite[Thm.~19.3]{TrefethenATAP} has a bound of this type with another constant.
\begin{remark}[printed versus internal enclosures]\label{B-rem:printed}
The radii printed in Lemma~\ref{lem:esscomp} and in the records of Section~\ref{sec:upper} are those of $25$-digit printouts. $\mathcal R_0(c)$ is evaluated from the internal $160$-bit balls (radii at most $3.2\cdot10^{-27}$ and $4.2\cdot10^{-28}$); the printed balls give the same twelve-digit $L_*$ (\docD{9}).
\end{remark}

\subsection{Records}\label{tab:certs}\label{tab:repo}
The computer-assisted statements refer to records R1--R12 and documents D1--D10 of release v1.1 of \cite{lbb3drepo}. Their files, with SHA-256 hashes, are listed on the page
\url{https://github.com/xfliu/lbb3d-verified/blob/v1.1/docs/records.md}. A record holds the exact input data, the verifier, the command line, the output and
logs, the software versions and the SHA-256 hashes of all inputs, corrector files included. A certificate is checked by re-running its verifier on the archived
data and correctors, without the floating-point design stages.

\subsection{Formal checks and trust base}\label{sec:formal}
Lean~4 \cite{LeanFour} with Mathlib \cite{Mathlib} (\recR{12}) checks: the step from the rational $\kappa_0$ to the decimal bound for $\beta$ in every certified
statement except Theorem~\ref{oth:thm:fich}, whose step is checked in exact rational arithmetic (Corollary~\ref{oth:cor:fich}); admissibility and the exact
quotient of the polynomial trial pairs of degrees $8$, $16$, $24$, $48$ on the cube; and the reliability theorem and the final matrix inequality of the planar
certificates of \cite{companion}. The extent of these checks, and how each numbered result is established, are listed in \docD{10}. Not formalised are the
three-dimensional reliability theorem and Theorem~\ref{thm:fv}. Not machine-checked, besides the results proved in the paper only, are the implemented formulas
and the software stack: CPython integers and rationals, IEEE arithmetic with the error bounds of the interval accumulation, FLINT and Arb.

\section{Conclusions}
Table~\ref{oth:tab:summary} collects the enclosures of $\beta$. \emph{Proved:} the two-sided bounds for the four domains; essential spectrum at or below
$0.116398$ from every convex rectangular vertex; near such a vertex no multiplier certificate certifies $\kappa_0<L_*=0.767204562046$
(Corollary~\ref{B-cor:obstr}). \emph{Computed:} the exactly verified multiplier bound and the enclosed angular trial value agree to $2.1\cdot10^{-6}$
(Theorem~\ref{thm:vexact}, Proposition~\ref{prop:UL}). \emph{Open:} (i) whether $U_*=L$, which would make the method sharp at the vertex; (ii) whether $\beta^2$
of the cube equals $\inf\sigma_{\rm ess}$ (an isolated eigenvalue below the essential spectrum is not excluded); (iii) a convergence proof in three dimensions,
where the planar argument of \cite{companion} needs a substitute for the corner fan at the vertices; (iv) finer interior designs for the other three domains,
whose certificates are limited by the interior design, not by edges or vertices (Remark~\ref{oth:rem:limit}).

\appendix
\section{Proofs of auxiliary results}\label{A:app:lit}
This appendix proves Lemmas~\ref{lem:weak}, \ref{lem:verif}, \ref{B-lem:gl}, \ref{oth:lem:lip} and~\ref{oth:lem:genclass}. Proofs of the
literature results used in Sections~\ref{sec:fv} and~\ref{sec:cert} are reproduced in \docD{2}.

\subsection{Proof of Lemma~\ref{lem:weak}}\label{B-app:weak}
Lemma~\ref{lem:weak} is a Persson-type criterion \cite[Prop.~2.1]{HKW2013}; we give the proof for the bounded self-adjoint operator $\mathcal S$ on $L^2_0(\Om)$.
\begin{proof}[Proof of Lemma~\ref{lem:weak}]
Let $\lambda<\lambda_{\rm e}$ (any real $\lambda$ if $\sigma_{\rm ess}(\mathcal S)=\emptyset$). The set $\sigma(\mathcal S)\cap(-\infty,\lambda]$ is compact and,
by the definition of $\lambda_{\rm e}$, consists of isolated points of $\sigma(\mathcal S)$, each an eigenvalue of finite multiplicity. Hence it
is finite, and the spectral projection $P:=\mathbf 1_{(-\infty,\lambda]}(\mathcal S)$ has finite rank. $P$ commutes with $\mathcal S$, $\mathcal S\ge0$ on $\operatorname{ran}P$
and $\mathcal S\ge\lambda$ on $\operatorname{ran}(I-P)$, so for $\|q\|=1$
\[(\mathcal Sq,q)=(\mathcal SPq,Pq)+(\mathcal S(I-P)q,(I-P)q)\ge\lambda\|(I-P)q\|^2\ge\lambda-|\lambda|\,\|Pq\|^2 .\]
A finite-rank operator maps weakly null sequences to norm null sequences, so $\|Pq_k\|\to0$ and $\liminf_k(\mathcal Sq_k,q_k)\ge\lambda$. If
$\sigma_{\rm ess}(\mathcal S)$ were empty, this would hold for every $\lambda$, contradicting $(\mathcal Sq_k,q_k)\le1$ (Lemma~\ref{A:lem:S}). Letting $\lambda\uparrow\lambda_{\rm e}$
gives the claim.
\end{proof}

\subsection{Proof of Lemma~\ref{lem:verif}}\label{B-app:verif}
\begin{proof}[Proof of Lemma~\ref{lem:verif}]
(a) The identity is the Sherman--Morrison formula (see \cite[Sect.~2.1]{GolubVanLoan}). Moreover $1+\phi b_{00}>0$, because $b_{00}>0$. The function $\phi\mapsto\phi/(1+\phi b_{00})$ is concave with value $0$ at $0$, so it is at least $\lambda_{\rm SM}\phi$ on $[0,\phi_{\max}]$. Subtracting a
larger multiple of $uu^{\mathsf T}\succeq0$ gives a smaller matrix.

(b) $(1+t)u^{\mathsf T}Mu+(1+1/t)v^{\mathsf T}Mv-(u+v)^{\mathsf T}M(u+v)=(\sqrt t\,u-v/\sqrt t)^{\mathsf T}M(\sqrt t\,u-v/\sqrt t)\ge0$.

(c) Since $V$ is invertible, $Ax=\lambda x$ with $x=Vz$ is equivalent to $Yz=\lambda\Xi z$, so the eigenvalues of $A$ are those of the pencil $(Y,\Xi)$. Let $\lambda:=\lambda_{\max}(A)\ge0$, let $z$ be a corresponding vector of the pencil, and let $|z_i|=\max_j|z_j|>0$. Row $i$ gives
$\lambda\Xi_{ii}z_i+\lambda\sum_{j\ne i}\Xi_{ij}z_j=Y_{ii}z_i+\sum_{j\ne i}Y_{ij}z_j$. Hence, using $\lambda\ge0$ and $Y_{ii}\ge0$ ($Y\succeq0$),
$\lambda\Xi_{ii}|z_i|\le\lambda\sum_{j\ne i}|\Xi_{ij}||z_i|+\big(Y_{ii}+\sum_{j\ne i}|Y_{ij}|\big)|z_i|$, i.e.\ $\lambda d_i\le Y_{ii}+\sum_{j\ne i}|Y_{ij}|$. This is the Gershgorin argument \cite[Thm.~6.1.1]{HornJohnson} applied to the pencil.
\end{proof}

\subsection{Gauss--Legendre quadrature for analytic integrands}\label{B-app:gl}
\begin{proof}[Proof of Lemma~\ref{B-lem:gl}]
(a) By \cite[Thm.~8.1]{TrefethenATAP}, $f=\sum_{m\ge0}c_mT_m$ uniformly on $[-1,1]$ with $|c_m|\le2M\varrho^{-m}$. The rule $\mathcal Q_{N_{\rm GL}}$ is exact
for degree at most $2N_{\rm GL}-1$ and has positive weights with sum $2$, so $|\int_{-1}^1T_m|\le2$ and $|\mathcal Q_{N_{\rm GL}}(T_m)|\le2$; the error is at most
$\sum_{m\ge2N_{\rm GL}}4\cdot2M\varrho^{-m}=8M\varrho^{1-2N_{\rm GL}}/(\varrho-1)$. This is the argument of \cite[Thm.~19.3]{TrefethenATAP} with a cruder constant.
On $[s_0,s_1]$ the affine change of variables multiplies the error by $(s_1-s_0)/2$.
(b) With $I_s,\mathcal Q_s$ and $I_t,\mathcal Q_t$ the integrals and rules in the two variables,
$I_sI_tf-\mathcal Q_s\mathcal Q_tf=I_t[(I_s-\mathcal Q_s)f]+\mathcal Q_s[(I_t-\mathcal Q_t)f]$. Apply (a) for each real $t$ and for each node $s$; the weights of
$\mathcal Q_s$ are positive with sum $s_1-s_0$.
\end{proof}
A detailed version of this proof is in \docD{2}.

\subsection{Proof of Lemma~\ref{oth:lem:lip}}\label{A:app:lip}
\begin{proof}[Proof of Lemma~\ref{oth:lem:lip}]
Write $x_t:=x+t(y-x)$. \emph{Step 1: $[x,y]\subset\overline\Om$.} For $\tau\in\mathfrak T$ the set $I_\tau:=\{t\in[0,1]:x_t\in\tau\}$ is a closed interval (possibly empty or a
point), since $\tau$ is convex and closed, and $\bigcup_\tau I_\tau=[0,1]$. The end points of all $I_\tau$ divide $[0,1]$ into finitely many closed subintervals; each of them is
contained in some $I_\tau$ (an interior point of the subinterval lies in some $I_\tau$, whose end points are among the division points). On such a subinterval
$t\mapsto m(x_t)$ is affine with derivative $Dm|_\tau(y-x)$, of norm at most $K_m|x-y|$. Since $t\mapsto m(x_t)$ is continuous, summing over the subintervals gives
$|m(x)-m(y)|\le K_m|x-y|$.

\emph{Step 2: general case.} Assume $[x,y]\not\subset\overline\Om$ and put $t_1:=\sup\{t\in[0,1]:x_s\in\overline\Om\ \text{for all }s\in[0,t]\}$ and
$t_2:=\inf\{t\in[0,1]:x_s\in\overline\Om\ \text{for all }s\in[t,1]\}$. Since $\overline\Om$ is closed, $p:=x_{t_1}$ and $q:=x_{t_2}$ lie in $\overline\Om$, and
$[x,p]$, $[q,y]\subset\overline\Om$. If $t_1>t_2$, these two segments would cover $[x,y]$; hence $t_1\le t_2$. If $p$ were in the open set
$\Om$, then $x_s\in\Om$ for all $s$ slightly larger than $t_1$, contradicting the definition of $t_1$ (note $t_1<1$); so $p\in\partial\Om$, and likewise $q\in\partial\Om$. Hence
$m(p)=m(q)=0$, and by Step~1
\[|m(x)-m(y)|\le|m(x)-m(p)|+|m(q)-m(y)|\le K_m\big(t_1+(1-t_2)\big)|x-y|\le K_m|x-y| .\qedhere\]
\end{proof}

\subsection{Proof of Lemma~\ref{oth:lem:genclass}}\label{A:app:genclass}
\begin{proof}[Proof of Lemma~\ref{oth:lem:genclass}]
(a) This is the first part of the proof of Proposition~\ref{oth:prop:classes}(a), which does not use $R^2=I$.

(b) The proof of Proposition~\ref{oth:prop:classes}(b),(c) is written with $\Lambda_R^{\mathsf T}$ and $R^{\mathsf T}$. It uses $R^2=I$ only through the properties of
$P_\chi$, which hold by (a), and through (iii), which is replaced here by $C'(Rx)=\Lambda_RC'(x)\Lambda_R^{\mathsf T}$. Hence it applies verbatim.

(c) Every point of $\overline\Om$ is of the form $Rx$ with $x\in\Om_{\mathcal G}$; if $Rx=R'x'$ with $x,x'\in\Om_{\mathcal G}$, then $x'=R''x$ with
$R'':=R'^{-1}R$, and the hypothesis applied to $x$ and $R''$ together with the homomorphism property of $R\mapsto\chi(R)\det R\,R$ and of $\chi$ shows that the two
definitions in \eqref{oth:eq:corrext} agree. So the extension is well defined; it is continuous, piecewise polynomial on the images of the tetrahedra, and zero on
$\partial\Om=\bigcup_RR(\partial\Om\cap\Om_{\mathcal G})$; hence $\xi\in H^1_0(\Om)^3\times H^1_0(\Om)$ (a continuous function that is polynomial on each tetrahedron of a
finite partition of the polyhedron $\overline\Om$ and vanishes on $\partial\Om$ lies in $H^1_0(\Om)$). By construction $\tilde\xi:=(\det R\,R^{\mathsf T}\varphi(R\,\cdot),
\vartheta(R\,\cdot))=\chi(R)\xi$ for every $R\in\mathcal G$. Finally $\mathcal D^*\tilde\xi=\mathcal T_R\mathcal D^*\xi$: by the chain rule
$\Div\tilde\varphi(x)=\det R\,(\Div\varphi)(Rx)$, by the identity $\operatorname{ax}(R^{\mathsf T}AR)=\det R\,R^{\mathsf T}\operatorname{ax}(A)$ of the proof of
Lemma~\ref{A:lem:sym} $\curl\tilde\varphi(x)=R^{\mathsf T}(\curl\varphi)(Rx)$, and $\nabla\tilde\vartheta(x)=R^{\mathsf T}(\nabla\vartheta)(Rx)$; compare \eqref{A:TR}.
Hence $\mathcal T_R\mathcal D^*\xi=\chi(R)\mathcal D^*\xi$, which is \eqref{A:parity}.
\end{proof}

\section{Design and verifier data}\label{app:design}
\setlength{\emergencystretch}{1.5em}%
This appendix lists the data that identify the certified computations; it contains no proof. The descriptions of the designs and their floating-point limits are in
\docD{7}, and the exact upper bounds for lower polynomial degrees in \docD{6}. The validation runs of the verifiers are not part of the proofs; they are in \docD{8}.

\emph{The cube.} Table~\ref{tab:lower} lists, for the best design per $n$, the value $\kappa_{\rm PD}$ and the smallest $\kappa_0$ for which all four symmetry classes pass the criterion of Theorem~\ref{thm:rel} in floating point (Gram integrals by
Gaussian quadrature). The $16$-sector design with $12$ balls, used for the $n=16$ certificates, passes at $\kappa_0=0.76851$; Section~\ref{sec:lower} explains why
Table~\ref{tab:certlow} nevertheless certifies $\kappa_0=0.7681$ with it.
The verifier refuses $\kappa_0=0.768$ for $n=16$ and $\kappa_0=0.7676$ for $n=32$. These values lie just below the exact floors
$\kappa_{\rm PD}^{\rm ex}=0.768012$ and $0.767645$ \eqref{oth:eq:kpd} of these designs, which lie above the enclosed lower bound $0.767204561906$ of Theorem~\ref{thm:upper}. The
exact floors lie slightly above the floating-point values of Table~\ref{tab:lower}, which are evaluated with $C_{\rm low}$, in agreement with
Lemma~\ref{oth:lem:floors}(a). Tests with corrupted correctors, a comparison with an independent rule, Monte Carlo estimates and the computing times are listed in
\docD{8}.
\begin{table}[tbp]
\caption{Cube, floating point: limit of positive definiteness $\kappa_{\rm PD}$ \eqref{B-eq:kPD} and smallest $\kappa_0$ passing the criterion in all four classes, for the best design per $n$ ($n=8$: the best of the $8$-sector designs tested in the floating-point design study, not used in any certificate; $n=16$: $25$ balls;
$n=32$: $22$ balls, the design of the certificate). Floating-point evaluation with $Q_1$ correctors and the cellwise weight; see Section~\ref{sec:lower}.}\label{tab:lower}
\begin{tabular}{rlll}
\toprule
sectors $n$ & $\kappa_{\rm PD}$ & $\kappa_0$ & $\beta\ge$ (if certified)\\
\midrule
8 & $0.7705$ & $0.7745$ & $0.3357$\\
16 & $0.76795$ & $0.76845$ & $0.3402$\\
32 & $0.76764$ & $0.76784$ & $0.3407$\\
\bottomrule
\end{tabular}
\end{table}

\emph{The box, the L-shaped prism and the Fichera domain.} Table~\ref{oth:tab:design} lists the data checked by the verifier in the proofs of
Theorems~\ref{oth:thm:box}, \ref{oth:thm:lp} and~\ref{oth:thm:fich}. For all three domains \docD{8} reports refusals, corrupted correctors and Monte Carlo estimates.
In addition it reports, for the box, twelve mutations of the data, the refusals that bracket $\kappa_{\rm PD}^{\rm ex}$, a Monte Carlo estimate by an independent
implementation and a quadrature cross-check; for the L-shaped prism, a test of the diagonal reduction, an independent evaluation of the upper bound and two weaker
certificates at $\kappa_0=\frac{47}{50}$ and $\frac{93}{100}$; for the Fichera domain, checks of the symmetry reduction, a regression against the verifier of the L-shaped
prism and three weaker certificates.
\begin{table}[t]
\caption{Data of the certificates of Theorems~\ref{oth:thm:box}, \ref{oth:thm:lp} and~\ref{oth:thm:fich}. Mesh: Kuhn subdivision of a tensor-product grid graded
geometrically towards the listed planes (first layer, ratio, largest step); numbers of tetrahedra, vertices and faces refer to the fundamental region; ``pieces'' are the
pieces on which the exact weight $B_P$ is evaluated (distinct matrices $B_P$ in parentheses); ``rules'' are the orders $s$ of the Grundmann--M\"oller rules and the
number of points per class; $p$: Taylor degree; $\kappa_0$: certified value.}\label{oth:tab:design}
\footnotesize
\begin{tabular}{@{}p{0.13\textwidth}p{0.25\textwidth}p{0.25\textwidth}p{0.25\textwidth}@{}}
\toprule
 & box $\Om_{\rm B}$ & L-prism $\Om_{\rm L}$ & Fichera $\Om_{\rm F}$\\
\midrule
region & octant $[0,1]\times[0,\frac12]^2$, $|\mathcal G|=8$ & $\{x_2\ge x_1\}\cap(\overline{\omega_{\rm L}}\times[0,\frac12])$, $|\mathcal G|=4$ & ${\rm Q}0$: $\Om_{\mathcal H_{\rm F}}$, $|\mathcal H_{\rm F}|=2$; ${\rm S}1$: $\Om_{\mathcal G_{\rm F}}$, $|\mathcal G_{\rm F}|=6$\\
mesh & outer faces; $\frac1{512}$, 2, $\le\frac1{32}$ & all boundary planes; $\frac1{256}$, 2, $\le\frac1{16}$ & planes $x_i\in\{-1,0,1\}$; $\FichRho$, 2, $\le\FichHmax$\\
tetrahedra (vertices) & $86\,400$ ($16\,317$) & $62\,208$ ($12\,025$) & ${\rm S}1$: $\FichTetsSix$ ($\FichVertsSix$); ${\rm Q}0$: $\FichTetsTwo$ ($\FichVertsTwo$)\\
faces (boundary) & $176\,480$ ($7\,360$) & $127\,584$ ($6\,336$) & ${\rm S}1$: $\FichFacesSix$ ($\FichBdSix$); ${\rm Q}0$: $\FichFacesTwo$ ($\FichBdTwo$)\\
balls & $16$ ($11$ orbits of $\Pi_{23}$) & $21$ orbit representatives ($40$ in $\omega_{\rm L}\times[0,\frac12]$, $58$ in $\Om$) & $\FichOrbits$ $\mathcal G_{\rm F}$-orbits ($\FichBallsTwo$ in $\Om_{\mathcal H_{\rm F}}$, $\FichBalls$ in $\Om$)\\
degree; radii & $p=4$; $R\in\{\frac12,\frac14\}$ & $p=3$; $R\in\{\frac12,\frac14,\frac{181}{512}\}$ & $p=\FichTaylor$; $R\in\FichRadii$\\
$N$ per class & $244$, $246$, $265$, $265$, $286$, $289$ & $208$, $212$, $252$, $256$ & ${\rm Q}0$: $\FichNQ$, ${\rm S}1$: $\FichNS$\\
correctors & P2 on the mesh, $490\,853$ dofs & P2 on the red refinement ($497\,664$ tetrahedra), $2\,756\,741$ dofs & P2 on the mesh, $\FichDofQ$ (${\rm Q}0$), $\FichDofS$ (${\rm S}1$) dofs\\
pieces (distinct) & $1\,782\,878$ ($697\,772$) & $5\,710\,774$ ($1\,404\,904$) & $\FichPiecesQ$ ($\FichDistinctQ$); $\FichPiecesS$ ($\FichDistinctS$)\\
rules; points & $s\le9$; $1.18\cdot10^9$ & $s\le8$; $2.52\cdot10^9$ & $s\le\FichRuleMax$; $\FichPointsQ$, $\FichPointsS$\\
$\kappa_0$ & $\frac{87}{100}$ & $\frac{37}{40}$ & $\FichKappaZero$\\
\bottomrule
\end{tabular}
\end{table}
\setlength{\emergencystretch}{0pt}
\section*{Code and data availability}
The verifiers, the exact input data and the records of all computer-assisted statements (Section~\ref{tab:certs}) are in release v1.1 of \cite{lbb3drepo}, \url{https://github.com/xfliu/lbb3d-verified}. The corrector files and Gram matrices ($0.1$--$3.3$\,GB each) are identified there by their SHA-256 hashes and are available from the author; they will be deposited with the release on Zenodo (DOI to be added on acceptance).

\end{document}